\documentclass{amsart}

\usepackage[utf8]{inputenc}
\usepackage[english]{babel}
\usepackage{graphicx}
\usepackage{amssymb}
\usepackage{amsmath}
\usepackage{tikz-cd}
\usepackage{adjustbox}
\usepackage{mathabx}
\usepackage{hyperref}
\usepackage[backend=bibtex,
style=alphabetic,
bibencoding=ascii
]{biblatex}
\usepackage{rotating}
\newcommand*{\isoarrow}[1]{\arrow[#1,"\rotatebox{90}{\(\sim\)}"
]}

\newtheorem{theorem}{Theorem}[subsection]
\newtheorem{lemma}[theorem]{Lemma}
\newtheorem{proposition}[theorem]{Proposition}
\newtheorem{corollary}[theorem]{Corollary}
\newtheorem{conjecture}[theorem]{Conjecture}
\theoremstyle{definition}
\newtheorem{definition}[theorem]{Definition}
\newtheorem{example}[theorem]{Example}

\newtheorem{question}[theorem]{Question}

\theoremstyle{remark}
\newtheorem{remark}[theorem]{Remark}

\title[The Cohomology of Rapoport-Zink Spaces of EL-Type]{The Cohomology of Unramified Rapoport-Zink Spaces of EL-Type and Harris's Conjecture}
\author{Alexander Bertoloni Meli} 
\address{
  \small{$^1$University of Michigan, 530 Church St, Ann Arbor, MI 48109, United States\\ email: \href{mailto:abertolo@umich.edu}{abertolo@umich.edu}}
  }
\begin{document}

\maketitle

\begin{abstract}
We study the $l$-adic cohomology of unramified Rapoport-Zink spaces of EL-type. These spaces were used in Harris and Taylor's proof of the local Langlands correspondence for $\mathrm{GL_n}$ and to show local-global compatibilities of the Langlands correspondence. In this paper we consider certain morphisms, $\mathrm{Mant}_{b, \mu}$, of Grothendieck groups of representations constructed from the cohomology of the above spaces, as studied by Harris and Taylor, Mantovan, Fargues, Shin, and others. Due to earlier work of Fargues and Shin we have a description of $\mathrm{Mant}_{b, \mu}(\rho)$ for $\rho$ a supercuspidal representation. In this paper, we give a conjectural formula for $\mathrm{Mant}_{b, \mu}(\rho)$ for $\rho$ an admissible representation and prove it when $\rho$ is essentially square integrable. Our proof works for general $\rho$ conditionally on a conjecture appearing in Shin's work. We show that our description agrees with a conjecture of Harris in the case of parabolic inductions of supercuspidal representations of a Levi subgroup.
\end{abstract}

\section{Introduction}
Our goal in this paper is to give a description of the $l$-adic cohomology of unramified Rapoport-Zink spaces of EL-type. These spaces are moduli spaces of $p$-divisible groups associated to unramified Weil-restrictions of general linear groups and can be thought of as generalizations of Lubin-Tate spaces.

This work generalizes, for these particular spaces, the Kottwitz conjecture stated in \cite[Conj 7.3]{RV1}. The Kottwitz conjecture describes the supercuspidal part of the $l$-adic cohomology of Rapoport-Zink spaces, and is known in the cases we consider by work of Shin \cite[Cor 1.3]{Shi1}. We prove our description of this cohomology is compatible with a conjecture of Harris \cite[Conj 5.4]{Har1}, generalizing the Kottwitz conjecture to parabolic inductions of supercuspidal representations. 

Our result describes the cohomology of these Rapoport-Zink spaces as a formal alternating sum (indexed by certain root theoretic data) of representation-theoretic constructions including the local Langlands correspondence, parabolic inductions, and Jacquet modules. 

We prove our result inductively using two formulas from the literature. The first of these is Shin's averaging formula \cite[Thm 7.5]{Shi1} which is proven using Mantovan's formula \cite[Thm 22]{Man2}. Mantovan's formula connects the cohomology of Rapoport-Zink spaces, Igusa varieties and Shimura varieties. The second formula is the Harris-Viehmann conjecture of \cite[Conj 8.4]{RV1} which relates the cohomology of so-called non-basic Rapoport-Zink spaces to a product of Rapoport-Zink spaces of lower dimension. A proof of this conjecture is expected to appear in a forthcoming paper of Scholze. 

To carry out our induction, we prove combinatorial analogues of the above formulas phrased purely in terms of root-theoretic data. Interestingly, we are able to prove these analogues for general quasisplit reductive groups, though at present we can only connect them to the cohomology of Rapoport-Zink spaces of unramified EL-type. To do so in other cases, one would need to generalize Shin's averaging formula. 

We now describe our main results more precisely.  We fix an algebraic closure $\overline{\mathbb{Q}_p}$ of $\mathbb{Q}_p$. We study Rapoport-Zink spaces of unramified EL-type which we denote $\mathbb{M}_{b, \mu}$. These are moduli spaces of $p$-divisible groups coming from an unramified EL-datum consisting of 
\begin{enumerate}
    \item a finite unramified extension $F \subset \overline{\mathbb{Q}_p}$ of $\mathbb{Q}_p$,
    \item a finite dimensional $F$ vector space $V$ which defines the group \\
    $G=\mathrm{Res}_{F/ \mathbb{Q}_p} \mathrm{GL}(V)$, 
    \item a $G_{\overline{\mathbb{Q}_p}}$-conjugacy class of cocharacters $\{\mu\}$, with  $\mu: \mathbb{G}_m \to G_{\overline{\mathbb{Q}_p}}$, and such that the weights of $\mu$ are elements of $\{0, 1 \}$.
    \item an element $b$ of a finite set $\mathbf{B}(G, \mu)$ which defines a group $J_b$ that is an inner twist of a Levi subgroup $M_b$ of $G$. 
\end{enumerate}
Roughly one can think of $b, \mu$ as specifying the Newton and Hodge polygons of a $p$-divisible group and $J_b$ as the automorphism group of the isocrystal $b$. 

Let $\mathbb{Q}^{ur}_p$ denote the maximal unramified extension of $\mathbb{Q}_p$ inside $\overline{\mathbb{Q}_p}$, and let $\widehat{\mathbb{Q}^{ur}_p}$ denote its completion. Then the spaces $\mathbb{M}_{ b, \mu}$ are formal schemes over $\widehat{\mathbb{Q}^{ur}_p}$. One constructs a tower of rigid spaces $\mathbb{M}^{rig}_{U, b, \mu}$ over the generic fiber $\mathbb{M}^{rig}_{b, \mu}$ of $\mathbb{M}_{b, \mu}$, where the index $U$ runs over compact open subgroups of $G(\mathbb{Q}_p)$. Associated to such a tower we have a cohomology space $[H^{\bullet}(G,b, \mu)]$ which is an element of the Grothendieck group $\mathrm{Groth}(G(\mathbb{Q}_p)\times J_b(\mathbb{Q}_p) \times W_{E_{\{\mu\}_G}})$ of admissible representations of $G(\mathbb{Q}_p), J_b(\mathbb{Q}_p)$ and $W_{E_{\{\mu\}_G}}$, where the latter group is the Weil group of the reflex field, $E_{\{\mu\}_G}$, of $\{\mu\}$. This construction can be thought of as an alternating sum of a direct limit over $U \subset G$ of $l$-adic cohomology groups with the actions of $G(\mathbb{Q}_p)$ and $J_b(\mathbb{Q}_p)$ arising from Hecke correspondences and isogenies of $p$-divisible groups, respectively. We refer to \S\ref{Mant defn} for a precise definition.

The cohomology object $[H^{\bullet}(G,b, \mu)]$ gives rise to a map of Grothendieck groups
\[
\mathrm{Mant}_{G,b,\mu}: \mathrm{Groth}({J_b}(\mathbb{Q}_p)) \to \mathrm{Groth}(G(\mathbb{Q}_p) \times W_{E_{\{\mu\}_G}})
\]
which maps a representation $\rho$ to the alternating sum of the $J_b(\mathbb{Q}_p)$-linear $\mathrm{Ext}$ groups of $[H^{\bullet}(G,b, \mu)]$ and $\rho$.  

The map $\mathrm{Mant}_{G,b,\mu}$ has been studied by many authors. Harris and Taylor \cite{HT1} used this construction to prove the local Langlands correspondence for general linear groups. It also appears naturally in Mantovan's work relating the cohomology of Shimura varieties, Igusa varieties, and Rapoport-Zink spaces \cite{Man2}. Fargues studied $\mathrm{Mant}_{G,b,\mu}$ for basic $b$ in some $EL$ and $PEL$-cases in \cite{Far1}. Shin combined Mantovan's formula with his trace formula description of the cohomology of Igusa varieties to prove instances of local-global Langlands compatibilities \cite{Shi2}.

In \cite{Shi1}, Shin proved an averaging formula for $\mathrm{Mant}_{G, b, \mu}$ which is key to our work. He defined a map 
\[
\mathrm{Red}_b: \mathrm{Groth}(G(\mathbb{Q}_p)) \to \mathrm{Groth}(J_b(\mathbb{Q}_p))
\]
which up to a character twist is given by composing the un-normalized Jacquet module 
\[
\mathrm{Jac}^G_{P^{op}_b}: \mathrm{Groth}(G(\mathbb{Q}_p)) \to \mathrm{Groth}(M_b(\mathbb{Q}_p))
\]
with the Jacquet-Langlands map of Badulescu [Bad1]
\[
\mathrm{LJ}: \mathrm{Groth}(M_b(\mathbb{Q}_p)) \to \mathrm{Groth}(J_b(\mathbb{Q}_p)) .
\]
Shin uses global methods and so necessarily works with a large but inexplicit class of representations which he denotes \emph{accessible}. This set loosely consists of those representations isomorphic to the $p$-component of an automorphic representation appearing in the cohomology of a certain unitary similitude group Shimura variety. In particular, the essentially square integrable representations in $\mathrm{Groth}(G(\mathbb{Q}_p))$ are accessible. 

In what follows  $r_{-\mu}$ is a finite dimensional representation of $\widehat{G} \rtimes W_{E_{\{\mu\}_G}}$ which restricts to the representation of highest weight $-\mu$ on $\widehat{G}$, and $LL$ is the semisimplifed local Langlands correspondence from \cite{HT1}. Shin shows the following result.
\begin{theorem}[Shin's Averaging Formula]
Assume $\pi$ is an accessible representation of $G(\mathbb{Q}_p)$. Then
\[
\sum\limits_{b \in \mathbf{B}(G ,\mu)} \mathrm{Mant}_{G, b, \mu}(\mathrm{Red}_b(\pi))=[\pi][r_{-\mu} \circ \mathrm{LL}(\pi)|_{W_{E_{\{\mu\}_G}}}],
\]
where the above formula is correct up to a Tate twist which we omit for clarity and $[\pi][\rho]$ is our notation for an element $\pi \boxtimes \rho \in \mathrm{Groth}(G(\mathbb{Q}_p) \times W_{E_{\{\mu\}_G}})$.
\end{theorem}

Additionally we have the conjecture of Harris and Viehmann which allows us to write $\mathrm{Mant}_{G, b, \mu}$ for non-basic $b$ ($b$ is basic when it corresponds to an isocrystal with a single slope) in terms of $\mathrm{Mant}_{G', b', \mu'}$ such that $G'$ is a general linear group of smaller rank than $G$. This conjecture was formulated in work of \cite{Har1} and \cite{RV1} and is expected to be proven in forthcoming work of Scholze. In what follows, $\mathrm{Ind}$ is the un-normalized parabolic induction functor.
\begin{conjecture}[Harris-Viehmann]
\[
\mathrm{Mant}_{G,b, \mu}= \sum\limits_{(M_b, \mu') \in \mathcal{I}^{G, \mu}_{M_b, b'}} \mathrm{Ind}^G_{P_b}(\otimes^k_{i=1} \mathrm{Mant}_{M_{b'_i}, b'_i, \mu'_i}),
\]
where we omit a Tate twist which we discuss at length in \S\ref{HV}. The finite set $\mathcal{I}^{G, \mu}_{M_b, b'}$ is described in Proposition \ref{I_M,b,mu defn}.
\end{conjecture}

Shin's averaging formula and the Harris Viehmann conjecture allow one to compute $\mathrm{Mant}_{G, b, \mu} \circ \mathrm{Red}_b$ recursively. The latter lets us compute $\mathrm{Mant}_{G,b,\mu}$ for non-basic $b$ given that we know $\mathrm{Mant}_{G',b', \mu'}$ for $G'$ of lower rank and the former lets us compute $\mathrm{Mant}_{G,b, \mu}$ for the unique basic $b \in \mathbf{B}(G, \mu)$ if we know it for all non-basic $b \in \mathbf{B}(G, \mu)$. One of our main results is to give a non-recursive description of $\mathrm{Mant}_{G,b,\mu} \circ \mathrm{Red}_b$ which we now describe.

Let $G= \mathrm{Res}_{F/\mathbb{Q}_p} \mathrm{GL}(V)$ as before, choose a rational Borel subgroup $B$ of $G$, and a rational maximal torus $T \subset B \subset G$. Then we consider pairs $(M_S, \mu_S)$ where $M_S \subset T$ is a Levi subgroup of a parabolic subgroup $P_S$ containing $B$, and $\mu_S \in X_*(T)$ is dominant as a cocharacter of $M_S$. We call a pair of the above form a \emph{cocharacter pair} for $G$.

We associate to a cocharacter pair $(M_S, \mu_S)$ the map of representations $[M_S, \mu_S]: \mathrm{Groth}(G(\mathbb{Q}_p)) \to \mathrm{Groth}(G(\mathbb{Q}_p) \times W_{E_{\{ \mu_S \}_{M_S}}})$, which up to a character twist is given by
\[
\pi \mapsto [(\mathrm{Ind}^G_{P_S} \circ [\mu_S] \circ \mathrm{Jac}^G_{P^{op}_S})(\pi)]
\]
and 
\[
[\mu_S]: \mathrm{Groth}(M_S(\mathbb{Q}_p)) \to \mathrm{Groth}(M_S(\mathbb{Q}_p) \times W_{E_{\{\mu_S\}_{M_S}}})
\]
given by 
\[
\pi \mapsto [\pi][r_{-\mu_S} \circ LL(\pi)]
\]
Then our main result, which follows from Theorem \ref{3.8} in this paper is
\begin{theorem} \label{1.3}
 Suppose $\mathrm{Mant_{G,b, \mu}}$ corresponds to a tower of unramified Rapoport-Zink spaces of EL-type. We assume that the Harris-Viehmann conjecture is true. Then if $\rho \in \mathrm{Groth}(G(\mathbb{Q}_p))$ is essentially square-integrable, we have 
 \[
\mathrm{Mant}_{G, b, \mu}( \mathrm{Red}_b(\rho))= \sum\limits_{(M_S, \mu_S) \in \mathcal{R}_{G, b, \mu}} (-1)^{L_{M_S, M_b}}[M_S, \mu_S](\rho),
 \]
where $\mathcal{R}_{G, b, \mu}$ is a collection of cocharacter pairs with a combinatorial definition and $(-1)^{L_{M_S, M_b}}$ is an easily determined sign. 
\end{theorem}
Shin conjectures (\cite[Conj 8.1]{Shi1}) that the averaging formula holds for all admissible representations of $G(\mathbb{Q}_p)$. If this were indeed the case, then our result would also immediately hold for all admissible representations of $G(\mathbb{Q}_p)$.
 
A crucial part of the proof of the above theorem is the following unconditional result, which is perhaps interesting in its own right.

\begin{theorem}[Imprecise version of Theorem \ref{combinatorial sum} and Corollary \ref{combinatorial induction} of our paper]
For general quasisplit $G$ and a cocharacter $\mu$ (not necessarily minuscule), combinatorial analogues of Shin's formula and the Harris-Viehmann conjecture hold true.
\end{theorem}
This result suggests that perhaps the combinatorics of cocharacter pairs is related to  $\mathrm{Mant}_{G, b, \mu}$ in cases more general than Rapoport-Zink spaces of unramified EL-type. However, we caution the reader that the existence of nontrivial $L$-packets and nontrivial endoscopy in more general groups will likely complicate the situation.

In \S\ref{section 4} of the paper, we use our combinatorial formula to prove the EL-type cases of a conjecture of Harris (\cite[Conj 5.4]{Har1}). This conjecture describes $\mathrm{Mant}_{G, b, \mu}(I^G_M(\rho))$ for $\rho$ a supercuspidal representation of $M(\mathbb{Q}_p)$ for $M$ a Levi subgroup of $G$. In this case, $I^G_M$ denotes normalized parabolic induction. In particular, we show the following result, which is stated as Conjecture \ref{Harris Conjecture} in our paper.
\begin{theorem}[Harris conjecture]
 We assume that Shin's averaging formula holds for all admissible representations of $G(\mathbb{Q}_p)$ and that the Harris-Viehmann conjecture is true. Let $\rho$ be a supercuspidal representation of $M(\mathbb{Q}_p)$. Then up to a precise character twist and sign which we omit for clarity,
 
 \[
\mathrm{Mant}_{G,b,\mu}(LJ (I^{M_b}_M(\rho)))=[I^G_M(\rho)] \left[\bigoplus_{(M, \mu') \in \mathrm{Rel}^{G, \mu}_{M, b}} r_{-\mu'} \circ LL(\rho) \right]
\]
for an explicit set of cocharacter pairs $\mathrm{Rel}^{G, \mu}_{M,b}$.
\end{theorem}
We prove our result for $I^G_M(\rho)$ not necessarily irreducible and $b$ not necessarily basic, which is a generalization of what Harris conjectured for the $G$ we consider.

Finally, in Appendix \ref{examples} we give an example to show that for general representations $\rho$, one cannot hope for an expression as simple as that in Harris's conjecture.

\subsection*{Acknowledgements}
I would like to thank Sug Woo Shin for suggesting I study the cohomology of Rapoport-Zink spaces and for countless helpful discussions on this topic. I thank Michael Harris for a fruitful conversation and for suggesting that my work might allow the verification of some cases of his conjecture. I thank Peter Scholze for a discussion in which he explained that the extra Tate twist in the Harris-Viehmann conjecture arises naturally in his work. This work is partially supported by NSF grant DMS-1646385 (RTG grant).

\section{Cocharacter Formalism}\label{2}
In this section we define and study the notion of a \emph{cocharacter pair}. This notation will be used in the third and fourth sections of this paper, where we describe the cohomology of certain Rapoport-Zink spaces in terms of cocharacter pairs. We endeavor to use a similar notation to \cite{Kot1}.

This section is divided into five subsections. These are structured so that the first contains the basic definitions and the fourth and fifth subsections contain the most important results. The second and third subsections prove a number of technical lemmas that the reader may want to skip at first and refer to as necessary.

\subsection{Notation and Preliminary Definitions}
For the remainder of this section, we fix $G$ a connected quasisplit reductive group defined over $\mathbb{Q}_p$. This is a significantly more general setting than we will need for applications in this paper. However, we choose to work in this generality because doing so is both conceptually clearer and potentially useful for future applications. The ideas in $\S 5$ of \cite{Kot1} might allow one to remove the quasisplit assumption, but we do not attempt this here as it is unnecessary for the applications. Moreover, Kottwitz's study of the set $\mathbf{B}(G)$ in that section relies on understanding the quasisplit case first. 

\begin{remark}
 The reader will notice that most of this section makes sense over an arbitrary field. The assumption that we work over $\mathbb{Q}_p$ is used in section \ref{CWI} when we connect cocharacter pairs to the set $\mathbf{B}(G)$ defined by Kottwitz. However, in $\S 5.1$ of \cite{Kot1}, Kottwitz shows that over $\mathbb{Q}_p$, the set $\mathbf{B}(G)$ is parametrized by a disjoint union of sets of the form $X^*(Z(\widehat{M_S})^{\Gamma})^+$ for $M_S$ a standard Levi subgroup of $G$. These latter sets make sense over general fields and one could make sense generally of all the results of this section by replacing $\mathbf{B}(G)$ with the sets parametrizing it.
\end{remark}

Since $G$ is quasisplit, we can pick a Borel subgroup $B \subset G$ defined over $\mathbb{Q}_p$ and a maximal split torus $A \subset B$ of $G$. We choose $T$ to be a maximal torus defined over $\mathbb{Q}_p$ satisfying $A \subset T \subset B$. We define $X^*(A)$ and  $X_*(A)$ respectively to be the character and cocharacter groups of $A_{\overline{\mathbb{Q}_p}}$. 

The group $G$ has a relative root datum $(X^*(A), \Phi^*(G,A), X_*(A), \Phi_*(G,A))$, where $\Phi^*(G,A)$ and $\Phi_*(G,A)$ respectively denote the set of relative roots and relative coroots of $G$ and the torus $A$. Our choice of Borel subgroup $B$ determines a decomposition $ \Phi^*(G,A)=\Phi^*(G,A)^+ \coprod \Phi^*(G,A)^-$ of positive and negative roots and a subset $\Delta \subset \Phi^*(G,A)^+$ of simple roots. Analogous statements are also true for the coroots. The set of parabolic subgroups $P \supset B$ defined over $\mathbb{Q}_p$ are called \emph{standard parabolic subgroups}. We define $P_S$ to be the unique standard parabolic subgroup such that $\Phi^*(P_S, A)= \Phi^*(G,A)^+ \cup ( \Phi_*(G,A)^- \cap \mathrm{Span}_{\mathbb{Z}}(S))$. There is an inclusion \emph{preserving} bijection between the set of standard parabolic subgroups and subsets of $\Delta$ given by $S \mapsto P_S$.

We let $N_S$ be the unipotent radical of the standard parabolic subgroup $P_S$. It is a standard result that there exists a connected reductive subgroup $M \subset P_S$ so that the natural map $M \to P_S/N_S$ is an isomorphism. In particular, this gives us a Levi decomposition $P_S=MN_S$ and the subgroup $M$ is called a Levi subgroup of $P_S$. The subgroup $M$ is not unique but any two Levi subgroups of $P_S$ are conjugate by an element of $N_S$. However, we have fixed a maximal torus $T$ and there is a unique Levi subgroup $M_S$ containing $T$. The subgroup $M_S$ is constructed explicitly as the centralizer $C_{G}(Z)$, where $Z \subset T$ is the connected component of the intersection of the kernels of the roots in $S$. We refer to the Levi subgroups $M_S$ that we produce in this way as \emph{standard Levi subgroups}.

Define 
\[
\mathfrak{A} := X_*(A).
\]
We have the closed rational Weyl chamber 
\[
\overline{C}_{\mathbb{Q}}= \{ x \in  \mathfrak{A}_{\mathbb{Q}} : \langle x, \alpha \rangle \geq 0, \alpha \in \Delta \}.
\]
We define  for each standard Levi subgroup, 
\[
\mathfrak{A}_{M_S, \mathbb{Q}} := \{ x \in \mathfrak{A}_{\mathbb{Q}}: \langle x, \alpha \rangle=0, \alpha \in S\},
\]
and denote the \emph{strictly dominant} elements of $\mathfrak{A}_{M_S, \mathbb{Q}}$ by
\[
\mathfrak{A}^+_{M_S, \mathbb{Q}}= \{ x \in \mathfrak{A}_{ \mathbb{Q}} : \langle x, \alpha \rangle  = 0, \alpha \in S, \langle x, \alpha \rangle > 0, \alpha \in \Delta \setminus S \},
\]
and we have 
\[
\coprod_{M_S} \mathfrak{A}^+_{M_S, \mathbb{Q}}=\overline{C}_{\mathbb{Q}}.
\]
There is a partial ordering of $\mathfrak{A}_{\mathbb{Q}}$ given by $\mu \preceq \mu'$ if $\mu' - \mu$ is a non-negative rational combination of simple roots.

\begin{definition}
We define a \emph{cocharacter pair} for a group $G$ (relative to some fixed choice of $T$ and $B$ defined over $\mathbb{Q}_p$) to be a pair $(M_S, \mu_S)$ such that $M_S \subset G$ is a standard Levi subgroup and $\mu_S \in X_*(T)$ satisfies $\langle \mu_S, \alpha \rangle \geq 0$ for each positive \emph{absolute} root $\alpha$ of $T$ in the Lie algebra of $M_{S, \overline{\mathbb{Q}_p}}$. Positivity for absolute roots is determined by the Borel subgroup $B$ which we have fixed.

We denote the set of cocharacter pairs for $G$ by $\mathcal{C}_G$.
\end{definition}
\begin{remark}
We caution the reader that the cocharacter $\mu_S$ need not be an element of $X_*(A)$, even though $M_S$ is defined over $\mathbb{Q}_p$.

We could define cocharacter pairs more canonically as the set of equivalence classes of pairs $(M, \mu)$ such that $M$ is a Levi subgroup of $G$ defined over $\mathbb{Q}_p$ and $\mu$ is a cocharacter of $M$. Two pairs $(M, \mu), (M', \mu')$ are equivalent if $M, M'$ are conjugate in $G_{\mathbb{Q}_p}$ and $\mu, \mu'$ are conjugate in $M_{\overline{\mathbb{Q}_p}}$. We choose not to do this as in practice we will often need to work with the unique dominant cocharacter in a conjugacy class relative to a fixed based root datum.
\end{remark}

Let $\Gamma= \mathrm{Gal}(\overline{\mathbb{Q}_p}/ \mathbb{Q}_p)$. Since we have assumed $T$ and $B$ are defined over $\mathbb{Q}_p$, $\Gamma$ acts on $T_{\overline{\mathbb{Q}_p}}$ and $B_{\overline{\mathbb{Q}_p}}$. This gives us a natural left action of $\Gamma$ on $X_*(T)$ given explicitly by $(\gamma \cdot \mu)(g)=\gamma(\mu(\gamma^{-1}(g))$ for $\mu \in X_*(T)$ and $\gamma \in \Gamma$. We get an analogous left action on $X^*(T)$ and one can easily check that the pairing $X^*(T) \times X_*(T) \to \mathbb{Z}$ is $\Gamma$ invariant under these actions.

We have 
\[
X_*(T)^{\Gamma}=\mathfrak{A}.
\]
Indeed, a $\Gamma$-invariant cocharacter $\mu$ factors through the identity component of $T^{\Gamma}$, where $T^{\Gamma}$ is the subscheme defined by $T^{\Gamma}(\overline{\mathbb{Q}_p})=T(\overline{\mathbb{Q}_p})^{\Gamma}$. But the identity component of $T^{\Gamma}$ is the torus $A$. Conversely any cocharacter of $A$ induces a $\Gamma$-invariant cocharacter via the natural inclusion $A \hookrightarrow T$.

Given $\mu \in X_*(T)$, we construct an element $\mu^{\Gamma}$ of $\mathfrak{A}_{\mathbb{Q}}$ as follows:
\[
\mu^{\Gamma}= \frac{1}{[\Gamma: \Gamma_{\mu}]} \sum\limits_{\gamma \in \Gamma/ \Gamma_{\mu}} \gamma(\mu)
\]
where $\Gamma_{\mu}$ is the stabilizer of $\mu$ in $\Gamma$. Then $\mu^{\Gamma} \in X_*(T)^{\Gamma}_{\mathbb{Q}} = \mathfrak{A}_{\mathbb{Q}}$. 

Given a standard Levi subgroup $M_S$, we let $W^{\mathrm{rel}}_{M_S}$ denote the relative Weyl group of $M_S$. The group $W^{\mathrm{rel}}_{M_S}$ is defined to be  the subgroup of the relative Weyl group, $W^{\mathrm{rel}}$, that is generated by the reflections corresponding to simple roots in $S$. 

\begin{definition}{\label{theta}}
We define a map 
\[
 \theta_{M_S}: X_*(T) \to  \mathfrak{A}_{\mathbb{Q}},
\]
given by
\[
 \theta_{M_S}( \mu)= \frac{1}{|W^{\mathrm{rel}}_{M_S}|} \sum\limits_{\sigma \in W^{\mathrm{rel}}_{M_S}} \sigma( \mu^{\Gamma}).
\]
\end{definition}

We are now ready to describe a formalism that will prove useful in studying the cohomology of certain Rapoport-Zink spaces. Crucial to everything that follows is a partial ordering on the set $\mathcal{C}_G$ of cocharacter pairs for $G$.

\begin{definition}\label{poset}
We define a partial ordering on $\mathcal{C}_G$ which we denote by the symbol $ \leq$. Unfortunately, our definition is somewhat indirect: we first define when $(M_{S_2}, \mu_{S_2}) \leq (M_{S_1}, \mu_{S_1})$ for $M_{S_2} \subset M_{S_1}$ (equivalently $S_2 \subset S_1$) and $S_1 \setminus S_2$ contains a single element (in other words, $M_{S_2}$ is a maximal proper Levi subgroup of $M_{S_1}$). We then extend the relation to all cocharacter pairs by taking the transitive closure.

Let $M_{S_2}, M_{S_1}$ be standard Levi subgroups of $G$ such that $M_{S_2} \subset M_{S_1}$ and $S_1 \setminus S_2$ is a singleton. For cocharacter pairs $(M_{S_2}, \mu_{S_2}), (M_{S_1}, \mu_{S_1}) \in \mathcal{C}_G$, we write $(M_{S_2}, \mu_{S_2}) \leq (M_{S_1}, \mu_{S_1})$ if  $\mu_{S_2}$ is conjugate to $\mu_{S_1}$ in ${M_{S_1}}_{\overline{\mathbb{Q}_p}}$ and  $\theta_{M_{S_2}}(\mu_{S_2}) \succ \theta_{M_{S_1}}(\mu_{S_1})$. We then take the transitive closure to extend to a partial ordering on $\mathcal{C}_G$.
\end{definition}
The following example shows that the above definition  depends on the assumption that $S_1 \setminus S_2$ is a singleton.
\begin{example}
Consider $G= \mathrm{GL}_4$ with $T$ the diagonal torus and $B$ the upper triangular matrices. We can pick a basis for $X_*(T)$ of cocharacters $\widehat{e_i}$ defined so that $\widehat{e_i}(g)$ is the diagonal matrix with $1$ in every position except for the $i$th, which equals $g$. Then we can identify an element of $X_*(T)$ with its coordinate vector in this basis. Finally, we use additional parenthesis to indicate the product structure of the standard Levi subgroup $M_S$. Using this notation, the set of cocharacter pairs that are less than or equal to $(\mathrm{GL}_4, (1^2, 0^2))$ is given in the diagram at the start of Appendix \ref{A}.

In particular, we see that $(\mathrm{GL}^4_1, (1)(1)(0)(0)) \leq (\mathrm{GL}_4, (1^2, 0^2))$ since we have a chain of cocharacter pairs where each Levi subgroup is maximal in the next: 
\begin{align*}
&(\mathrm{GL}^4_1, (1)(1)(0)(0)) \leq (\mathrm{GL}_1 \times \mathrm{GL}_2 \times \mathrm{GL}_1, (1)(1,0)(0))\\ &\leq (\mathrm{GL}_3 \times \mathrm{GL}_1, (1^2,0)(0)) \leq (\mathrm{GL}_4, (1^2, 0^2)).
\end{align*}
However, it is not the case that $(\mathrm{GL}^4_1, (1)(0)(1)(0)) \leq (\mathrm{GL}_4, (1^2,0^2))$ even though $\theta_{\mathrm{GL}^4_1}((1,0,1,0)) \succ \theta_{\mathrm{GL}_4}( (1,1, 0,0) )$ and the cocharacters are conjugate in $G$.

Finally, we remark that the fact that all the related  cocharacter pairs in the above example have equal (as opposed to just conjugate) cocharacters is very much a result of us choosing a fairly small group $G$. Even for $G=\mathrm{GL}_5$, this is not the case. 
\end{example}

\begin{definition}
We define a cocharacter pair $(M_S, \mu_S)$  for $G$ to be \emph{strictly decreasing} if $\theta_{M_S}(\mu_S) \in \mathfrak{A}^+_{M_S, \mathbb{Q}}$. We denote by $\mathcal{SD} \subset \mathcal{C}_G$ the strictly decreasing elements of $\mathcal{C}_G$ and by $\mathcal{SD}_{\mu}$ (for dominant $\mu \in X_*(T)$) the strictly decreasing elements $(M_S, \mu_S) \in \mathcal{C}_G$ such that $(M_S, \mu_S) \leq (G, \mu)$.
\end{definition}
\begin{remark}
 The $\theta_{M_S}$ map can be thought of as associating a tuple of slopes to a cocharacter pair. Then the strictly decreasing cocharacter pairs with Levi subgroup $M_S$ are the ones whose slope tuple lies in the image of the Newton map $\nu: \mathbf{B}(G)_{M_S} \to \mathfrak{A}_{M_S, \mathbb{Q}}$. The above statement is made precise by Proposition \ref{2.16}.
\end{remark}

\subsection{An Alternate Characterization of the Averaging Map} 

The following two subsections consist of a collection of lemmas developing the theory of the map $\theta_{M_S}$ and the set of strictly decreasing elements $\mathcal{SD}$ of $\mathcal{C}_G$. 

In this section, we give an alternate description of the map $\theta_{M_S}$. To do so, we will need several properties of cocharacters and root data which we record in the following lemma. For this lemma only, we consider $T$ and $G$ defined over a more general class of fields so that these results also apply to the complex dual groups $\widehat{T}$ and $\widehat{G}$.
\begin{lemma} \label{2.11}
 Let $F \supset \mathbb{Q}$ be a field and $\overline{F}$ an algebraic closure. Let $G$ be a connected quasisplit reductive group defined over $F$. Suppose that $T \subset G$ is a maximal torus defined over $F$ and that the group scheme $T_{\overline{F}}$ admits an action defined over $\overline{F}$ by a finite group $\Lambda$. Let $X^*(T^{\Lambda})$ denote the characters of the subgroup scheme of $\Lambda$-fixed points of $T_{\overline{F}}$. The anti-equivalence of categories between tori and finitely generated free Abelian groups given by $T_{\overline{F}} \mapsto X^*(T)$ induces an action of $\Lambda$ on $X^*(T)$. We then have the following.
\begin{enumerate}
\item There is a unique isomorphism $X^*(T^{\Lambda}) \cong X^*(T)_{\Lambda}$ such that the following diagram commutes.

\begin{tikzcd}
X^*(T) \arrow[r, "res"] \arrow[dr, two heads,swap, "proj" ] & X^*(T^{\Lambda}) \arrow[d, leftrightarrow] \\
& X^*(T)_{\Lambda}
\end{tikzcd}
\item Let $M_S \subset G$ be a standard Levi subgroup. Let $W^{\mathrm{abs}}_{M_S}, W^{\mathrm{rel}}_{M_S}$ denote the absolute and relative Weyl groups of $M_S$ and let $\Gamma= \mathrm{Gal}(\overline{F}/ F)$.  Then $W_{M_S, \mathrm{rel}}$ acts on $X_*(T)^{\Gamma}$ via its natural identification with $\mathfrak{A}$ and $\Gamma$ acts on $X_*(T)^{W_{M_S,\mathrm{abs}}}$ since for $w \in W_{M_S,\mathrm{abs}}$, and $\gamma \in \Gamma$, and $\mu \in X_*(T)^{W_{M_S,\mathrm{abs}}}$, we have $w(\gamma(\mu))=\gamma(\gamma^{-1}(w)(\mu))=\gamma(\mu)$. Then the identity map on $X_*(T)$ induces an isomorphism of groups
\[
(X_*(T)^{W_{M_S, \mathrm{abs}}})^{\Gamma} \cong (X_*(T)^{\Gamma})^{W_{M_S, \mathrm{rel}}}
\]
\item The natural map $X_*(T)^{\Lambda}_{\mathbb{Q}} \hookrightarrow X_*(T)_{\mathbb{Q}} \twoheadrightarrow X_*(T)_{\mathbb{Q}, \Lambda}$ induces an isomorphism $X_*(T)^{\Lambda}_{\mathbb{Q}} \cong X_*(T)_{\Lambda, \mathbb{Q}}$.
\end{enumerate}
\end{lemma}
\begin{proof}

The functor $T \mapsto X^*(T)$ is an anti-equivalence between the categories of diagonalizable groups over $\overline{F}$ and finitely generated Abelian groups. The diagram for the universal property for $\Lambda$-invariants is that of $\Lambda$-coinvariants but with all the arrows reversed. Thus, there must exist a unique isomorphism between $X^*(T^{\Lambda})$ and $X^*(T)_{\Lambda}$ that makes the diagram
\[
\begin{tikzcd}
X^*(T) \arrow[r, "res"] \arrow[dr, two heads, swap, "proj" ] & X^*(T^{\Lambda}) \arrow[d, leftrightarrow] \\
& X^*(T)_{\Lambda}
\end{tikzcd}
\]
 commute. This proves $(1)$.

In \cite[Lem 1.1.3]{Kot3}, Kottwitz proves that the identity map on $X_*(T)$ induces an isomorphism \[(X_*(T)^{\Gamma})/W^{\mathrm{rel}}_{M_S} \cong (X_*(T)/ W^{\mathrm{abs}}_{M_S})^{\Gamma}.\] Thus, to prove $(2)$, we need only show that this isomorphism gives a bijection of the singleton orbits. This will give an isomorphism of groups (not just sets) between $(X_*(T)^{W_{M_S, \mathrm{abs}}})^{\Gamma}$ and $ (X_*(T)^{\Gamma})^{W_{M_S, \mathrm{rel}}}$ that is induced from the identity map on $X_*(T)$.

Kottwitz's isomorphism maps the $W^{\mathrm{rel}}_{M_S}$-orbit of $\mu \in X_*(T)^{\Gamma}$ to its $W^{\mathrm{abs}}_{M_S}$ orbit in $X_*(T)$. Thus, it suffices to show that if $\mu \in X_*(T)^{\Gamma}$ is invariant by $W^{\mathrm{rel}}_{M_S}$ then it is also invariant by $W^{\mathrm{abs}}_{M_S}$. If $\mu$ is invariant by $W^{\mathrm{rel}}_{M_S}$, then the pairing of $\mu$ with each relative root of $M_S$ is $0$. Thus the image of $\mu$ lies in the intersection of the kernels of the relative roots of $M_S$ which is $Z(M_S) \cap A$. Therefore, $\mu$ is invariant under the action of $W^{\mathrm{abs}}_{M_S}$.

Finally, we note that the proof of Kottwitz uses the fact that the intersection of the absolute Weyl chamber $\overline{C}^{\mathrm{abs}}_{\mathbb{Q}}$ with the image of $X_*(A)$ in $X_*(T)$ gives the relative Weyl chamber $\overline{C}_{\mathbb{Q}}$. Indeed, this follows easily from the fact that the restriction of the set of absolute simple roots $\Delta^{\mathrm{abs}}$ relative to our choice of $B$ and $T$ equals the set of relative simple roots $\Delta$ (see Proposition \ref{rootrestrict}). An analogous fact is known for the Weyl chambers in the character group $X^*(T)$ (see Proposition \ref{chamberresprop}) but this seems to be much more subtle.

For $(3)$, we need to construct an inverse to the map
\[
X_*(T)^{\Lambda}_{\mathbb{Q}} \hookrightarrow X_*(T)_{\mathbb{Q}} \twoheadrightarrow X_*(T)_{\mathbb{Q}, \Lambda}.
\]
Take $[\mu] \in X_*(T)_{\mathbb{Q}, \Lambda}$ for $\mu \in X_*(T)_{\mathbb{Q}}$. Then 
\[
\frac{1}{\Lambda} \sum\limits_{\lambda \in \Lambda} \lambda(\mu) \in X_*(T)^{\Lambda}_{\mathbb{Q}} 
\]
is independent of the choice of lift of $[\mu]$ to $X_*(T)_{\mathbb{Q}}$ and gives an inverse to the map above.
\end{proof}

Let $A_{M_S}$ be the maximal split torus in the center of $M_S$. Then 
\[
X_*(A_{M_S})_{\mathbb{Q}} \cong \mathfrak{A}_{M_S, \mathbb{Q}}.
\]
We now prove a lemma that we will need to use to describe the alternate characterization of $\theta_{M_S}$.
\begin{lemma}\label{lemma 2.12}
%% This first isomorphism is not so obvious it comes from the fact that the weyl invariants of $\hat{T}$ and $Z(G)$ have the same connected component (but are not always equal) Moreover I think the iso between x_*(T)_Q and X_*(A)_Q also requires the fact that we're tensoring with Q. 
\begin{enumerate}
    \item There is a natural isomorphism $X^*(Z(\widehat{M_S})^{\Gamma})_{\mathbb{Q}} \cong \mathfrak{A}_{M_S, \mathbb{Q}}$ defined via a series of canonical identifications.\\
    \item The isomorphism in $(1)$ coincides with the one constructed in \S 4.4.3 of \cite{Kot1}.
\end{enumerate}
\end{lemma}
\begin{proof}
We prove $(1)$ first. By Lemma \ref{2.11}, we have the following isomorphisms.
\begin{align*}
X^*(\widehat{T}^{W^{\mathrm{abs}}_{M_S}, \Gamma})_{\mathbb{Q}} &\cong  X^*(\widehat{T})_{\mathbb{Q}, W^{\mathrm{abs}}_{M_S}, \Gamma}
=X_*(T)_{\mathbb{Q}, W^{\mathrm{abs}}_{M_S}, \Gamma}\\
&\cong X_*(T)^{W^{\mathrm{abs}}_{M_S}, \Gamma}_{\mathbb{Q}} \,\ \,\ \cong  X_*(T)^{\Gamma, W^{\mathrm{rel}}_{M_S}}_{\mathbb{Q}}\\ & \cong X_*(A_{M_S})_{\mathbb{Q}} \,\ \,\ \,\ \cong \mathfrak{A}_{M_S, \mathbb{Q}}.
\end{align*}
We explicate the isomorphism $X_*(T)^{\Gamma, W^{\mathrm{rel}}_{M_S}}_{\mathbb{Q}} \cong  X_*(A_{M_S})_{\mathbb{Q}}$. This follows from the isomorphism $X_*(A)^{W^{\mathrm{rel}}_{M_S}} \cong X_*(A_{M_S})$ which we now describe. Suppose we have $\mu \in X_*(A)^{W^{\mathrm{rel}}_{M_S}}$. Equivalently, for each relative root $\alpha$ of $\mathrm{Lie}(M_S)$, we have $\sigma_{\alpha}(\mu)=\mu$ (where $\sigma_{\alpha}$ is the reflection in the Weyl group corresponding to $\alpha$). Since $\sigma_{\alpha}(\mu)=\mu-\langle \mu, \alpha \rangle \check{\alpha}$, this is equivalent to $\langle \mu, \alpha \rangle = 0$ for all relative roots $\alpha$ of $\mathrm{Lie}(M_S)$, which in turn is equivalent to the statement that $\mathrm{im}(\mu) \subset \bigcap\limits_{\alpha} \mathrm{ker} \alpha$. Finally, this is equivalent to $\mathrm{im}(\mu) \subset Z(M_S) \cap A$. Since the image of a cocharacter is connected, we in fact have that $\mu \in X_*(A_{M_S})$.

To finish the argument, we need to construct an isomorphism
\[
X^*(Z(\widehat{M_S})^{\Gamma})_{\mathbb{Q}} \cong X^*(\widehat{T}^{W^{\mathrm{abs}}_{M_S}, \Gamma})_{\mathbb{Q}}.
\]
Note that it is necessary to take the tensor product with $\mathbb{Q}$ here as $Z(\widehat{M_S})$ and $\widehat{T}^{W^{\mathrm{abs}}_{M_S}}$ need not be isomorphic.

It suffices to show that 
\[
X^*(Z(\widehat{M_S}))_{\mathbb{Q}} \cong X^*(\widehat{T}^{W^{\mathrm{abs}}_{M_S}})_{\mathbb{Q}}.
\]
The group $Z(\widehat{M_S})$ is equal to the intersection of the kernels of the roots of $\widehat{M_S}$ and so $X^*(Z(\widehat{M_S}))$ is identified with $X^*(\widehat{T})/R$ where $R$ is the $\mathbb{Z}$-module spanned by the roots of $\widehat{M_S}$. By Lemma \ref{2.11}, $X^*(\widehat{T}^{W^{\mathrm{abs}}_{M_S}}) \cong X^*(\widehat{T})_{W^{\mathrm{abs}}_{M_S}}=X^*(\widehat{T})/D$ where $D$ is the $\mathbb{Z}$ module spanned by $w(\mu)-\mu$ for every $w \in W^{\mathrm{abs}}_{M_S}$ and $\mu \in X^*(\widehat{T})$. Since $Z(\widehat{M_S}) \subset \widehat{T}^{W^{\mathrm{abs}}_{M_S}}$, we have a natural surjection
\[
X^*(\widehat{T}^{W^{\mathrm{abs}}_{M_S}}) \twoheadrightarrow X^*(Z(\widehat{M_S})).
\]
By our previous discussion, the kernel of this map is $R/D$. Thus, to prove our claim, it suffices to show that $R/D$ is finite. But if $\alpha$ is a root of $\widehat{M_S}$, then $\sigma_{\alpha}(\alpha)-\alpha=-2 \alpha$. Thus $2R \subset D$ and so we have the desired result.

We now show $(2)$. The map in \cite[\S 4.4.3]{Kot1} is defined as follows:
\[
\mathfrak{A}_{M_S, \mathbb{Q}} \to X_*(T)_{\mathbb{Q}} =X^*(\widehat{T})_{\mathbb{Q}} \xrightarrow{res} X^*(Z( \widehat{M_S})^{\Gamma})_{\mathbb{Q}},
\]
where the final map is restriction of characters. By Lemma \ref{2.11} (1), this last map is the same as the composition
\[
X^*(\widehat{T})_{\mathbb{Q}} \to X^*(\widehat{T})_{\mathbb{Q}, W^{\mathrm{abs}}_{M_S}, \Gamma} \cong X^*(\widehat{T}^{W^{\mathrm{abs}}_{M_S}, \Gamma})_{\mathbb{Q}} \cong X^*(Z( \widehat{M_S})^{\Gamma})_{\mathbb{Q}},
\]
Thus, by applying Lemma \ref{2.11} and the proof of Lemma \ref{lemma 2.12}, we get that the entire map is given by
\[
\mathfrak{A}_{M_S, \mathbb{Q}} \cong X_*(T)^{\Gamma, W^{\mathrm{rel}}_{M_S}}_{\mathbb{Q}} \cong X_*(T)^{W^{\mathrm{abs}}_{M_S}, \Gamma}_{\mathbb{Q}} \cong  X_*(T)_{\mathbb{Q}, W^{\mathrm{abs}}_{M_S}, \Gamma},
\]
\[
\cong X^*(\widehat{T}^{W^{\mathrm{abs}}_{M_S}, \Gamma})_{\mathbb{Q}} \cong X^*(Z( \widehat{M_S})^{\Gamma})_{\mathbb{Q}}.
\]
We observe that this is the inverse of what we wrote down above. 
\end{proof}
We are now ready to give our alternate characterization of the map $\theta_{M_S}$.
\begin{proposition}\label{2.14}[Alternate Characterization of $\theta_{M_S}$]
The map $\theta_{M_S}$ that was introduced in Definition \ref{theta} is equal to the composition
\begin{align*}
X_*(T)=X^*(\widehat{T}) \xrightarrow{res} X^*(Z(\widehat{M_S})^{\Gamma}) \to  X^*(Z(\widehat{M_S})^{\Gamma})_{\mathbb{Q}} \cong \mathfrak{A}_{M_S, \mathbb{Q}} \subset \mathfrak{A}_{\mathbb{Q}},
\end{align*}
where the final isomorphism is the one described in Lemma \ref{lemma 2.12}.
\end{proposition}
\begin{proof}
We recall Definition \ref{theta} where $\theta_{M_S}$ is defined to be the composition
\[
X_*(T) \to X_*(T)^{\Gamma}_{\mathbb{Q}} \to X_*(T)^{\Gamma, W^{\mathrm{rel}}_{M_S}}_{\mathbb{Q}} \subset \mathfrak{A}_{\mathbb{Q}},
\]
where both maps are averages over the relevant group. As we now show, this is the same as the composition
\[
X_*(T) \to X_*(T)^{W^{\mathrm{abs}}_{M_S}}_{\mathbb{Q}} \to X_*(T)^{W^{\mathrm{abs}}_{M_S}, \Gamma}_{\mathbb{Q}} \cong  X_*(T)^{\Gamma, W^{\mathrm{rel}}_{M_S}}_{\mathbb{Q}} \subset \mathfrak{A}_{\mathbb{Q}},
\]
where the first two maps are averages and the third is as in Lemma \ref{2.11} (2). Indeed for $\mu \in X_*(T)$,
\[
\frac{1}{|W^{\mathrm{rel}}_{M_S}|}\sum\limits_{w \in W^{\mathrm{rel}}_{M_S}} \sum\limits_{\gamma \in \Gamma} w(\gamma(\mu)),
\]
is invariant by $W^{\mathrm{abs}}_{M_S}$ by Lemma \ref{2.11} (2) and so equals (keeping in mind that $W^{\mathrm{rel}}_{M_S} \subset W^{\mathrm{abs}}_{M_S}$ by Corollary \ref{absrelweyl})
\[
\frac{1}{|W^{\mathrm{abs}}_{M_S}|}\sum\limits_{w \in W^{\mathrm{abs}}_{M_S}} \sum\limits_{\gamma \in \Gamma} w(\gamma(\mu))=\frac{1}{|W^{\mathrm{abs}}_{M_S}|}\sum\limits_{w \in W^{\mathrm{abs}}_{M_S}} \sum\limits_{\gamma \in \Gamma} \gamma(w)(\gamma(\mu))
\]
\[
=\frac{1}{|W^{\mathrm{abs}}_{M_S}|}\sum\limits_{w \in W^{\mathrm{abs}}_{M_S}} \sum\limits_{\gamma \in \Gamma} \gamma(w(\mu))=\frac{1}{|W^{\mathrm{abs}}_{M_S}|}\sum\limits_{\gamma \in \Gamma}\sum\limits_{w \in W^{\mathrm{abs}}_{M_S}} \gamma(w(\mu)).
\]

Now, we consider the following commutative diagram.
\[
\begin{adjustbox}{max size={.95\textwidth}{.8\textheight}}
\begin{tikzcd}
X_*(T)_{\mathbb{Q}} \arrow[rr, "avg" ] \arrow[rd] && X_*(T)^{W^{\mathrm{abs}}_{M_S}}_{\mathbb{Q}} \arrow[rr, "avg" ] \arrow[rd] && X_*(T)^{W^{\mathrm{abs}}_{M_S}, \Gamma}_{\mathbb{Q}} \\
& X_*(T)_{\mathbb{Q}, W^{\mathrm{abs}}_{M_S}} \arrow[ru, "avg"] \arrow[rd] & & X_*(T)^{W^{\mathrm{abs}}_{M_S}}_{\mathbb{Q}, \Gamma} \arrow[ru, "avg"]& \\
&& X_*(T)_{\mathbb{Q}, W^{\mathrm{abs}}_{M_S}, \Gamma} \arrow[ru, "avg"] &&
\end{tikzcd}
\end{adjustbox}
\]

The commutativity essentially follows from the definition of the averaging maps. The benefit of this is that now we can write $\theta_{M_S}$ as the composition of 
\begin{align*}
X_*(T) & \to X_*(T)_{W^{\mathrm{abs}}_{M_S}} \to X_*(T)_{W^{\mathrm{abs}}_{M_S}, \Gamma} \to X_*(T)_{\mathbb{Q}, W^{\mathrm{abs}}_{M_S}, \Gamma} \\
& \to X^*(T)^{W^{\mathrm{abs}}_{M_S}}_{\mathbb{Q}, \Gamma} \to X_*(T)^{W^{\mathrm{abs}}_{M_S}, \Gamma}_{\mathbb{Q}} \cong X_*(T)^{\Gamma, W^{\mathrm{rel}}_{M_S}}  \subset \mathfrak{A}_{\mathbb{Q}}
\end{align*}
where we no longer need to base change the first three spaces to $\mathbb{Q}$ because denominators are not introduced in the maps until later.

Using the equality between cocharacters of $T$ and characters of $\widehat{T}$, we rewrite this as
\begin{align*}
X_*(T) &=X^*(\widehat{T}) \to X^*(\widehat{T})_{W^{\mathrm{abs}}_{M_S}} \to X^*(\widehat{T})_{W^{\mathrm{abs}}_{M_S}, \Gamma} \to X^*(\widehat{T})_{\mathbb{Q},W^{\mathrm{abs}}_{M_S}, \Gamma}\\ & \to X^*(\widehat{T})^{W^{\mathrm{abs}}_{M_S}}_{\mathbb{Q}, \Gamma} \to X^*(\widehat{T})^{W^{\mathrm{abs}}_{M_S}, \Gamma}_{\mathbb{Q}}= X_*(T)^{W^{\mathrm{abs}}_{M_S}, \Gamma}_{\mathbb{Q}} \cong X_*(T)^{\Gamma, W^{\mathrm{rel}}_{M_S}} \subset \mathfrak{A}_{\mathbb{Q}}.
\end{align*}
Now we invoke Lemma \ref{2.11} $(1)$ to get that the above composition is equal to
\begin{align*}
X_*(T) & =X^*(\widehat{T}) \xrightarrow{res} X^*(\widehat{T}^{W^{\mathrm{abs}}_{M_S}, \Gamma}) \to X^*(\widehat{T}^{W^{\mathrm{abs}}_{M_S}, \Gamma})_{\mathbb{Q}} \cong X^*(\widehat{T})_{\mathbb{Q}, W^{\mathrm{abs}}_{M_S}, \Gamma}\\ & \to X^*(\widehat{T})^{W^{\mathrm{abs}}_{M_S}}_{\mathbb{Q}, \Gamma} \to X^*(\widehat{T})^{W^{\mathrm{abs}}_{M_S}, \Gamma}_{\mathbb{Q}}= X_*(T)^{W^{\mathrm{abs}}_{M_S}, \Gamma}_{\mathbb{Q}} \cong X_*(T)^{\Gamma, W^{\mathrm{rel}}_{M_S}} \subset \mathfrak{A}_{\mathbb{Q}}.
\end{align*}
The final step is to observe that we have a commutative diagram
\[
\begin{tikzcd}
X^*(\widehat{T}^{W^{\mathrm{abs}}_{M_S}, \Gamma }) \arrow[d, "res"] \arrow[r] & X^*(\widehat{T}^{W^{\mathrm{abs}}_{M_S}, \Gamma })_{\mathbb{Q}}  \isoarrow{d}\\
X^*(Z(\widehat{M_S})^{\Gamma}) \arrow[r] & X^*(Z(\widehat{M_S})^{\Gamma})_{\mathbb{Q}}.
\end{tikzcd}
\]
Thus, the previous expression equals
\begin{align*}
X_*(T) & =X^*(\widehat{T}) \xrightarrow{res} X^*(\widehat{T}^{W^{\mathrm{abs}}_{M_S}, \Gamma}) \xrightarrow{res} X^*(Z(\widehat{M_S})^{\Gamma}) \to X^*(Z(\widehat{M_S})^{\Gamma})_{\mathbb{Q}}\\ & \cong X^*(\widehat{T}^{W^{\mathrm{abs}}_{M_S}, \Gamma})_{\mathbb{Q}} \cong X^*(\widehat{T})_{\mathbb{Q}, W^{\mathrm{abs}}_{M_S}, \Gamma} \to X^*(\widehat{T})^{W^{\mathrm{abs}}_{M_S}}_{\mathbb{Q}, \Gamma}\\ & \to X^*(\widehat{T})^{W^{\mathrm{abs}}_{M_S}, \Gamma}_{\mathbb{Q}}= X_*(T)^{W^{\mathrm{abs}}_{M_S}, \Gamma}_{\mathbb{Q}} \cong X_*(T)^{\Gamma, W^{\mathrm{rel}}_{M_S}} \subset \mathfrak{A}_{\mathbb{Q}}.
\end{align*}
comparing with Lemma \ref{lemma 2.12}, we can rewrite $\theta_{M_S}$ as 
\begin{align*}
X_*(T)=X^*(\widehat{T}) \xrightarrow{res} X^*(Z(\widehat{M_S})^{\Gamma}) \to X^*(Z(\widehat{M_S})^{\Gamma})_{\mathbb{Q}} \cong \mathfrak{A}_{M_S ,\mathbb{Q}} \subset \mathfrak{A}_{\mathbb{Q}}
\end{align*}
as desired.

\end{proof}
We record the following useful corollary of the ideas discussed in the above argument.
\begin{corollary}{\label{conjinv}}
Suppose that $\mu, \mu' \in X_*(T)$ are conjugate in $M_{S, \overline{\mathbb{Q}_p}}$. Then $\theta_{M_S}(\mu)=\theta_{M_S}(\mu')$.
\end{corollary}
\begin{proof}
By the observation at the start of Proposition \ref{2.14}, $\theta_{M_S}$ is equivalently defined as the composition
\[
X_*(T) \to X_*(T)^{W^{\mathrm{abs}}_{M_S}}_{\mathbb{Q}} \to X_*(T)^{W^{\mathrm{abs}}_{M_S}, \Gamma}_{\mathbb{Q}} \cong  X_*(T)^{\Gamma, W^{\mathrm{rel}}_{M_S}}_{\mathbb{Q}} \subset \mathfrak{A}_{\mathbb{Q}}.
\]
In particular, $\mu$ and $\mu'$ are mapped to the same element under the first map in the above composition.
\end{proof}
\subsection{Strictly Decreasing Cocharacter Pairs}
 In this section, we prove a number of properties of strictly decreasing cocharacter pairs and their relation to the partial order we defined in Definition \ref{poset}. As always, we let $\sigma_{\alpha}$ denote the reflection in the relative Weyl group corresponding to the relative root $\alpha$.
\begin{lemma}\label{2.8}
If $x \in \mathfrak{A}_{\mathbb{Q}}$ is dominant, then 
\[
y=\frac{1}{|W^{\mathrm{rel}}_{M_S}|} \sum\limits_{\sigma \in W^{\mathrm{rel}}_{M_S}} \sigma(x)
\]
is also dominant. If in addition, $\langle x, \alpha \rangle > 0$  for some $\alpha \in \Delta \setminus S$, then we also have $\langle y, \alpha \rangle >0$.
\end{lemma}
\begin{proof} 
For the first part of the lemma, we claim that if we can show that $\langle \sigma(x), \alpha \rangle  \geq 0$ for each $\sigma \in W^{\mathrm{rel}}_{M_S}$ and $\alpha \in \Delta \setminus S$, then we are done. This follows because if a collection of cocharacters pair non-negatively with $\alpha$, then so will their average. Thus for $\alpha \in \Delta \setminus S$, we get  $\langle y, \alpha \rangle \geq 0$. For $\alpha \in S$, we automatically have $\langle y, \alpha \rangle =0$ since $0=y-\sigma_{\alpha}(y)=\langle y, \alpha \rangle \check{\alpha}$. 

Pick $\alpha \in \Delta \setminus S$. Then the root group of $\alpha$ is contained in the unipotent radical $N_S$ of $P_S$. The group $N_S$ is normalized by $M_S$. In particular, for any $\sigma \in W^{\mathrm{rel}}_{M_S}$, the root group of $\sigma^{-1}(\alpha)$ is contained in $N_S$ and hence $\sigma^{-1}(\alpha)$ is also a positive root. Thus $\langle \sigma(x), \alpha \rangle= \langle x, \sigma^{-1}(\alpha) \rangle \geq 0$ as desired.

To prove the second part, we notice since $\langle x, \alpha \rangle >0$, the term in $y$ corresponding to $\sigma=1$ has positive pairing with $\alpha$. Since all the other terms have non-negative pairing with $\alpha$, we must have that $\langle y, \alpha \rangle > 0$.
\end{proof}
\begin{lemma} \label{2.9}
If $x$ as in the previous lemma is dominant, then 
\[
\frac{1}{|W^{\mathrm{rel}}_{M_S}|} \sum\limits_{\sigma \in W^{\mathrm{rel}}_{M_S}} \sigma(x) \preceq x
\]
\end{lemma}
\begin{proof}
It suffices to show that for any $\sigma \in W^{\mathrm{rel}}_{M_S}$, we have $\sigma(x) \preceq x$. This is a standard fact (\cite[Ch6 1.6.18, p. 158]{Bou1}).
\end{proof}
\begin{corollary}{\label{sdupclsd}}
Let $(M_S, \mu_S) \in \mathcal{SD}$ be a strictly decreasing cocharacter pair and let $(M_{S'}, \mu_{S'}) \in \mathcal{C}_G$ and suppose that $(M_S, \mu_S) \leq (M_{S'}, \mu_{S'})$. Then $(M_{S'}, \mu_{S'}) \in \mathcal{SD}$.
\end{corollary}
\begin{proof}
We need to show that for each $\beta \in \Delta \setminus S'$, that $\langle \theta_{M_{S'}}(\mu_{S'}), \beta \rangle >0$. By \ref{conjinv}, $\theta_{M_{S'}}(\mu_{S'})=\theta_{M_{S'}}(\mu_S)$. Further, we observe that
\begin{align}{\label{kapavg}}
\theta_{M_{S'}}(\mu_{S})= \frac{1}{|W^{\mathrm{rel}}_{M_{S'}}|} \sum\limits_{\sigma \in W^{\mathrm{rel}}_{M_{S'}}} \sigma( \theta_{M_S}( \mu_S)).
\end{align}
Since $\theta_{M_S}(\mu_S)$ is dominant by assumption and satisfies $\langle \theta_{M_S}(\mu_S), \beta \rangle > 0$, we can apply \ref{2.8} to get the desired result.
\end{proof}

The following easy uniqueness result is quite useful.
\begin{lemma} \label{uniqueness}
Let $(M_{S_1}, \mu_{S_1}), (M_{S_2}, \mu_{S_2}), (M_{S'_2}, \mu_{S'_2}) \in \mathcal{C}_G$. Suppose further that $(M_{S_1}, \mu_{S_1}) \leq (M_{S_2}, \mu_{S_2})$, that $(M_{S_1}, \mu_{S_1}) \leq (M_{S'_2}, \mu_{S'_2})$. If $M_{S_2}=M_{S'_2}$, then $(M_{S_2}, \mu_{S_2}) =(M_{S'_2}, \mu_{S'_2})$.
\end{lemma}
\begin{proof}
By definition, $\mu_{S_1}, \mu_{S_2}, \mu_{S'_2}$ are all conjugate in $M_{S_2}$. But also, $\mu_{S_2}$ and $\mu_{S'_2}$ are dominant in the absolute root system. Thus they are equal.
\end{proof}

We now define the notion of a cocharacter pair being strictly decreasing relative to a Levi subgroup.
\begin{definition}
Let $M_S \subsetneq M_{S'}$ be standard Levi subgroups of $G$. We say $(M_S, \mu_S)$ is strictly decreasing relative to $M_{S'}$ if $ \langle \theta_{M_S}( \mu_S) , \alpha \rangle > 0$ for $\alpha  \in S' \setminus S$.
\end{definition}
\begin{remark}
Recall that by construction, $\langle \theta_{M_S}(\mu_S) , \alpha \rangle =0$ for $\alpha \in S$. Thus, $(M_S, \mu_S) \in \mathcal{SD}$ exactly when it is strictly decreasing relative to $G$.
\end{remark}
\begin{lemma} \label{2.22}
Let $(M_{S_1}, \mu_{S_1}), (M_{S'_1}, \mu_{S'_1}) \in \mathcal{C}_G$ be cocharacter pairs such that $(M_{S_1}, \mu_{S_1}) \leq (M_{S'_1}, \mu_{S'_1})$. Let $M_{S_2} \supset M_{S_1}$ be a standard Levi subgroup of $G$ and suppose $(M_{S_1}, \mu_{S_1})$ is strictly decreasing relative to $M_{S_2}$. Then $(M_{S'_1}, \mu_{S'_1})$ is strictly decreasing relative to $M_{S'_1 \cup S_2}$.
\end{lemma}
\begin{proof}
We first reduce to the case where $M_{S_1}$ is a maximal Levi subgroup of $M_{S'_1}$ (i.e. $S'_1=S_1 \cup \{\alpha\}$ for some $\alpha \in \Delta \setminus S_1$). To do so, we recognize that the relation $(M_{S_1}, \mu_{S_1}) \leq (M_{S'_1}, \mu_{S'_1})$ definitionally implies that there is a finite sequence of cocharacter pairs
\[
(M_{S_1}, \mu_{S_1})=(M_{S^0}, \mu_{S^0}) \leq ... \leq (M_{S^k}, \mu_{S^k})=(M_{S'_1}, \mu_{S'_1})
\]
where each $M_{S^i}$ is a maximal Levi subgroup of $M_{S^{i+1}}$. Thus, if we prove the lemma in the maximal Levi subgroup case, we can inductively prove it in the general case.

We now assume that $M_{S_1} \subset M_{S'_1}$ is a maximal Levi subgroup so that $S'_1=S_1 \cup \{ \alpha \}$ for some $\alpha \in \Delta \setminus S_1$. We need to show that $\langle \theta_{M_{S'_1}}( \mu_{S'_1}), \beta \rangle >0$ for each $\beta \in S'_1 \cup S_2 \setminus S'_1$. First note that any such $\beta$ is an element of $S_2 \setminus S_1$. By Corollary \ref{conjinv}, since $\mu_{S_1}$ and  $\mu_{S'_1}$ are conjugate in $M_{S'_1}$, we have $\theta_{M_{S'_1}}(\mu_{S_1})=\theta_{M_{S'_1}}(\mu_{S'_1})$. Thus we are reduced to showing $\langle  \theta_{M_{S'_1}}( \mu_{S_1}), \beta \rangle >0$ for $\beta \in S_2 \setminus S_1$.

Note that since $(M_{S_1}, \mu_{S_1})$ is strictly decreasing relative to $M_{S_2}$, we have \\
$\theta_{M_{S_1}}(\mu_{S_1})$ is dominant relative to the root datum of $M_{S_2}$ and $\langle \theta_{M_{S_1}}(\mu_{S_1}), \beta \rangle >0$. Therefore, by Equation \eqref{kapavg} and Lemma \ref{2.8}, $\langle \theta_{M_{S'_1}}(\mu_{S_1}), \beta \rangle >0$ as desired.
\end{proof}
\begin{proposition}{\label{ext}}
Let $(M_S, \mu_S) \in \mathcal{C}_G$ and suppose it is strictly decreasing relative to some standard Levi subgroup $M_{S'} \supset M_{S}$. Then there is a unique $(M_{S'}, \mu_{S'}) \in \mathcal{C}_G$ such that $(M_S, \mu_S) \leq (M_{S'}, \mu_{S'})$. We call $(M_{S'}, \mu_{S'})$ the \emph{extension} of $(M_S, \mu_S)$ to $M_{S'}$.

In the case where $S'= S \cup \{\alpha\}$ for $\alpha \in \Delta \setminus S$, the converse is true. Specifically, if $(M_S, \mu_S) \in \mathcal{C}_G$ and there exists $(M_{S'}, \mu_{S'}) \in \mathcal{C}_G$ satisfying $(M_{S'}, \mu_{S'}) \geq (M_S, \mu_S)$ with $S'=S \cup \{\alpha \}$, then $(M_S, \mu_S)$ is strictly decreasing relative to $M_{S'}$.
\end{proposition}
\begin{proof}
We begin by proving the first statement. Uniqueness follows from Lemma \ref{uniqueness}. For existence, we first reduce to the case where  $M_S$ is a maximal Levi subgroup of $M_{S'}$. Suppose we have proven the proposition in this reduced case. We might then try to prove the general case by iteratively applying the reduced case of the proposition to a chain of standard Levi subgroups $M_S=M_{S_0} \subset ... \subset M_{S_k}=M_{S'}$ such that each is maximal in the next. Such a chain clearly exists, but to apply the reduced case of the proposition we need to show that if we have constructed a cocharacter pair $(M_{S_i}, \mu_{S_i}) \geq (M_S, \mu_S)$ then $(M_{S_i}, \mu_{S_i})$ is strictly decreasing relative to $M_{S'}$. This follows from Lemma \ref{2.22}.

Now, we let $\mu_{S'}$ be the unique  conjugate of $\mu_S$ which is dominant in $M_{S'}$. If we can show that $\theta_{M_{S'}}( \mu_{S'}) \prec \theta_{M_S}(\mu_S)$, then $(M_{S'}, \mu_{S'})$ will satisfy the conditions of the proposition. By Corollary \ref{conjinv} and Equation \eqref{kapavg},
\[
\theta_{M_{S'}}( \mu_{S'})=\theta_{M_{S'}}(\mu_S)=\frac{1}{|W_{M_{S'}}|} \sum\limits_{\sigma \in W_{M_{S'}}} \sigma( \theta_{M_S}(\mu_S)),
\] 
so we can reduce to showing that 
\[
\frac{1}{|W_{M_{S'}}|} \sum\limits_{\sigma \in W_{M_{S'}}} \sigma( y) \prec y,
\]
for any $y$ satisfying  $\langle y , \alpha \rangle > 0$ for $\alpha \in S' \setminus S$ and $\langle y, \alpha \rangle =0$ for $\alpha \in S$. Any such $y$ is dominant in the root datum of $M_{S'}$ and so by Lemma \ref{2.9}, 
\[
\frac{1}{|W_{M_{S'}}|} \sum\limits_{\sigma \in W_{M_{S'}}} \sigma( y) \preceq y.
\]
Further, the above equation cannot be an equality because  $y$ has positive pairing with each root of $S' \setminus S$ while $\frac{1}{|W_{M_{S'}}|} \sum\limits_{\sigma \in W_{M_{S'}}} \sigma( y)$ has $0$ pairing with these roots.

To prove the converse, suppose that $(M_S, \mu_S) \leq (M_{S'}, \mu_{S'})$ and $S'=S \cup \{\alpha \}$ for some $\alpha \in \Delta \setminus S$. Then by Corollary \ref{conjinv}
\[
\theta_{M_{S'}}(\mu_{S'}) = \theta_{M_{S'}}(\mu_{S})= \frac{\theta_{M_S}(\mu_S)+\sigma_{\alpha}( \theta_{M_S}(\mu_S))}{2},
\]
and so
\[
\theta_{M_S}(\mu_S)- \theta_{M_{S'}}(\mu_{S'})=\frac{\theta_{M_S}(\mu_S) - \sigma_{\alpha}(\theta_{M_S}(\mu_S))}{2}=\frac{1}{2} \langle \theta_{M_S}(\mu_S), \alpha \rangle \check{\alpha}.
\]
Since by assumption $\theta_{M_{S'}}(\mu_{S'}) \prec \theta_{M_S}(\mu_S)$, it follows that $\langle \theta_{M_S}(\mu_S), \alpha \rangle > 0$.
\end{proof}
\begin{remark}
Note that the converse of the above proposition is not true in the general case.
\end{remark}
\begin{corollary}{\label{sqr}}
Fix a standard Levi subgroup $M_S$ and roots $\alpha_1, \alpha_2 \in \Delta \setminus S$. Suppose we have cocharacter pairs $(M_S, \mu_S), (M_{S \cup \{\alpha_1\}}, \mu_{S \cup \{\alpha_1\}}), (M_{S \cup \{\alpha_1, \alpha_2\}}, \mu_{S \cup \{\alpha_1, \alpha_2\}}) \in \mathcal{C}_G$ satisfying 
\[
(M_S, \mu_S) \leq (M_{S \cup \{\alpha_1\}}, \mu_{S \cup \{\alpha_1\}}) \leq (M_{S \cup \{\alpha_1, \alpha_2\}}, \mu_{S \cup \{\alpha_1, \alpha_2\}})
\]
and that $(M_S, \mu_S)$ is strictly decreasing relative to $M_{S \cup \{\alpha_2\}}$. 

Then the extension of $(M_S, \mu_S)$ to $M_{S \cup \{\alpha_2\}}$, which we denote $(M_{S \cup \{\alpha_2\}}, \mu_{S \cup \{ \alpha_2 \}})$, satisfies
\[
(M_S, \mu_S) \leq (M_{S \cup \{\alpha_2\}}, \mu_{S \cup \{\alpha_2\}}) \leq (M_{S \cup \{\alpha_1, \alpha_2\}}, \mu_{S \cup \{\alpha_1, \alpha_2\}})
\]
\end{corollary}
\begin{proof}
By the second statement of Proposition \ref{ext}, we have that $(M_S, \mu_S)$ is strictly decreasing relative to $M_{S \cup \{\alpha_1\}}$. Then by Lemma \ref{2.22}, $(M_{S \cup \{\alpha_2 \}}, \mu_{S \cup \{ \alpha_2 \}})$ is strictly decreasing relative to $M_{S \cup \{\alpha_1, \alpha_2\}}$. Thus by Proposition \ref{ext}, we have $(M_{S \cup \{\alpha_2\}}, \mu_{S \cup \{\alpha_2\}}) \leq (M_{S \cup \{\alpha_1, \alpha_2\}}, \mu_{S \cup \{\alpha_1, \alpha_2\}})$ as desired.
\end{proof}
\begin{proposition} \label{2.26}
Let $S \subset S_1 \subset S_2$ be subsets of $\Delta$ and suppose $(M_S, \mu_S), (M_{S_2}, \mu_{S_2}) \in \mathcal{C}_G$ with
\[
(M_S, \mu_S) \leq (M_{S_2}, \mu_{S_2})
\]
and $(M_S, \mu_S)$ is strictly decreasing relative to  $M_{S_1}$. Then the unique extension $(M_{S_1}, \mu_{S_1})$ of $(M_S, \mu_S)$ to $M_{S_1}$ satisfies 
\[
(M_{S_1}, \mu_{S_1}) \leq (M_{S_2}, \mu_{S_2}).
\]
\end{proposition}
\begin{proof}
Since $(M_S, \mu_S) \leq (M_{S_2}, \mu_{S_2})$, there is an increasing chain of cocharacter pairs $(M_S, \mu_S)=(M_{S^0}, \mu_{S^0}) \leq ... \leq (M_{S^k}, \mu_{S^k})=(M_{S_2}, \mu_{S_2})$ such that each standard Levi subgroup is maximal in the next. The content of this proposition is that we can pick a chain such that $(M_{S_1}, \mu_{S_1})$ appears. By Lemma \ref{2.22}, we can assume that $M_S$ is maximal in $M_{S_1}$. Let $\alpha$ be the unique element of $S_1 \setminus S$.

Pick a chain of cocharacter pairs $(M_S, \mu_S)=(M_{S^0}, \mu_{S^0}) \leq ... \leq (M_{S^k}, \mu_{S^k})=(M_{S_2}, \mu_{S_2})$ as above. Chains of cocharacter pairs are determined by an ordering on the roots in $S_2 \setminus S=\{\alpha_1, ..., \alpha_k\}$, such that the $S^i=S \cup \{\alpha_1, ..., \alpha_i \}$. The root $\alpha$ appears in this chain so $\alpha=\alpha_i$ for some $i$. If $i=1$ we are done. Otherwise, we consider $(M_{S^{i-2}}, \mu_{S^{i-2}}) \leq (M_{S^{i-1}}, \mu_{S^{i-1}}) \leq (M_{S^i}, \mu_{S^i})$. By Lemma \ref{2.22}, $(M_{S^{i-2}}, \mu_{S^{i-2}})$ is strictly decreasing relative to $M_{S^{i-2} \cup \{ \alpha \}}$ and so by Corollary \ref{sqr} (applied so that $(M_{S^{i-2}}, \mu_{S^{i-2}})$ takes the place of $(M_S, \mu_S)$ in Corollary \ref{sqr}), we get a new chain of cocharacter pairs between $(M_S, \mu_S)$ and $(M_{S_2}, \mu_{S_2})$ where we switch the positions of $\alpha, \alpha_{i-1}$ in the corresponding ordering of $S_2 \setminus S$. By repeating this argument, we can construct a chain where $\alpha=\alpha_1$, which is what we need.
\end{proof}
The preceding propositions give us the following picture. Given a cocharacter pair $(M_S, \mu_S)$ we check which simple roots $\alpha$ satisfy $\langle \theta_{M_S}(\mu_S) , \alpha \rangle > 0$. Suppose there are $n$ such simple roots. Then we get $2^n$ standard Levi subgroups containing $M_S$ corresponding to adding different subsets of these simple roots. The cocharacter pair $(M_S, \mu_S)$ has a unique extension to each of the Levi subgroups and the poset lattice of these co-character pairs can be thought of as the graph of an $n$ dimensional cube in the following way.  The vertices of the cube are the $2^n$ cocharacter pairs extending $(M_S, \mu_S)$ that we have just constructed. For two such pairs $(M_{S_1}, \mu_{S_1}), (M_{S_2}, \mu_{S_2})$, we draw an edge between the two corresponding vertices if either $S_1 \subset S_2$ and $|S_2 \setminus S_1|=1$, or $S_2 \subset S_1$ and $|S_1 \setminus S_2|=1$. We can upgrade this graph to a directed graph by stipulating that an edge between $(M_{S_1}, \mu_{S_1})$ and $(M_{S_2}, \mu_{S_2})$ is directed from $(M_{S_1}, \mu_{S_1})$ to $(M_{S_2}, \mu_{S_2})$ if $(M_{S_2},\mu_{S_2}) < (M_{S_1}, \mu_{S_1})$. 

Finally, note that for any two pairs $(M_{S_1}, \mu_{S_1})$ and $(M_{S_2}, \mu_{S_2})$ corresponding to vertices in the above cube, we have $(M_{S_2}, \mu_{S_2}) \leq (M_{S_1}, \mu_{S_1})$ if and only if there is a directed path in the cube travelling from the vertex of $(M_{S_1}, \mu_{S_1})$ to that of $(M_{S_2}, \mu_{S_2})$.

\subsection{Connection With Isocrystals}{\label{CWI}}

We now investigate the relation between strictly decreasing cocharacter pairs and Kottwitz's theory of isocrystals with additional structure. See \cite{Kot1} for omitted details on the theory of isocrystals.

An isocrystal is a pair $(V, \Phi)$ where $V$ is a finite dimensional $\widehat{\mathbb{Q}_p^{ur}}$ vector space and $\Phi: V \to V$ is an additive transformation satisfying $\Phi(av)=\sigma(a) \Phi(v)$ for $a \in \widehat{\mathbb{Q}_p^{ur}}, v \in V$ and $\sigma$ the arithmetic Frobenius morphism. As before, let $G$ be a connected quasisplit reductive group defined over $\mathbb{Q}_p$ and consider the set of isomorphism classes of exact $\otimes$-functors from $\mathrm{Rep}(G)$ to $\mathrm{Isoc}$, the category of isocrystals. Such isomorphism classes are classified by $H^1(W_{\mathbb{Q}_p}, G(\overline{\widehat{\mathbb{Q}_p^{ur}}}))$ which we denote $\mathbf{B}(G)$ (where $W_{\mathbb{Q}_p}$ is the Weil group of $\mathbb{Q}_p$). 

In $\S 4.2$ of \cite{Kot1}, Kottwitz constructs the Newton map $\nu: \mathbf{B}(G) \to \overline{C}_{\mathbb{Q}}$ and the Kottwitz map $\kappa: \mathbf{B}(G) \to X^*(Z(\widehat{G})^{\Gamma})$. An element of $\mathbf{B}(G)$ is uniquely determined by its image under these maps.

We say that the standard Levi subgroup $M_S$ is associated to $b \in \mathbf{B}(G)$ if $\nu(b) \in \mathfrak{A}^+_{M_S, \mathbb{Q}}$. Henceforth, we will often denote the standard Levi subgroup associated to $b$ by $M_b$. Notice that many elements of $\mathbf{B}(G)$ could be associated to the same Levi subgroup. We call $b$ \emph{basic} if $M_b=G$. We write 
\[
\mathbf{B}(G)= \coprod_{S \subset \Delta} \mathbf{B}(G)_{M_S}
\]
such that $\mathbf{B}(G)_{M_S}$ consists of those $b \in \mathbf{B}(G)$ associated to $M_S$. We denote by $\mathbf{B}(M_S)^+$ the maximal subset of $\mathbf{B}(M_S)$ such that $\nu(\mathbf{B}(M_S)^+) \subset \overline{C}_{\mathbb{Q}}$. In $\S 5.1$ of \cite{Kot1}, Kottwitz uses the Kottwitz map for $M_S$ to construct canonical bijections 
\begin{align}{\label{BGiso}}
\mathbf{B}(G)_{M_S} \cong \mathbf{B}(M_S)^+_{M_S} \cong X^*(Z(\widehat{M_S})^{\Gamma})^+
\end{align}
where Kottwitz constructs a canonical isomorphism 
\begin{align}{\label{kotiso}}
X^*(Z(\widehat{M_S})^{\Gamma})_{\mathbb{Q}} \cong \mathfrak{A}_{M_S, \mathbb{Q}}
\end{align}
and $X^*(Z(\widehat{M_S})^{\Gamma})^+$ denotes the subset of $X^*(Z(\widehat{M_S})^{\Gamma})$ mapping to $\mathfrak{A}_{M_S, \mathbb{Q}}^+$. In fact, Kottwitz shows that the composition of the above isomorphisms gives the Newton map
\[
\mathbf{B}(G)_{M_S} \to \mathfrak{A}^+_{M_S, \mathbb{Q}} \hookrightarrow \overline{C}_{\mathbb{Q}}.
\]
For a further discussion of Equation \eqref{kotiso}, we refer the reader to Lemma \ref{lemma 2.12}.

We now prove an important lemma that will be used to relate the set $\mathbf{B}(G)$ to the strictly decreasing elements of $\mathcal{C}_G$.
\begin{lemma}{\label{kotimg}}
Fix a standard Levi subgroup $M_S$ of $G$ and let $(M_S, \mu_S ) \in \mathcal{SD}$.  Then $\theta_{M_S}(\mu_S) \in \nu(\mathbf{B}(G)_{M_S})$.
\end{lemma}
\begin{proof}
We first describe the set $\nu(\mathbf{B}(G)_{M_S})$. By Equations \eqref{BGiso} and \eqref{kotiso}, the set $\nu(\mathbf{B}(G)_{M_S})$ is equal to the image of $X^*(Z(\widehat{M_S})^{\Gamma})^+$ in $\mathfrak{A}_{M_S, \mathbb{Q}}$. Thus, to prove this lemma, it suffices to show that $\theta_{M_S}$ factors through the map $X^*(Z(\widehat{M_S})^{\Gamma}) \hookrightarrow X^*(Z(\widehat{M_S})^{\Gamma})_{\mathbb{Q}} \cong \mathfrak{A}_{M_S, \mathbb{Q}}$ where the isomorphism is as in Equation \eqref{kotiso} or Lemma \ref{lemma 2.12}. Then, since $(M_S, \mu_S)$ is strictly decreasing, the factoring of $\theta_{M_S}$ will map $\mu_S$ to an element of $X^*(Z(\widehat{M_S})^{\Gamma})^+$ as desired. That $\theta_{M_S}$ factors in this way follows from the alternate characterization of $\theta_{M_S}$ given in Proposition \ref{2.14}.
\end{proof}

\begin{definition}
Fix $\mu \in X_*(T)$. Then we recall the following definition of Kottwitz \cite[\S 6.2]{Kot1}:
\[
\mathbf{B}(G, \mu) :=\{b \in \mathbf{B}(G): \nu(b) \preceq \theta_T(\mu), \kappa(b)=\mu|_{Z(\widehat{G})^{\Gamma}}\}.
\]
\end{definition}
Now we prove the key result of this section, which permits us to associate an element of $\mathbf{B}(G)$ to each strictly decreasing cocharacter pair.
\begin{proposition} \label{2.16}
 We have a natural map 
\[
\mathcal{T}: \mathcal{SD} \to \mathbf{B}(G)
\]
defined as follows. Let $(M_S, \mu_S) \in \mathcal{SD}$. Then there exists a $b \in \mathbf{B}(G)$ so that $\kappa(b)= \mu_S |_{Z(\hat{G})^{\Gamma}}$ and $\nu(b)=\theta_{M_S}(\mu_S)$. We note that by construction, $b$ is unique. Then we define $\mathcal{T}((M_S, \mu_S))=b$. Furthermore, we show that
\[
\mathcal{T}(\mathcal{SD}_{\mu}) \subset \mathbf{B}(G, \mu).
\]

\end{proposition}

\begin{proof}
We first define $b$. Note that since $(M_S , \mu_S)$ is strictly decreasing, $\theta_{M_S}(\mu_S) \in \mathfrak{A}^+_{M_S, \mathbb{Q}}$. By Proposition \ref{2.14}, it follows that $\mu_S|_{Z(\widehat{M_S})^{\Gamma}} \in X^*(Z(\widehat{M_S})^{\Gamma})^+$ and so we can define $b$ to be the element of $\mathbf{B}(G)$ corresponding to $\mu_S|_{Z(\widehat{M_S})^{\Gamma}}$ under the isomorphism $B(G)_{M_S} \cong X^*(Z(\widehat{M_S})^{\Gamma})^+$ of Equation \eqref{BGiso}. Recall that the composition of this isomorphism with Equation \eqref{kotiso} induces the Newton map restricted to $\mathbf{B}(G)_{M_S}$. Thus, we have $\theta_{M_S}(\mu_S)=\nu(b)$. Equation (4.9.2) of \cite{Kot1}  implies that  $\kappa(b)= \mu_S |_{Z(\hat{G})^{\Gamma}}$.

It remains to show that if $(M_S, \mu_S) \in \mathcal{SD}_{\mu}$ then the element $b \in \mathbf{B}(G)$ that we have constructed lies in the set $\mathbf{B}(G, \mu)$. For this, we need to show that $\nu(b)= \theta_{M_S}(\mu_S) \preceq \theta_T(\mu)$.  

We claim that $\theta_T(\mu) \succeq \theta_T(\mu_S)$. After all, by (\cite[Ch6 1.6.18, p. 158]{Bou1}), we have  $\mu \succeq \mu_S$. Then the claim follows from Corollary \ref{posres}.

Now we claim that $\theta_T(\mu_S)$ is dominant in the relative root system of $M_S$. To prove the claim, we first observe that $\mu_S$ is dominant relative to the absolute root system of $M_S$. As above, the Galois group $\Gamma$ preserves the Weyl chamber corresponding to the positive absolute roots given by $B$. Thus, $\gamma( \mu_S)$ is dominant for each $\gamma \in \Gamma$, and so $\theta_T(\mu_S)$ is dominant relative to the absolute roots of $M_S$. The intersection of the closed positive Weyl chamber for the absolute root datum of $M_S$ with $\mathfrak{A}_{\mathbb{Q}}$ is the Weyl chamber for relative root datum of $M_S$ (cf. proof of Lemma \ref{2.11} $(2)$ ). Thus, $\theta_T(\mu_S)$ is dominant with respect to the relative roots as desired.

Finally, we apply Lemma \ref{2.9} and Equation \eqref{kapavg} to get
\[
\theta_T(\mu_S) \succeq \theta_{M_S}(\mu_S),
\]
which finishes the proof.
\end{proof}
\begin{question}
Can one describe the image
\[
\mathcal{T}(\mathcal{SD}_{\mu}) \subset \mathbf{B}(G, \mu)?
\]
\end{question}
Fix $G=GL_n$ with $T$ and $B$ the diagonal maximal torus and upper triangular Borel subgroup respectively. Suppose $\mu$ has weights $1$ and $0$. Then we claim $\mathcal{T}(\mathcal{SD}_{\mu}) = \mathbf{B}(G, \mu)$. Indeed, pick any $b \in \mathbf{B}(GL_n, \mu)$. Then without loss of generality, $\nu_b = ((a_1/b_1)^{x_1b_1}, ..., (a_r/b_r)^{x_rb_r})$ for some $a_i, b_i \in \mathbb{N}$ such that $a_i/b_i$ is written in reduced form. Then let $M$ be the standard Levi subgroup isomorphic to $GL_{x_1b_1} \times ... \times GL_{x_rb_r}$ and embedded diagonally. Since $b \in \mathbf{B}(GL_n, \mu)$, we must have that $\mu=(1^{\sum\limits^r_{i=1} x_ia_i}, 0^{n- \sum\limits^r_{i=1} x_ia_i})$. Finally, we define $\mu' \in X_*(T)$ by $\mu'=(1^{x_1a_1}, 0^{x_1b_1-x_1a_1}, ... , 1^{x_ra_r}, 0^{x_rb_r-x_ra_r})$. Then we note that $\mu'$ is dominant in the root system of $M$ so that $(M, \mu') \in \mathcal{C}_G$. Moreover, $\theta_M(\mu')=\nu_b$ so that $(M, \mu') \in \mathcal{SD}$. Then since $\mu'$ and $\mu$ are conjugate in $GL_n$, it is easy to see that $(M,\mu') \leq (GL_n, \mu)$. In conclusion, we have shown that $(M', \mu') \in \mathcal{SD}_{\mu}$ and $\mathcal{T}( (M', \mu'))=b$ as desired.

On the other hand for different choices of $\mu$, we can have $\mathcal{T}(\mathcal{SD}_{\mu}) \subsetneq \mathbf{B}(G, \mu)$. For instance, let $G=GL_3$, let $\mu=(2,0,0)$, and let $b \in \mathbf{B}(G,\mu)$ be such that $\nu_b=(1, 1/2, 1/2)$. Then it is easy to check that $\mathcal{T}(\mathcal{SD}_{\mu})$ does not contain $b$.
\subsection{The Induction and Sum Formulas}
We are now ready to prove our main theorems on cocharacter pairs. We begin by defining some key subsets of $\mathcal{C}_G$, the set of cocharacter pairs for $G$. In this section we fix a dominant $\mu \in X_*(T)$ and $b \in \mathbf{B}(G, \mu)$.

\begin{definition}
We define the sets $\mathcal{T}_{G,b,\mu}$ and $\mathcal{R}_{G,b,\mu}$ as follows:
\[
\mathcal{T}_{G,b,\mu} :=\mathcal{T}^{-1}(b) \cap \mathcal{SD}_{\mu}
\]
and
\[
\mathcal{R}_{G, b, \mu}=\{ (M_{S_1}, \mu_{S_1}) \in \mathcal{C}_G: (M_{S_1}, \mu_{S_1}) \leq (M_{S_2}, \mu_{S_2}) \,\ \text{for some} \,\ (M_{S_2}, \mu_{S_2}) \in \mathcal{T}_{G, b, \mu}\}.
\]
\end{definition}

\begin{definition}{\label{MGBmu}}
Let $\mathbb{Z} \langle \mathcal{C}_G \rangle$ denote the free Abelian group generated by the set of cocharacter pairs for $G$.

We define $\mathcal{M}_{G, b, \mu} \in \mathbb{Z} \langle \mathcal{C}_G \rangle$ by
\[
\mathcal{M}_{G, b, \mu}=\sum\limits_{(M_S, \mu_S) \in \mathcal{R}_{G, b, \mu}}(-1)^{L_{M_S, M_b}} (M_S, \mu_S)
\]
such that for $M_{S_1} \subset M_{S_2}$, $L_{M_{S_1}, M_{S_2}}$ is defined to be $| S_2 \setminus S_1 |$.
\end{definition}
\begin{remark}
We observe that for $(M_S, \mu_S) \in \mathcal{SD}$, if $\mathcal{T}((M_S, \mu_S))=b$, then $M_S=M_b$. 
\end{remark}

We will show in Theorem \ref{3.8} that at least in the case where $G$ is an unramified restriction of scalars of a general linear group, $\mathcal{M}_{G, b, \mu}$ is related to the cohomology of Rapoport-Zink spaces for $G$. Thus one expects there to be a combinatorial analogue of the Harris-Viehmann conjecture (Conjecture \ref{3.4}). We call this combinatorial analogue the \emph{induction formula}. Perhaps the more surprising result is that there is also an analogue of Shin's averaging formula (which we call the \emph{sum formula}) \cite[Thm 7.5]{Shi1}. We first prove the sum formula.

\begin{theorem}[Sum Formula]\label{combinatorial sum}

The following holds in $\mathbb{Z} \langle \mathcal{C}_G \rangle$: 
\[
\sum\limits_{b \in B(G, \mu)} \mathcal{M}_{G, b, \mu} = (G, \mu).
\]
\end{theorem}
\begin{proof}
We need to show that 
\[
\sum\limits_{b \in B(G, \mu)} \mathcal{M}_{G, b, \mu} = (G, \mu),
\]
or equivalently
\[
\sum\limits_{b \in B(G, \mu)} \sum\limits_{(M_S, \mu_S) \in \mathcal{R}_{G, b, \mu}} (-1)^{L_{M_S, M_b}} (M_S, \mu_S)=(G, \mu).
\]
We prove this equality by counting how many times a given cocharacter pair shows up on the left-hand side. The pair $(G, \mu)$ shows up exactly once in the left-hand sum as an element of $\mathcal{R}_{G,b,\mu}$ for $b$ the unique basic element of  $\mathbf{B}(G, \mu)$. Suppose$(M_S, \mu_S) \in \mathcal{C}_G$ is some other cocharacter pair. Then define 
\[
Y_{(M_S, \mu_S)}=\{b \in \mathbf{B}(G, \mu): (M_S, \mu_S) \in \mathcal{R}_{G, b, \mu} \}.
\]
 We are reduced to showing 
\begin{align}{\label{keyeqn}}
\sum\limits_{b \in Y_{(M_S, \mu_S)}} (-1)^{L_{M_S, M_b}} = 0.
\end{align}
Our general strategy will be to show that the left-hand side of equation \ref{keyeqn} vanishes for each $(M_S, \mu_S) < (G, \mu)$ by inducting on the size of $\Delta \setminus S$. However, in the case that $(M_S, \mu_S) \in \mathcal{SD}_{\mu}$, we can prove the vanishing without an inductive argument. We show this first before discussing the induction. 

Suppose now that $(M_S, \mu_S) \in \mathcal{SD}_{\mu}$. By Corollary \ref{sdupclsd}, every pair  $(M_{S'}, \mu_{S'}) \in \mathcal{C}_G$ satisfying $(M_S, \mu_S) \leq (M_{S'}, \mu_{S'}) \leq (G, \mu)$ is strictly decreasing and thus by Proposition \ref{2.16}, we have $\mathcal{T}((M_{S'}, \mu_{S'})) \in \mathbf{B}(G, \mu)$. These are precisely the elements $b \in \mathbf{B}(G, \mu)$ so that $(M_S, \mu_S) \in \mathcal{R}_{G, b, \mu}$. By the discussion after Proposition \ref{2.26}, we can associate the graph of a cube to the set of $(M_{S'}, \mu_{S'})$ such that each cocharacter pair is a vertex. To the vertex associated to $(M_{S'}, \mu_{S'})$ we attach the sign $(-1)^{L_{M_S, M_S'}}$. We note that adjacent vertices in this graph will have opposite signs since if $(M_{S'}, \mu_{S'})$ and $(M_{S''}, \mu_{S''})$ have adjacent vertices, then the cardinality of $S'$ and $S''$ differs by $1$. Now, it is a standard fact that if we associate an element of $\{ 1, -1 \}$ to each vertex of the graph of an $n$-dimensional cube for $n \geq 1$ so that adjacent vertices have opposite signs, then the sum of all the signs is $0$. This implies that the left-hand side of Equation \eqref{keyeqn} vanishes in the strictly decreasing case.

Now we discuss the inductive argument. The base case will be for pairs $(M_S, \mu_S)<(G, \mu)$ satisfying $|\Delta \setminus S|=1$. The second statement of Proposition \ref{ext} implies that in this case $(M_S, \mu_S)$ is strictly decreasing relative to $G$, which means that $(M_S, \mu_S) \in \mathcal{SD}_{\mu}$. Thus, the base case is proven by the previous paragraph.

We now discuss the inductive step. Suppose $(M_S, \mu_S)< (G, \mu)$. If $(M_S, \mu_S)$ is strictly decreasing, then we are done by the above. Suppose now that $(M_S, \mu_S)$ is not strictly decreasing. We claim that $(M_S, \mu_S)$ must be strictly decreasing with respect to at least some standard Levi subgroup of $G$ that properly contains $M_S$.  After all, since $(M_S, \mu_S) < (G, \mu)$, there must exist at least some $\alpha \in \Delta \setminus S$ and $(M_{S \cup \{\alpha\}}, \mu_{S \cup \{\alpha\}}) \in \mathcal{C}_G$ so that $(M_S, \mu_S) \leq (M_{S \cup \{\alpha\}}, \mu_{S \cup \{\alpha\}})$. Then by Proposition \ref{ext}, this implies that $(M_S, \mu_S)$ is strictly decreasing relative to $M_{S \cup \{\alpha\}}$.  

Thus, let $M_{S'}$ be the maximal standard Levi subgroup of $G$ such that $(M_S, \mu_S)$ is strictly decreasing relative to $M_{S'}$. We can write $S'=S \cup \{ \alpha_1, ..., \alpha_n \}$ where $\alpha_i \neq \alpha_j$ for $i \neq j$ and each $\alpha_i \in \Delta \setminus S$.  We denote by $X$ the $n$-cube of cocharacter pairs above $(M_S, \mu_S)$ as in the discussion after Proposition \ref{2.26}.

We claim that
\[
\sum\limits_{b \in Y_{(M_S, \mu_S)}} (-1)^{L_{M_S, M_b}} 
\]
\[
= -\sum\limits_{(M_{S'}, \mu_{S'}) \in X \setminus \{(M_S, \mu_S)\}} \,\ \sum\limits_{b \in  Y_{(M_{S'},\mu_{S'})}} (-1)^{L_{M_{S'}, M_b}}.
\]
Given this claim, we see that to finish the proof, it suffices to show that the right-hand side is identically $0$. However, the right-hand side consists of a sum of a number of terms similar to the left-hand side but for pairs $(M_{S'}, \mu_{S'})$ in place of $(M_S, \mu_S)$. Note that each $S'$ is strictly larger than $S$ and thus we are done by induction.

We now prove the claim. Moving all the terms to one side, we need only show that 
\[
\sum\limits_{(M_{S'}, \mu_{S'}) \in X} \,\ \sum\limits_{b \in Y_{ (M_{S'},\mu_{S'})}} (-1)^{L_{ M_{S'}, M_b}}=0.
\]

Fix $b \in \mathbf{B}(G, \mu)$. Then it suffices to show the contribution from $b$ in the above formula vanishes. Thus, we must show
\begin{align}{\label{finaleqn}}
\sum\limits_{(M_{S'}, \mu_{S'}) \in X \cap \mathcal{R}_{G, b, \mu}} (-1)^{L_{M_{S'},M_b}} =0.   
\end{align}

We examine the structure of $X \cap \mathcal{R}_{G, b, \mu}$ when it is nonempty. If we can show that the cocharacter pairs in this set form a sub-cube of $X$ of positive dimension, then we will be done by the standard fact that if we place alternating signs on the vertices of a cube and add up all the signs we get $0$. 

Clearly, any $(M_{S'}, \mu_{S'}) \in X \cap \mathcal{R}_{G, b, \mu}$ must satisfy $M_S \subset M_{S'} \subset M_b$. The subset of $X$ satisfying this latter property forms a sub-cube of $X$ since its elements are indexed by subsets of $S_b \setminus S$, where $S_b$ is the subset of $\Delta$ corresponding to $M_b$ in the standard way (note that by Lemma \ref{uniqueness}, there is at most one element of $X \cap \mathcal{R}_{G, b, \mu}$ for each standard Levi $M_{S'}$). Moreover, this latter set cannot form a cube of dimension $0$ for then we would have $M_S=M_b$ and so $X \cap \mathcal{R}_{G,b,\mu}=\{ (M_S, \mu_S)\}$ which would imply that $(M_S ,\mu_S)$ is strictly decreasing contrary to assumption.

Thus to finish the proof, we need only show that every $(M_{S'}, \mu_{S'})$ such that 
\begin{enumerate}
\item $M_S \subset M_{S'} \subset M_b$,\\
\item $(M_S, \mu_S) \leq (M_{S'}, \mu_{S'})$, \\
\item $(M_S, \mu_S)$ is strictly decreasing relative to $M_{S'}$,\\
\end{enumerate}
satisfies $(M_{S'}, \mu_{S'}) \leq (M_b, \mu_b)$ for some $(M_b, \mu_b) \in \mathcal{T}_{G, b, \mu}$. Since we assumed that $X \cap \mathcal{R}_{G, b, \mu} \neq \emptyset$, then in fact there is an $(M_b, \mu_b) \in \mathcal{T}_{G, b, \mu}$ with $(M_S, \mu_S) \leq (M_b, \mu_b)$. Then the desired result follows from Proposition \ref{2.26}.
\end{proof}

We now turn to the induction formula.
 Fix a standard Levi subgroup $M_S$ of $G$. Then our choice of maximal torus $T$ and Borel subgroup $B$ of $G$ provides us with natural choices $B \cap M_S$ and $T$ of a Borel subgroup and maximal torus of $M_S$. This allows us to define the set $\mathcal{C}_{M_S}$ of cocharacter pairs for $M_S$. There is a natural inclusion
 \begin{align}{\label{inclusion}}
     i^G_{M_S}:  \mathcal{C}_{M_S} \hookrightarrow \mathcal{C}_G.
 \end{align}
 The image of this inclusion is precisely the set of cocharacter pairs $(M_{S'}, \mu_{S'})$ where $S' \subset S$. This inclusion preserves the partial ordering of cocharacter pairs. The strictly decreasing elements of $\mathcal{C}_{M_S}$ map to the elements of $\mathcal{C}_G$ which are strictly decreasing relative to $M_S$.
 
Now choose a $b \in \mathbf{B}(G, \mu)$ and rational Levi $M_S$ such that $M_b \subset M_S \subset G$. We have a unique $b' \in \mathbf{B}(M_b)^+_{M_b}$ corresponding to $b$ under the isomorphism given by Equation \eqref{BGiso}. The inclusion $M_b \subset M_S$ induces a map 
\[
\mathbf{B}(M_b) \to \mathbf{B}(M_S).
\]
Let $b_S$ be the image of $b'$ under this map. 

The following definition will be important in relating cocharacter pairs of a group $G$ to those of a standard Levi. Compare with \cite[Equation (8.1)]{RV1}.
\begin{definition}{\label{I_M,b,mu defn}}
Let $M_S$ be a standard Levi subgroup of $G$, let $\mu \in X_*(T)$ be a dominant cocharacter and choose $b \in \mathbf{B}(G, \mu)$. We take $b_S \in \mathbf{B}(M_S)$ as constructed in the previous paragraph and define the set 

$\mathcal{I}^{G, \mu}_{M_S, b_S}=\{ (M_S, \mu_S) \in \mathcal{C}_{M_S} : b_S \in \mathbf{B}(M_S, \mu_S), \mu_S$ is conjugate to $\mu$ in $G \}$.

\end{definition}
We first check the following transitivity property of $\mathcal{I}^{G, \mu}_{M_S, b_S}$.
\begin{proposition}{\label{Itrans}}
Fix $(G, \mu) \in \mathcal{C}_G$ and $b \in \mathbf{B}(G, \mu)$. Suppose $M_{S_2}$ and $M_{S_1}$ are standard Levi subgroups of $G$ such that $M_b \subset M_{S_2} \subset M_{S_1}$. Then
\[
\mathcal{I}^{G, \mu}_{M_{S_2}, b_{S_2}} =\{ (M_{S_2}, \mu_{S_2}) \in \mathcal{C}_{M_{S_2}} : (M_{S_2}, \mu_{S_2}) \in \mathcal{I}^{M_{S_1},  \mu_{S_1}}_{M_{S_2}, b_{S_2}} \,\ \text{for some} \,\ (M_{S_1}, \mu_{S_1}) \in \mathcal{I}^{G, \mu}_{M_{S_1}, b_{S_1}} \}.
\]
\end{proposition}
\begin{proof}
We show each set is a subset of the other. Take $(M_{S_2}, \mu_{S_2}) \in \mathcal{I}^{G, \mu}_{M_{S_2}, b_{S_2}}$. Let $\mu_{S_1}$ be the unique dominant cocharacter conjugate to $\mu_{S_2}$ in $M_{S_1}$. Then we consider $(M_{S_1}, \mu_{S_1})$ as an element of $\mathcal{C}_{M_{S_1}}$ and just need to show that $b_{S_1} \in \mathbf{B}(M_{S_1}, \mu_{S_1})$ since we already know that $b_{S_2} \in \mathbf{B}(M_{S_2}, \mu_{S_2})$ by assumption. Thus, we need only show that $\nu(b_{S_1}) \leq \theta_T(\mu_{S_1})$ and $\kappa(b_{S_1})=\mu_{S_1}|_{Z(\widehat{M_{S_1}})^{\Gamma}}$. 

We prove the inequality first. By assumption, $\nu(b_{S_2}) \preceq \theta_T(\mu_{S_2})$ and by Equations \eqref{BGiso} and \eqref{kotiso}, $\nu(b_{S_1})=\nu(b)=\nu(b_{S_2})$. Since $\mu_{S_1}$ and $\mu_{S_2}$ are conjugate in $M_{S_1}$ and $\mu_{S_1}$ is dominant, it follows from \cite[Ch6 1.6.18, p. 158]{Bou1} that $\mu_{S_2} \preceq \mu_{S_1}$. Then, by Corollary \ref{posres} it follows that $\theta_T(\mu_{S_2}) \preceq \theta_T(\mu_{S_1})$ in the relative root system. Combining all this data, we get
\[
\nu(b_{S_1})=\nu(b_{S_2}) \preceq \theta_T(\mu_{S_2}) \preceq \theta_T(\mu_{S_1}),
\]
as desired.

To prove $\kappa(b_{S_1})=\mu_{S_1}|_{Z(\widehat{M_{S_1}})^{\Gamma}}$, we note that by Equation (4.9.2) of \cite{Kot1} and the fact that $b_{S_2} \in \mathbf{B}(M_{S_2}, \mu_{S_2})$, we have 
\[
\kappa(b_{S_1})= \mu_{S_2}|_{Z(\widehat{M_{S_1}})^{\Gamma}}.
\]
Then $\mu_{S_1}$ and $\mu_{S_2}$ are conjugate in $M_{S_1}$ so there exists a $w \in W^{\mathrm{abs}}_{M_{S_1}}$ so that $w(\mu_1)=\mu_2$. This implies that $\mu_1$ and $\mu_2$ are conjugate in $\widehat{M_{S_1}}$ and in particular equal when restricted to $Z(\widehat{M_{S_1}})$. This implies the desired equality.

To show the converse inclusion, we start with $(M_{S_2}, \mu_{S_2}) \in \mathcal{I}^{M_{S_1}, \mu_{S_1}}_{M_{S_2}, b_{S_2}}$ for some $(M_{S_1}, \mu_{S_1}) \in \mathcal{I}^{G, \mu}_{M_{S_1}, b_{S_1}}$ and need to show that $b_{S_2} \in \mathbf{B}(M_{S_2}, \mu_{S_2})$ and that $\mu_{S_2}$ is conjugate to $\mu$ in $G$. But $(M_{S_2}, \mu_{S_2}) \in \mathcal{I}^{M_{S_1}, \mu_{S_1}}_{M_{S_2}, b_{S_2}}$ implies that $b_{S_2} \in \mathbf{B}(M_{S_2}, \mu_{S_2})$ and also that $\mu_{S_2}$ is conjugate to $\mu_{S_1}$ in $M_{S_1}$. Further, $(M_{S_1}, \mu_{S_1}) \in \mathcal{I}^{G, \mu}_{M_{S_1}, b_{S_1}}$ implies that $\mu_{S_1}$ is conjugate to $\mu$ in $G$. Thus, $\mu_{S_2}$ is conjugate to $\mu$ in $G$ as desired.
\end{proof}

The set $\mathcal{I}^{G, \mu}_{M_S, b_S}$ will primarily be useful because it allows us to relate the set $\mathcal{T}_{G, b, \mu}$ to analogous constructions in $M_S$. This is encapsulated in the following proposition.
\begin{proposition}{\label{relTbu}}
Fix $M_S$, $\mu$ and $b$ as in Definition \ref{I_M,b,mu defn}. The natural inclusion $i^G_{M_S} : \mathcal{C}_{M_S} \hookrightarrow \mathcal{C}_G$ of Equation \eqref{inclusion} induces a bijection 
\[
\coprod_{(M_S, \mu_S) \in \mathcal{I}^{G, \mu}_{M_S, b_S}} \mathcal{T}_{M_S, b_S, \mu_S} \cong \mathcal{T}_{G, b, \mu}
\]
\end{proposition}
\begin{proof}
We first show that 
\[
i^G_{M_S}(\coprod_{(M_S, \mu_S) \in \mathcal{I}^{G, \mu}_{M_S, b_S}} \mathcal{T}_{M_S, b_S, \mu_S}) \supset \mathcal{T}_{G, b, \mu}.
\]
Since $M_b \subset M_S$, it follows from the discussion after Equation \eqref{inclusion} that \[
\mathcal{T}_{G, b, \mu} \subset i^G_{M_S}(\mathcal{C}_{M_S}).
\]
Thus, pick an arbitrary element of $\mathcal{T}_{G, b, \mu}$ of the form $i^G_{M_S}(M_b ,\mu_b)$ for $(M_b, \mu_b) \in \mathcal{C}_{M_S}$. The cocharacter pair $i^G_{M_S}(M_b, \mu_b)$ is strictly decreasing, and therefore so is $(M_b, \mu_b) \in \mathcal{C}_{M_S}$. By Proposition \ref{ext} we can find $(M_S, \mu_S) \in \mathcal{C}_{M_S}$ such that $(M_b, \mu_b) \leq (M_S, \mu_S)$. Observe that since $i^G_{M_S}(M_b, \mu_b) \leq (G, \mu)$, the cocharacter $\mu_b$ is conjugate to $\mu$ in $G$ and therefore $\mu_S$ must be as well by construction. If we can show that $\mathcal{T}((M_b, \mu_b))=b_S$, then we will be done because by Proposition \ref{2.16}, this implies that $b_S \in \mathbf{B}(M_S, \mu_S)$ and so therefore that $(M_S, \mu_S) \in \mathcal{I}^{G, \mu}_{M_S, b_S}$ and $(M_b, \mu_b) \in \mathcal{T}_{M_S, b_S, \mu_S}$.

By assumption,  $\mathcal{T}(i^G_{M_S}(M_b, \mu_b))=b \in \mathbf{B}(G, \mu)$. Recall that the map $\mathcal{T}$ is defined so that a strictly decreasing $(M_b, \mu_b) \in \mathcal{C}_G$ which satisfies $(M_b, \mu_b) \leq (G, \mu)$ is mapped first to the element $\mu_b|_{Z(\widehat{M_b})^{\Gamma}} \in X^*(Z(\widehat{M_b})^{\Gamma})^+$. Then, this element is identified  with an element of $\mathbf{B}(G)$ via the isomorphisms of Equation \eqref{BGiso}:
\[
X^*(Z(\widehat{M_b}))^{\Gamma})^+ \cong \mathbf{B}(M_b)^+_{M_b} \cong \mathbf{B}(G)_{M_b},
\]
where the second isomorphism above is induced by the inclusion $M_b \hookrightarrow G$. We have the commutative diagram
\[
\begin{tikzcd}
\mathbf{B}(M_b) \arrow[r] \arrow[rd] & \mathbf{B}(M_S) \arrow[d]  \\
& \mathbf{B}(G)
\end{tikzcd}\\
\]
where each map is induced from the inclusion of groups. By definition, the element $b' \in \mathbf{B}(M_b)^+$ maps to $b \in \mathbf{B}(G)$ and $b_S \in \mathbf{B}(M_S)$ respectively. Thus, we see that by construction, $\mathcal{T}((M_b, \mu_b)) =b_S$.

Conversely, suppose $(M_b, \mu_b) \in \mathcal{T}_{M_S, b_S, \mu_S}$ for some $(M_S, \mu_S) \in \mathcal{I}^{G, \mu}_{M_S, b_S}$. Since $b' \in \mathbf{B}(M_b)^+_{M_b}$, it follows from the definition of $b_S$ and $\mathcal{T}_{M_S, b_S, \mu_S}$ that $\mu_b |_{Z(\widehat{M_b})^{\Gamma}}$ is an element of $X^*(Z(\widehat{M_b})^{\Gamma})^+$. This implies that $i^{G}_{M_S} (M_b, \mu_b) \in \mathcal{SD}$. By Proposition \ref{ext}, we have an extension of $i^{G}_{M_S} (M_b, \mu_b)$ to $G$, and since $\mu_b$ and $\mu$ are conjugate in $G$ by assumption, it follows that this extension is $(G, \mu)$ so that $i^G_{M_S}(M_b, \mu_b) \leq (G, \mu)$. It follows from these facts that $i^G_{M_S}(M_b, \mu_b) \in \mathcal{T}_{G, b, \mu}$. 

Finally, we remark that for distinct $(M_S, \mu_S), (M_S, \mu'_S) \in \mathcal{I}^{G, \mu}_{M_S, b_S}$ the sets $\mathcal{T}_{M_S, b_S, \mu_S}$ and $\mathcal{T}_{M_S, b_S, \mu'_S}$ are indeed disjoint by Lemma \ref{uniqueness}.
\end{proof}
As a corollary of this result, we have the induction formula.
\begin{corollary}[Induction Formula]\label{combinatorial induction}
We continue using the notation of the previous proposition. The natural map
\[
i^G_{M_S}: \mathcal{C}_{M_S} \hookrightarrow \mathcal{C}_G,
\]
induces a map
\[
i^G_{M_S}: \mathbb{Z} \langle \mathcal{C}_{M_S} \rangle \hookrightarrow \mathbb{Z} \langle \mathcal{C}_G \rangle,
\]
which gives an equality
\[
\sum\limits_{(M_S, \mu_S) \in \mathcal{I}^{G, \mu}_{M_S, b_S}} i^G_{M_S}(\mathcal{M}_{M_S, b_S, \mu_S})=\mathcal{M}_{G, b, \mu}.
\]
\end{corollary}
\begin{proof}
It follows from Proposition \ref{relTbu} that the map $i^G_{M_S}$ induces a bijection
\[
\coprod_{(M_S, \mu_S) \in \mathcal{I}^{G, \mu}_{M_S, b_S}} \mathcal{R}_{M_S, b_S, \mu_S} \cong \mathcal{R}_{G, b, \mu}.
\]
We remark that for distinct $(M_S, \mu_S), (M_S, \mu'_S) \in \mathcal{I}^{G, \mu}_{M_S, b_S}$ we have $\mathcal{R}_{M_S, b_S, \mu_S} \cap \mathcal{R}_{M_S, b_S, \mu'_S}= \emptyset$ by Lemma \ref{uniqueness}.

The corollary then follows from the definition of $\mathcal{M}_{G, b, \mu}$.
\end{proof}
This result can be thought of as an analogue of the $\emph{Harris-Viehmann}$ conjecture which we discuss in the next section.

In the cases we are interested in, we will also need a description of how cocharacter pairs behave with respect to products.

Suppose $G=G_1 \times ... \times G_k$ and $T=T_1 \times ... \times T_k$ such that $T_i$ is a maximal torus for $G_i$. Then
\[
X_*(T) \cong X_*(T_1) \oplus ... \oplus X_*(T_k),
\]
and any standard Levi subgroup admits a product decomposition 
\[
M_S \cong M_{S_1} \times ... \times  M_{S_k},
\]
such that $T_i \subset M_{S_i} \subset G_i$. Then any cocharacter pair $(M_S, \mu_S)$ of $G$ corresponds to a tuple of cocharacter pairs 
\[
((M_{S_1}, \mu_{S_1}), ..., (M_{S_k}, \mu_{S_k})) \in \mathcal{C}_{G_1} \times ... \times \mathcal{C}_{G_k},
\]
in the obvious way. The pair $(M_S, \mu_S) \in \mathcal{C}_G$ is strictly decreasing if and only if each pair $(M_{S_i}, \mu_{S_i}) \in \mathcal{C}_{G_i}$ is, and if $\mathcal{T}((M_S, \mu_S))= b \in \mathbf{B}(G, \mu)$, then we also have $\mathcal{T}_i((M_{S_i}, \mu_{S_i}))=b_i \in \mathbf{B}(G_i, \mu_i)$ where $\mathcal{T}_i$ is the map $\mathcal{T}$ defined for the group $G_i$. Thus, $b \mapsto (b_1, ..., b_k)$ under the natural bijection 
\[
\mathbf{B}(G) \cong \mathbf{B}(G_1) \times ... \times \mathbf{B}(G_k).
\]
We record the following proposition
\begin{proposition}{\label{prodcochar}}
We use the notation of the previous two paragraphs. 

The natural bijection 
\[
\mathcal{C}_G \cong \mathcal{C}_{G_1} \times ... \times \mathcal{C}_{G_k},
\]
induces bijections
\[
\mathcal{T}_{G, b, \mu} \cong \mathcal{T}_{G_1, b_1, \mu_1} \times ... \times \mathcal{T}_{G_k, b_k, \mu_k},
\]
and
\[
\mathcal{R}_{G, b, \mu} \cong \mathcal{R}_{G_1, b_1, \mu_1} \times ... \times \mathcal{R}_{G_k, b_k, \mu_k}.
\]
Further, under the natural isomorphism $\mathbb{Z} \langle \mathcal{C}_G \rangle \cong \mathbb{Z} \langle \mathcal{C}_{G_1} \rangle \otimes  ... \otimes \mathbb{Z} \langle \mathcal{C}_{G_k} \rangle$ we have
\[
\mathcal{M}_{G,b,\mu}=\mathcal{M}_{G_1, b_1, \mu_1} \otimes ... \otimes \mathcal{M}_{G_k, b_k, \mu_k}.
\]
\end{proposition}

\section{Cohomology of Rapoport-Zink spaces and the Harris-Viehmann Conjecture}
In this section, we define the Rapoport-Zink spaces we will work with and show how we can describe their cohomology using the language developed in the previous section. We also give a statement of the Harris-Viehmann conjecture, and explain the necessity of a small correction to the conjecture. We follow \cite{Far1}, \cite{Shi1}, and \cite{RV1}.

The theory necessarily involves several choices of signs. This is often a point of confusion, so we describe our conventions here. We choose the cocharacter $\mu$ appearing in the definition of Rapoport-Zink spaces to have non-negative weights, in agreement with most authors. In this paper, we use the contravariant Dieudonne functor, which means that our $p$-divisible groups will have isocrystals in the set $\mathbf{B}(G, \mu)$ (as opposed to $\mathbf{B}(G, -\mu)$ for the covariant theory). This convention agrees with that of \cite{Far1} and \cite{RV1}, but \cite{Shi1} uses the opposite convention. We use the local Langlands correspondence for $\mathrm{GL}_n(\mathbb{Q}_p)$ as in \cite[pg. 2]{HT1}. In particular, we normalize the local Artin map so that uniformizers correspond to geometric Frobenius elements.

\subsection{Rapoport-Zink Spaces of EL-Type}

We fix the following notation. Suppose $G$ is a reductive group defined over a field $k$ and $\mu \in X_*(G)$. Then if $H$ is a subgroup of $G$ such that $\mu$ factors through the inclusion $X_*(H) \hookrightarrow X_*(G)$, we denote by $\{\mu\}_H$ the $H(\overline{k})$ conjugacy class of $\mu$ and by $E_{\{\mu\}_H}$ the field of definition of $\{\mu\}_H$ (i.e the smallest extension of $k$ so that each element of $\mathrm{Gal}(\overline{k}/E_{\{\mu\}_H})$ stabilizes $\{\mu\}_H$).

Now we define the Rapoport-Zink data we consider.
\begin{definition}{\label{RZdat}}
An \emph{unramified Rapoport-Zink datum of EL type} is a tuple\\ $(F, V, \{\mu \}_G, b)$ where

\begin{enumerate}
\item $F$ is a finite unramified extension of $\mathbb{Q}_p$,
\item $V$ is a finite dimensional $F$ vector space,
\item $G := \mathrm{Res}_{F/ \mathbb{Q}_p}(\mathrm{GL}_F(V))$,
\item $\mu: \mathbb{G}_{m, \overline{\mathbb{Q}_p}} \to G_{\overline{\mathbb{Q}_p}}$ is a cocharacter inducing a weight decomposition $V \otimes \widehat{\mathbb{Q}^{ur}_p} \cong V_0 \oplus V_1$ where $\mu(z)$ acts by $z^i$ on $V_i$,
\item $b \in \mathbf{B}(G, \mu)$.
\end{enumerate}
\end{definition}
We fix a Borel subgroup $B \subset G$ defined over $\mathbb{Q}_p$, a $\mathbb{Q}_p$-split torus $A \subset G$ of maximal rank in $G$ and such that $A \subset B$, and a maximal torus $T \subset B$ containing $A$ and defined over $\mathbb{Q}_p$. We can choose $\mu$ in the above definition so that it is dominant relative to $B$. 

Let $\mathbb{X}$ be a $p$-divisible group defined over $\overline{\mathbb{F}_p}$ with an action of $\mathcal{O}_F$ and such that the isocrystal attached to $\mathbb{X}$ by the contravariant Dieudonne functor is isomorphic to $(V_F, b \sigma)$. We consider the moduli functor $\mathbb{M}_{b, \mu}$ such that for $S$ a scheme over $\mathcal{O}_{\widehat{\mathbb{Q}^{ur}_p}}$ with $p$ locally nilpotent, $\mathbb{M}_{b, \mu}(S)=\{(X, i, \rho)\}/ \sim$. Where $X$ is a $p$-divisible group defined over $S$, $i: \mathcal{O}_F \to \mathrm{End}_F(X)$, and $\rho: \mathbb{X} \times_{\overline{\mathbb{F}_p}} \overline{S} \to \overline{X}$ is a quasi-isogeny ($\overline{S}, \overline{X}$ are the reductions modulo $p$).

By work of Rapoport and Zink \cite[Thm 3.25]{RZ1}, the above moduli problem is represented by a formal scheme over $\mathcal{O}_{\widehat{\mathbb{Q}^{ur}_p}}$ which we also denote by $\mathbb{M}_{b, \mu}$. We have the generic fiber $\mathbb{M}^{rig}_{b, \mu}$ which is a rigid analytic space over $\widehat{\mathbb{Q}^{ur}_p}$. Further, we get a tower of coverings $\mathbb{M}^{rig}_{b, \mu, U}$ of $\mathbb{M}^{rig}_{b, \mu}$ for each compact open subgroup  $U \subset G(\mathbb{Q}_p)$. 

For a fixed prime $l \neq p$, we denote by $ H^j_c(\mathbb{M}^{rig}_{b, \mu, U} \times \overline{\widehat{\mathbb{Q}^{ur}_p}}, \overline{\mathbb{Q}_l})$ the etale cohomology with compact supports. This is a $\overline{\mathbb{Q}_l}$ vector space which is a smooth representation of $J_b(\mathbb{Q}_p) \times W_{E_{\{\mu\}_G}}$, where $J_b$ is the inner form of the standard Levi subgroup $M_b$ associated to $b$ (as constructed in $\S 3.3 $ of \cite{Kot1}) and $W_{E_{\{\mu\}_G}}$ is the Weil group of $E_{\{\mu\}_G}$ (for example see \cite[Prop 6.1]{RV1}).

We use the notation $\mathrm{Groth}( \cdot )$ for the Grothendieck group of admissible representations of topological groups. See $\S I.2$ of \cite{HT1} for the precise definition of these Grothendieck groups.

Let $P_b$ be the standard parabolic subgroup with Levi factor $M_b$ and denote the opposite parabolic by $P^{op}_b$. 
We define $J^G_P, \mathrm{Jac}^G_P$ to be the normalized and un-normalized Jacquet module functors, and we define $I^G_P, \mathrm{Ind}^G_P$ to be the normalized and un-normalized parabolic induction functors. Often, if $M \subset P$ is the standard Levi subgroup of $P$ and we are taking $I^G_P$ or $I^G_{P^{op}}$ to be a map of Grothendieck groups, we will write $I^G_M$ to remind the reader that these maps do not depend on choice of $P, P^{op}$ when considered as maps of Grothedieck groups.

\label{Mant defn} In \cite{Man2}, Mantovan considers the following construction (see also \cite{Shi1}).  We define a map 
\[
\mathrm{Mant}_{G, b, \mu}: \mathrm{Groth}(J_b(\mathbb{Q}_p)) \to \mathrm{Groth}(G(\mathbb{Q}_p) \times W_{E_{\{\mu\}_G}}),
\]
by 
\[
\mathrm{Mant}_{G,b, \mu}(\rho)= \sum\limits_{i, j \geq 0} (-1)^{i+j} \varinjlim\limits_{U \subset G(\mathbb{Q}_p)} \mathrm{Ext}^i_{J_b(\mathbb{Q_p})}(H^j_c(\mathbb{M}^{rig}_{b, \mu, U} \times \overline{\hat{\mathbb{Q}^{ur}_p}}, \overline{\mathbb{Q}_l}), \rho)(- \mathrm{dim} \mathbb{M}^{rig}_{b, \mu, U}).
\]
In $\S 6.2$ of \cite{Shi1} and $\S 2.4$ of \cite{Shi2}, Shin considers a map 
\[
\mathrm{Red}_b: \mathrm{Groth}(G(\mathbb{Q}_p)) \to \mathrm{Groth}(J_b(\mathbb{Q}_p)).
\]
We follow the construction given in \cite{Shi2}\footnote{We believe the construction given before Lemma 6.2 of \cite{Shi1} has a slight typo. There, $\mathrm{Red}_b$ is defined by $\pi \mapsto e(J_b) ( \mathrm{LJ} \circ \mathrm{Jac}^G_{P^{op}_b}(\pi))$. As maps of Grothendieck groups, $\mathrm{Jac}^G_{P^{op}_b}=\mathrm{J}^G_{P^{op}_b} \otimes \delta^{\frac{1}{2}}_{P^{op}_b}=\mathrm{J}^G_{P^{op}_b} \otimes \delta^{-\frac{1}{2}}_{P_b}$. But this is not equal to $\mathrm{J}^G_{P^{op}_b}(\pi)\otimes \delta^{\frac{1}{2}}_{P_b}$, which is the construction given in \cite{Shi2} that is compatible with \cite{HT1}.}.  We define $\mathrm{Red}_b$ by
\[
\pi \mapsto e(J_b) ( \mathrm{LJ} \circ \mathrm{J}^G_{P^{op}_b}(\pi)\otimes \delta^{\frac{1}{2}}_{P_b}),
\]
where 
\[
LJ: \mathrm{Groth}(M_b(\mathbb{Q}_p)) \to \mathrm{Groth}(J_b(\mathbb{Q}_p)),
\]
is the map defined by Badulescu extending the inverse Jacquet-Langlands correspondence (see \cite[Prop 3.2]{Bad1})  and $e(J_b)$ is the Kottwitz sign as defined in \cite{Kot4}.

We now describe the main result of \cite{Shi1}. The cocharacter $\mu$ of $G$ is a map $\mu: \mathbb{G}_{m, \overline{\mathbb{Q}_p}} \to \prod\limits_{\tau \in \mathrm{Hom}(F, \overline{\mathbb{Q}_p})} GL_{n,  \overline{\mathbb{Q}_p}}$ such that the weights in each $\mathrm{GL}_n$ factor are $1$s or $0$s. Thus we let $p_{\tau}, q_{\tau}$ denote the number of $1$ and $0$ weights respectively in the factor corresponding to $\tau$. 

The following formula is the main theorem in \cite[Thm 7.5]{Shi1}.
\begin{theorem}[Shin]\label{3.2}
We have the following equality for accessible representations in  $\mathrm{Groth}(G(\mathbb{Q}_p) \times W_{E_{\{\mu\}_G}})$.

\[
\sum\limits_{b \in B(G, \mu)} \mathrm{Mant}_{b, \mu}( \mathrm{Red}_b( \pi)) = [ \pi ] [ r_{-\mu} \circ \mathrm{LL}( \pi)|_{W_{E_{\{\mu\}_G}}} \otimes | \cdot |^{\mathrm{- \sum_{\tau} p_{\tau} q_{\tau}/2}}].
\]
\end{theorem}
Loosely speaking, accessible representations in Shin's paper are character twists of the local components of global representations that can be found within the cohomology of Shimura varieties. Shin shows that all essentially square-integrable representations are accessible.

In this case $\mathrm{LL}$ is the semisimplified local Langlands correspondence (known by the work of \cite{HT1} for instance). The map $r_{-\mu}$ is the algebraic representation of $\widehat{G} \rtimes W_{E_{\{\mu\}_G}} \subset {}^LG$ defined by Kottwitz (\cite[Lem 2.1.2]{Kot3}). It is characterized by the fact that $r_{-\mu}|_{\widehat{G}}$ is the irreducible representation of extreme weight $-\mu$ and if we take a $\Gamma$-invariant splitting of $\widehat{G}$, then the subgroup $W_{E_{\{\mu\}_G}}$ of $^LG$ acts trivially on the highest weight vector of $r_{-\mu}$ associated with this splitting.

\begin{remark}
The Tate twist appearing on the right-hand side of the above formula comes from the dimension formula for Shimura varieties and is equal to $-\langle \rho_G, \mu \rangle$ where $\rho_G$ is the half sum of the positive roots in $G$. 
\end{remark}

The above theorem is analogous to the sum formula for cocharacter pairs (Theorem \ref{combinatorial sum}). The induction formula (Corollary \ref{combinatorial induction}) is related to the Harris-Viehmann conjecture (Conjecture \ref{3.4} in this document). A proof of this conjecture is expected to appear in forthcoming work of Scholze.

\subsection{Harris-Viehmann Conjecture}{\label{HV}}

We now state the Harris-Viehmann conjecture following Rappoport and Viehmann in \cite{RV1}. In this subsection, we return to the notation of \S \ref{2} so that in particular, $G$ is a connected, quasisplit reductive group defined over $\mathbb{Q}_p$. 

Choose a dominant minuscule $\mu \in X_*(T)$ (where we can consider $\mu$ as a cocharacter of $G$ since $T \subset G$) and a $b \in \mathbf{B}(G, \mu)$ . Associated to $b$, we have the standard Levi subgroup $M_b$. Suppose we have a standard rational Levi subgroup $M_S$ so that $M_b \subset M_S \subset G$. We define $b', b_S$ as we did before Definition \ref{I_M,b,mu defn}.

In \cite[Equation (6.2)]{RV1}, the authors associate a cohomological construction to the triple $(G, b, \mu)$ which they denote $H^{\bullet}((G, [b], \{ \mu \}))$. This construction is a map of Grothendieck groups: $H^{\bullet}((G, [b], \{ \mu \})): \mathrm{Groth}(J_b(\mathbb{Q}_p)) \to \mathrm{Groth}(G(\mathbb{Q}_p) \times W_{E_{\{\mu\}}})$ and agrees with $\mathrm{Mant_{G, b, \mu}}$ in the case above. We will denote this construction $H^{\bullet}(G, b, \mu)$ since we deal with dominant cocharacters instead of conjugacy classes. Then they have the following conjecture.

\begin{conjecture}[Harris-Viehmann]\label{3.4}
For $\rho \in \mathrm{Groth}(J_b(\mathbb{Q}_p))$, we have the equality
\[
H^{\bullet}(G, b, \mu )[\rho]=\sum\limits_{ (M_S, \mu_S) \in \mathcal{I}^{G, \mu}_{M_S, b_S}} (\mathrm{Ind}^G_{P_S} H^{\bullet}(M_S, b_S, \mu_S)[\rho])\otimes [1][|\cdot|^{\langle \rho_G, \mu_S \rangle - \langle \rho_G, \mu \rangle}],
\]
in $\mathrm{Groth}(G(\mathbb{Q}_p) \times W_{E_{\{\mu\}_G}})$. The parabolic induction only modifies the $\mathrm{Groth}(G(\mathbb{Q}_p))$ parts of these representations. 
\end{conjecture}

\begin{remark}{\label{Harris-Viehmann}}
 We need to explain several things in the above conjecture. First we explain why the right-hand side is a representation of $W_{E_{\{\mu\}_G}}$, second we check that the conjecture satisfies a transitivity property, and third we give an example justifying the extra character twist appearing in our formulation. This twist is not present in the original formulation of the conjecture.  
\end{remark}
We first explain why the right-hand side is a representation of $W_{E_{\{ \mu\}_G}}$. We start with a general lemma.
\begin{lemma}{\label{indact}}
Suppose a group $\Lambda$ acts on a finite set $S$. Suppose further that for each $s \in V$, we attach a vector space $V_s$ and for each $\lambda \in \Lambda$ and $s \in S$ we have an isomorphism
\[
i(s, \lambda): V_{s} \to V_{\lambda(s)}.
\]
We suppose further that $i(s, 1)$ is the identity map and that $i(\lambda_1(s),\lambda_2) \circ i(s, \lambda_1)=i(s, \lambda_2\lambda_1)$. Then $\bigoplus\limits_{s \in S} V_s$ is naturally a representation of $\Lambda$. 

Let $\{s_1, ..., s_k\} \subset S$ be a set of one representative from each $\Lambda$-orbit in $S$. Then
\[
\bigoplus\limits_{s \in S} V_s \cong \bigoplus\limits^k_{i=1} \mathrm{Ind}^{\Lambda}_{\mathrm{stab}(s_i)} V_{s_i},
\]
where $\mathrm{Ind}$ refers to the induced representation (\emph{not} parabolic induction).
\end{lemma}
\begin{proof}
The proof is clear from the definition of induced representation.
\end{proof}
Moreover, we record the following transitivity property for later use.
\begin{lemma}{\label{indtrans}}
Suppose that $\Lambda$ acts on $S$ as before. Let $S_1 \coprod ... \coprod S_k=S$ be a partition of $S$ so that $\Lambda$ acts on $\{S_1, ..., S_k\}$. Suppose we have for each $s \in S$ a vector space $V_s$ and isomorphisms $i(s,\lambda)$ as above. Then by Lemma \ref{indact} we can consider the $\mathrm{stab}(S_i) \subset \Lambda$ representation $V_{S_i}=\bigoplus\limits_{s \in S_i} $. For each $\lambda \in \Lambda$, we get isomorphisms $i(S_i, \lambda): V_{S_i} \to V_{\lambda(S_i)}$. Thus, again by Lemma \ref{indact}, we get a $\Lambda$ representation $\bigoplus\limits_i V_{S_i}$. This representation is isomorphic to the $\Lambda$  representation $\bigoplus\limits_{s \in S} V_s$ we get from applying Lemma \ref{indact} to $S$. 
\end{lemma}

Now we discuss the $W_{E_{\{\mu\}_G}}$-action in the Harris-Viehmann conjecture. Observe that for $\mu \in X_*(G)$, if $\gamma \in W_{E_{\{ \mu\}_G}}$ stabilizes $\{\mu\}_{M_S}$ then it also stabilizes $\{\mu\}_G$ so that $W_{E_{\{\mu\}_{M_S}}} \subset W_{E_{\{\mu \}_G}}$. 

Now we claim that $W_{E_{\{\mu\}_G}}$ acts on $\mathcal{I}^{G, \mu}_{M_S, b_S}$ and that the stabilizer of $(M_S, \mu_S)$ under this action is $W_{E_{\{\mu\}_{M_S}}}$. To prove the first part of the claim, we pick $\gamma \in W_{E_{\{\mu\}_G}}$ and observe that since $M_S$ and $P_S$ are defined over $\mathbb{Q}_p$, we have $\gamma(M_S)=M_S$ and $\gamma(\mu_S)$ is dominant in $M_S$. Thus $(M_S, \gamma(\mu_S)) \in \mathcal{C}_{M_S}$ so we need only check that $b_S \in \mathbf{B}(M_S, \gamma(\mu_S))$ and $\gamma(\mu_S) \sim_G \mu$. The first check follows from the fact that 
\[
\theta_T(\mu_S)=\theta_T(\gamma(\mu_S)),
\]
and 
\[
\mu_S|_{Z(\widehat{M_S})^{\Gamma}}=\gamma(\mu_S)|_{Z(\widehat{M_S})^{\Gamma}}.
\]
The second check follows because $\gamma$ stabilizes $\{\mu\}_G$.

To prove the second part of the claim, we note that if $\mu_S=\gamma(\mu_S)$ then $\gamma$ stabilizes $\{\mu_S\}_{M_S}$. Conversely, if $\gamma$ stabilizes $\{\mu_S\}_{M_S}$ then since it maps dominant elements relative to $M_S$ to dominant elements, we must have $\gamma(\mu_S)=\mu_S$.

We observe that we have now shown that $W_{E_{\{\mu\}_G}}$ acts on the collection of Rapoport-Zink data $(M_S, b_S, \mu_S)$ for $(M_S, \mu_S) \in \mathcal{I}^{G, \mu}_{M_S, b_S}$. By \cite[Proposition 5.3.iv]{RV1}, these actions induce morphisms of the corresponding towers of rigid spaces and therefore the spaces $H^{\bullet}(M_S, b_S, \mu_S)[\rho]$ for $\rho \in \mathrm{Groth}(J_b(\mathbb{Q}_p))$. Thus by Lemma \ref{indact} we get an action of $W_{E_{\{\mu\}_G}}$ on the sum of vector spaces
\[
\sum\limits_{(M_S, \mu_S) \in \mathcal{I}^{G, \mu}_{M_S, b_S}} H^{\bullet}(M_S, b_S, \mu_S)[\rho],
\]
and therefore on
\[
\sum\limits_{(M_S, \mu_S) \in \mathcal{I}^{G, \mu}_{M_S, b_S}} \mathrm{Ind}^G_{P_S} (H^{\bullet}(M_S, b_S, \mu_S)[\rho]).
\]
We remark that the character twist by $- \mathrm{dim} \mathcal{M}^{rig}_{b,\mu, U}$ in the definition of $H^{\bullet}(M_S, b_S, \mu_S)$ is not an obstacle to defining the $W_{E_{\{\mu\}_G}}$-action as the dimensions of the spaces associated to $(M_S, b_S, \mu_s)$ and $(M_S, b_S, \gamma(\mu_S))$ are the same (for $\gamma \in W_{E_{\{\mu\}_G}}$). Also we observe that the twist by  $[1][|\cdot|^{\langle \rho_G, \mu_S \rangle - \langle \rho_G, \mu \rangle}]$ is harmless as it is constant over orbits of $W_{E_{\{\mu_S\}_G}}$. This concludes our discussion of the $W_{E_{\{\mu\}_G}}$ action.

We now check that the Harris-Viehmann conjecture is transitive. By this, we mean that if we have standard Levi subgroups $M_{S_1}$ and $M_{S_2}$ of $G$ such that $M_b \subset M_{S_2} \subset M_{S_1} \subset G$, then first applying the conjecture to $(G,b,\mu)$ and the inclusion $M_{S_1} \subset G$ and then applying the conjecture to each resulting $(M_{S_1}, b_{S_1}, \mu_{S_1})$ for $(M_{S_1}, \mu_{S_1}) \in \mathcal{I}^{G, \mu}_{M_{S_1}, b_{S_1}}$ and the inclusion $M_{S_2} \subset M_{S_1}$ should be the same as applying the conjecture to $(G, b, \mu)$ and the inclusion $M_{S_2} \subset G$. 

We need to check that the character twists match, that 
\[
\mathcal{I}^{G, \mu}_{M_{S_2}, b_{S_2}} =\{ (M_{S_2}, \mu_{S_2}) \in \mathcal{C}_{M_{S_2}} : (M_{S_2}, \mu_{S_2}) \in \mathcal{I}^{M_{S_1},  \mu_{S_1}}_{M_{S_2}, b_{S_2}} \,\ \text{for some} \,\ (M_{S_1}, \mu_{S_1}) \in \mathcal{I}^{G, \mu}_{M_{S_1}, b_{S_1}} \}.
\]
and that the $W_{E_{\{\mu\}_G}}$ actions are the same.

To check the characters match, it suffices to check that for $(M_{S_1}, \mu_{S_1}), (M_{S_2}, \mu_{S_2}) \in \mathcal{C}_G$ such that $(M_{S_2}, \mu_{S_2}) \leq (M_{S_1}, \mu_{S_1}) \leq (G, \mu)$, we have
\[
\langle \rho_G, \mu_{S_2} \rangle - \langle \rho_G, \mu \rangle = (\langle \rho_G, \mu_{S_1} \rangle - \langle \rho_G, \mu \rangle ) + (\langle \rho_{M_{S_1}}, \mu_{S_2} \rangle - \langle \rho_{M_{S_1}}, \mu_{S_1} \rangle).
\]
This reduces to showing the equality
\begin{align}{\label{transtwist}}
\langle \rho_{G \setminus M_{S_1}}, \mu_{S_1} \rangle= \langle \rho_{G \setminus M_{S_1}}, \mu_{S_2} \rangle,
\end{align}
where $\rho_{G \setminus M_{S_1}}$ is the half-sum of the absolute roots of $G$ that are not roots of $M_{S_1}$.  Since $\mu_{S_2}$ and $\mu_{S_1}$ are conjugate in $M_{S_1}$, there exists a $w \in W^{\mathrm{abs}}_{M_{S_1}}$ so that $w (\mu_1)= \mu_2$. Then the desired equality follows from the fact that the pairing $\langle \cdot, \cdot \rangle$ is $W^{\mathrm{abs}}_{M_{S_1}}$-invariant and that $W^{\mathrm{abs}}_{M_{S_1}}$ stabilizes the set of positive absolute roots in $G$ but not $M_{S_1}$. To prove this second fact, note that $M_{S_1}$ normalizes the unipotent radical $U_{S_1}$ of $P_{S_1}$ and that the roots of $\mathrm{Lie}(U_{S_1})$ are precisely the positive absolute roots of $G$ that are not contained in $M_{S_1}$. 

The second check is precisely Proposition \ref{Itrans}, and the third check follows from Proposition \ref{Itrans} and Lemma \ref{indtrans}.

Now we compute an example to illustrate the necessity of the extra Tate twist in our statement of Conjecture \ref{3.4}. The following example is also discussed in \cite[\S 8.3]{Shi1}

\begin{example}{\label{HVex}}
Let $n_1<n_2$ be coprime positive integers and let $G=\mathrm{GL}_{n_1+n_2}$. Fix $T$ the standard maximal torus of diagonal matrices and $B$ the Borel subgoup of upper triangular matrices. Let $\mu$ be the minuscule cocharacter with weight vector $(1^2,0^{n_1+n_2-2})$ and $b \in \mathbf{B}(G, \mu)$ satisfying $\nu_b=({(1/n_1)}^{n_1}, {(1/n_2)}^{n_2})$. Let $\rho_1, \rho_2$ be supercuspidal representations of $\mathrm{GL}_{n_1}(\mathbb{Q}_p), \mathrm{GL}_{n_2}(\mathbb{Q}_p)$ respectively. Define the standard Levi subgroup $M_b=\mathrm{GL}_{n_1} \times \mathrm{GL}_{n_2}$, and consider the representation $\pi=I^G_{M_b}(\rho_1 \boxtimes \rho_2)$. We will be interested in computing $\mathrm{Mant}_{G,b, \mu}(\mathrm{Red}_b(\pi))$.

The key point is that we can use Shin's formula (Theorem \ref{3.2} of this paper) and known cases of the Harris-Viehmann conjecture due to Mantovan (\cite{Man1}) to do this computation, even though the Harris-Viehmann conjecture is not known to be true in the case of $M_b$ since $b$ is not of Hodge-Newton type.

We observe that there are only $3$ elements $b' \in \mathbf{B}(G, \mu)$ that satisfy
\[
\mathrm{Mant}_{G,b',\mu}(\mathrm{Red}_{b'}(\pi)) \neq 0.
\]
After all, the fact that $\rho_1, \rho_2$ are supercuspidal and the geometric lemma of Bernstein-Zelevinski ($\S 2.11$ of \cite{Ber1}) forces $M_{b'}$ to be one of $G, \mathrm{GL}_{n_1} \times \mathrm{GL}_{n_2}, \mathrm{GL}_{n_2} \times \mathrm{GL}_{n_1}$. In the case where $M_{b'}=G$, we also get $0$ since $LJ(\pi)=0$. Thus, if we write out Shin's formula for our $\pi$, the only elements of $\mathbf{B}(G, \mu)$ whose terms contribute to the left-hand side are $b, b_1, b_2$ where $b$ is as before and $b_1, b_2$ are defined by 
\[
\nu_{b_1}=({(2/n_1)}^{n_1},0^{n_2}), \nu_{b_2}=({(2/n_2)}^{n_2}, 0^{n_1}).
\]
Thus, we have $M_{b_1}=M_b=\mathrm{GL}_{n_1} \times \mathrm{GL}_{n_2}$ and  $M_{b_2}=\mathrm{GL}_{n_2} \times \mathrm{GL_{n_1}}$. Note that $b_1, b_2$ are both of Hodge-Newton type so that we can apply the results of Mantovan.

We have 
\[
\mathrm{Mant}_{G, b_1, \mu}( \mathrm{Red}_{b_1} (\pi))=\mathrm{Mant}_{G, b_1, \mu}(LJ(\delta^{\frac{1}{2}}_{P_{b_1}} \otimes J^G_{P^{op}_{b_1}} I^G_{M_{b_1}} (\rho_1 \boxtimes \rho_2))).
\]
By the geometric lemma of Bernstein-Zelevinski ($\S 2.11$ of \cite{Ber1}) we have that the above equals
\[
\mathrm{Mant}_{G, b_1, \mu_1}(LJ ((\rho_1 \boxtimes \rho_2) \otimes \delta^{\frac{1}{2}}_{P_{b_1}})).
\]
We recall that $\delta_{P_{b_1}} = (| \cdot |^{n_2} \circ  \det )\boxtimes (| \cdot |^{-n_1} \circ \det) $ and henceforth use the notation $\rho(n)$ to mean $(|\cdot|^{n} \circ \det )\otimes \rho$. Thus, we can rewrite the above formula as
\[
\mathrm{Mant}_{G, b_1, \mu_1}(LJ (\rho_1(n_2/2)) \boxtimes LJ(\rho_2(-n_1/2))).
\]
Then applying the Harris-Viehmann formula we get that the above equals
\begin{align}{\label{3.1}}
    \mathrm{Ind}^G_{M_b}(\mathrm{Mant}_{\mathrm{GL}_{n_1}, (1^2, 0^{n_1-2})}(LJ (\rho_1(n_2/2))) \boxtimes \mathrm{Mant}_{\mathrm{GL}_{n_2}, (0^{n_2})}(LJ( \rho_2(-n_1/2)))).
\end{align}
Since $\rho_1$ and $\rho_2$ are supercuspidal, we can compute (by an easy application of Shin's formula for instance) that
\[
\mathrm{Mant}_{\mathrm{GL}_{n_1}, (1^2, 0^{n_1-2})}(LJ (\rho_1(n_2/2)))=[\rho_1(n_2/2)][r_{(-1^2, 0^{n_1-2})} \circ LL(\rho_1(n_2/2))\otimes |\cdot |^{2-n_1}],
\]
and so Equation \eqref{3.1} becomes equal to
\begin{align*}
    [\pi][r_{(-1^2, 0^{n_1-2})} \circ LL(\rho_1(n_2/2))\otimes | \cdot |^{2-n_1} \otimes r_{(0^{n_2})} \circ LL(\rho_2(-n_1/2))].
\end{align*}
Pulling the twists through the $r_{-\mu}$ maps, we get
\begin{align*}
    [\pi][(r_{(-1^2, 0^{n_1-2})} \boxtimes r_{(0^{n_2})}) \circ (LL(\rho_1) \oplus LL(\rho_2))\otimes | \cdot |^{2-n_1-n_2}].
\end{align*}
Repeating this computation for the $b_2$ term, we get 
\[
\mathrm{Mant}_{G, b_2, \mu}(\mathrm{Red}_{b_2}(\pi))
\]
\[
=[\pi][(r_{(-1^2, 0^{n_2-2})} \boxtimes r_{(0^{n_1})}) \circ (LL(\rho_2) \oplus LL(\rho_1))\otimes | \cdot |^{2-n_1-n_2}].
\]

We now compare these terms to the righthand side of Shin's formula. There the term is
\[
[\pi][r_{-\mu} \circ LL(\pi) \otimes |\cdot|^{2-n_1-n_2}].
\]
Now $LL(\pi)=LL(\rho_1) \oplus LL(\rho_2)$. Thus, we can restrict $r_{-\mu}$ to $\widehat{M_b} \subset \widehat{G}$ (we have been ignoring the Galois part of $^LG$ in this example since $G$ is a split group). Using the theory of weights, we get
\[
r_{-\mu}|_{\widehat{M}}=[r_{(-1^2,0^{n_1-2})} \boxtimes r_{(0^{n_2})} ]\oplus [r_{(-1,0^{n_1-1})} \boxtimes r_{(-1,0^{n_2-1})} ] \oplus [r_{(0^{n_1})} \boxtimes r_{(-1^2,0^{n_2-2})}],
\]
and so we see that the contributions for $b_1, b_2$ which we computed above will cancel terms on the righthand side of Shin's formula leaving us with
\[
\mathrm{Mant_{G,b,\mu}}(\mathrm{Red}_b(\pi))=[\pi][(r_{(-1,0^{n_1-1})} \boxtimes r_{(-1,0^{n_2-1})}) \circ (LL(\rho_1) + LL(\rho_2)) \otimes | \cdot |^{2-n_1-n_2}].
\]
However, if the Harris-Viehmann conjecture without the extra Tate twist were to hold for $b$, we would get
\[
\mathrm{Mant}_{G,b,\mu}(\mathrm{Red}_b(\pi))=\mathrm{Mant}_{G,b,\mu}(LJ(\rho_1(n_2/2)) \boxtimes LJ(\rho_2(-n_1/2)))
\]
\[
=[\pi][r_{(-1, 0^{n_1-1})} \boxtimes r_{(-1,0^{n_2-1})} \circ (LL(\rho_1)+LL(\rho_2))|\cdot|^{1-n_2}].
\]
Thus, we see the Tate twists do not agree. 
\end{example}

In general, the righthand side of Shin's formula has a twist of $-\langle \rho_G, \mu \rangle $ where $\rho_G$ is the half sum of the positive roots of $G$.  Suppose now that $b \in \mathbf{B}(G, \mu)$ and $b' \in \mathbf{B}(M_b)^+$ corresponds to $b$ under Equation \eqref{BGiso}. Then for any $(M_b, \mu') \in \mathcal{I}^{G, \mu}_{M_b, b'}$, we would expect the Galois part of $\mathrm{Mant}_{M_b, b', \mu'}(\rho)$ for $\rho \in \mathrm{Groth}(J_b(\mathbb{Q}_p))$ to come with a twist of $-\langle \rho_{M_b}, \mu' \rangle $. Then the Galois part of $\mathrm{Mant}_{G, b, \mu} (\mathrm{Red}_b(\pi))$ for $\pi \in \mathrm{Groth}(G(\mathbb{Q}_p))$ would carry an extra twist of $-\langle \frac{\det(Ad_{N_b}(M_b))|_T}{2}, \mu' \rangle$ corresponding to twisting $J^G_{P^{op}_b}(\pi)$ by $\delta^{\frac{1}{2}}_{P_b}$ in the definition of $\mathrm{Red}_b$ . We note that 
\[
\langle \rho_{M_b}, \mu' \rangle+\langle \frac{\det(Ad_{N_b}(M_b))|_T}{2}, \mu' \rangle= \langle \rho_G, \mu' \rangle,
\]

Thus, we see that the difference between these Tate twists is 
\[
\langle \rho_G, \mu' \rangle - \langle \rho_G, \mu \rangle.
\]
which is the twist in Conjecture \ref{3.4}

\begin{remark}
We note that in the Hodge-Newton case studied by Mantovan, $\mu=\mu'$ (as in the notation of the previous paragraph) so that this extra twist vanishes, agreeing with Mantovan's results (\cite[Corollary 5]{Man1}, cf. \cite[Theorem 8.8]{RV1}).
\end{remark}

We now give an alternate version of the Harris-Viehmann conjecture that we will use in numerous arguments in this paper. Suppose that $G ,b, \mu$ are as in Theorem \ref{3.2}. The standard Levi subgroup $M_b$ has a natural product decomposition 
\[
M_b=M_1 \times... \times M_k
\]
so that under the natural isomorphism 
\[
\mathbf{B}(M_b) \cong \mathbf{B}(M_1) \times \mathbf{B}(M_k), b' \mapsto (b'_1, ..., b'_k),
\]
each $\nu(b_i)$ has a single slope.  Now pick $(M_b, \mu_b) \in \mathcal{I}^{G, \mu}_{M_b, b'}$. Then the  local Shimura variety datum $(M_b, b', \mu_b)$ decomposes into a collection $(M_1, b'_1, \mu_{b,1}),...,(M_k, b'_k, \mu_{b,k})$. In $\S 5.2.(ii)$ of \cite{RV1}, the authors show that the local Shimura variety associated to $(M_b, b', \mu_b)$ is the product of those associated to $(M_i, b'_i, \mu_{b,i})$. Furthermore using the Kunneth formula  (as in \cite[15]{Man1}), we get that for $\rho_1 \boxtimes ... \boxtimes \rho_k \in \mathrm{Groth}(M_1(\mathbb{Q}_p) \times ... \times M_k(\mathbb{Q}_p))$,
\[
\mathrm{Mant}_{M_b, b' , \mu_b}(\rho_1 \boxtimes ... \boxtimes \rho_k)= \boxtimes^k_{i=1} \mathrm{Mant}_{G_i, b'_i, \mu_{b,i}}(\rho_i),
\]
as a representation of $M_b \times W_{E_{\{\mu_b\}_{M_b}}}$ (the group $W_{E_{\{\mu_b\}_{M_b}}}$ acts diagonally through the product $W_{E_{\{\mu_{b,1}\}_{M_1}}} \times ... \times W_{E_{\{\mu_{b,k}\}_{M_k}}}$).

Thus, we have the following alternate form of the Harris-Viehmann conjecture for the Rapoport-Zink spaces we consider.
\begin{conjecture}[Alternate Form of Harris-Viehmann Conjecture]
We use the notation of the previous paragraphs so that in particular, $(G, b , \mu)$ comes from an unramified Rapoport-Zink space of EL-type as in Definition \ref{RZdat}. Then for any $\rho \in \mathrm{Groth}(J_b(\mathbb{Q}_p))$, we have the following equality in $\mathrm{Groth}(G(\mathbb{Q}_p) \times W_{E_{\{\mu\}_G}})$:
\[
\mathrm{Mant}_{G,b, \mu}(\rho)= \sum\limits_{(M_b, \mu_b) \in \mathcal{I}^{G, \mu}_{M_b, b'}} \mathrm{Ind}^G_{P_b}(\boxtimes^k_{i=1} \mathrm{Mant}_{M_b, b'_i, \mu_{b,i}}(\rho_i)) \otimes [1][|\cdot|^{\langle \rho_G, \mu_b \rangle - \langle \rho_G, \mu \rangle}].
\]
\end{conjecture}
\subsection{Proof of Theorem \ref{1.3}}
The combination of the Harris-Viehmann conjecture and sum formula allows us to relate the cohomology of Rapoport-Zink spaces to the cocharacter pairs studied in \S2. To do so, we attach a map of Grothendieck groups to each cocharacter pair. We return to the notation of \S3.1.

Fix a cocharacter pair $(G, \mu) \in \mathcal{C}_G$. Suppose $(M_S, \mu_S) \in \mathcal{C}_G$ and satisfies $\mu_S \sim_G \mu$. We associate $(M_S, \mu_S)$ to a map of representations 
\[
[M_S, \mu_S]: \mathrm{Groth}(G(\mathbb{Q}_p)) \to \mathrm{Groth}(G(\mathbb{Q}_p) \times W_{E_{\{\mu_S\}_{M_S}}}),
\]
given by 
\[
\pi \mapsto (\mathrm{Ind}^G_{P_S} \circ [\mu_S] \circ ( \delta_{P_S} \otimes \mathrm{Jac}^G_{P^{op}_S} ))(\pi)\otimes [1][|\cdot|^{\langle \rho_G, \mu_S \rangle - \langle \rho_G, \mu \rangle}],
\]
with
\[
[\mu_S]: \mathrm{Groth}(M_S(\mathbb{Q}_p)) \to \mathrm{Groth}(M_S(\mathbb{Q}_p) \times W_{E_{\{ \mu\}_{M_S}}}),
\]
given by 
\[
\pi \mapsto [\pi][r_{-\mu_S} \circ LL(\pi)|_{W_{E_{\{\mu_S\}_{M_S}}}} \otimes | \cdot | ^{-\langle \rho_{M_S}, \mu_S \rangle }].
\]
\begin{remark}
We note that the map $[M_S, \mu_S]$ is only defined relative to a cocharacter pair $(G, \mu)$.
\end{remark}
\begin{remark}{\label{twister}}
We observe an interesting property of the maps $[M_S, \mu_S]$. Fix $(G, \mu)$ and consider $(M_S, \mu_S)$ such that $\mu_S \sim_G \mu$. Since the normalized Jacquet module and parabolic induction functors behave better with respect to the local Langlands correspondence, it makes sense to rewrite $[M_S, \mu_S]$ in terms of these maps. We get
\[
[M_S, \mu_S]= (I^G_{M_S} \otimes \delta^{-\frac{1}{2}}_{P_S} \circ [\mu_S] \circ (\delta^{\frac{1}{2}}_{P_S} \otimes J^G_{P^{op}_S})) \otimes [1][|\cdot|^{\langle \rho_G, \mu_S - \mu \rangle}].
\]
Note that the twists by the modular character cancel in the admissible part but do not cancel in the Galois part. Thus, the total Tate twist of the Galois part is 
\[
\langle \rho_G, \mu_S - \mu \rangle - \langle \rho_{M_S}, \mu_S \rangle - \langle \frac{\det (Ad_{N_S}(M_S))|_T}{2}, \mu_S \rangle 
\]
\[
=- \langle \rho_G, \mu \rangle.
\]
This twist does not depend on $(M_S, \mu_S)$ but rather only on $(G, \mu)$. Thus, as we will see in the computations of the next section, it is possible for large cancellations to occur in computations of $\mathrm{Mant}_{G, b, \mu}(\rho)$ for various $\rho$.
\end{remark}

We now prove some lemmas relating to these maps before tackling the main theorem.
\begin{lemma}{\label{combind}}
Let $M_{S_1}, M_{S_2}$ be standard Levi subgroups of $G$ satisfying $M_{S_2} \subset M_{S_1}$. Consider the natural map 
\[
i^G_{M_{S_1}}: \mathcal{C}_{M_{S_1}} \to \mathcal{C}_G,
\]
as defined in Equation \eqref{inclusion}. Let $(M_{S_2}, \mu_{S_2}) \in \mathcal{C}_{M_{S_1}}$.  Suppose further that we have fixed pairs $(M_{S_1}, \mu_{S_1}) \in \mathcal{C}_{M_{S_1}}$ and $(G, \mu) \in \mathcal{C}_G$ so that $\mu_{S_2} \sim_{M_{S_1}} \mu_{S_1}$ and $\mu_{S_2} \sim_G \mu$. Then for $\pi \in \mathrm{Groth}(G_{\mathbb{Q}_p})$,
\[
i^G_{M_{S_1}}([M_{S_2}, \mu_{S_2}])(\pi)= (\mathrm{Ind}^G_{P_{S_1}} \circ [M_{S_2}, \mu_{S_2}] \circ ( \delta_{P_{S_1}} \otimes \mathrm{Jac}^G_{P^{op}_{S_1}}))(\pi) \otimes [1][| \cdot |^{\langle \rho_G, \mu_{S_1} \rangle - \langle \rho_G, \mu \rangle} ],
\]
where we write
\[
i^G_{M_{S_1}}([M_{S_2}, \mu_{S_2}]): \mathrm{Groth}(G(\mathbb{Q}_p)) \to \mathrm{Groth}(G(\mathbb{Q}_p) \times W_{E_{\{\mu_{S_2}\}_{M_{S_2}}}}),
\]
to denote the map associated to $i^G_{M_{S_1}}((M_{S_2}, \mu_{S_2}))$ in the manner above.
\end{lemma}
\begin{proof}
We first note that by transitivity of the Jacquet module and modulus character constructions, we have
\begin{equation*}
    \delta_{P_{S_2}} \otimes \mathrm{Jac}^G_{P^{op}_{S_2}} = (\delta_{P_{S_2} \cap M_1} \otimes \mathrm{Jac}^{M_{S_1}}_{P^{op}_{S_2}}) \circ (\delta_{P_{S_1}} \otimes \mathrm{Jac}^G_{P^{op}_{S_1}}).
\end{equation*}
Hence, we just need to check that the twists on the Galois parts of both sides match. By Remark \ref{twister}, both twists are by $-\langle \rho_G, \mu \rangle$
\end{proof}
\begin{lemma} \label{3.7}
Suppose we are in the situation of Proposition \ref{prodcochar} so that $G=G_1 \times ... \times G_k$ is a connected reductive group with standard Levi subgroup $M_S=M_{S_1} \times ... \times M_{S_k}$. Fix cocharacter pairs $(M_S, \mu_S), (G, \mu) \in \mathcal{C}_G$ with $\mu_S \sim_G \mu$. The bijection $\mathcal{C}_G \cong \mathcal{C}_{G_1} \times ... \mathcal{C}_{G_k}$ takes $(M_S, \mu_S)$ to $((M_{S_1}, \mu_{S_1}), ..., (M_{S_k}, \mu_{S_k}))$ and $(G, \mu)$ to $((G_1, \mu_1), ..., (G_k, \mu_k))$ and we have $\mu_{S_i} \sim_{G_i} \mu_i$. Then we define
\[
\boxtimes^k_{i=1} [M_{S_i}, \mu_{S_i}] : \mathrm{Groth}(G(\mathbb{Q}_p)) \to \mathrm{Groth}(G(\mathbb{Q}_p) \times W_{E_{\{\mu_S\}_{M_S}}})
\]
by
\[
\pi_1 \boxtimes ... \boxtimes \pi_k \mapsto [M_{S_1}, \mu_{S_1}](\pi_1) \boxtimes ... \boxtimes [M_{S_k}, \mu_{S_k}](\pi_k).
\]
Then we have the following equality of homomorphisms of Grothendieck groups:
\[
\boxtimes^k_{i=1} [M_{S_i}, \mu_{S_i}] = [M_S, \mu_S]
\]
\end{lemma}
\begin{proof}
We have 
\begin{align*}
\boxtimes^k_{i=1} [M_{S_i}, \mu_{S_i}] &= \boxtimes^k_{i=1} \mathrm{Ind}^{G_i}_{P_{S_i}} \circ [\mu_{S_i}] \circ (\delta_{P_{S_i}}\otimes \mathrm{Jac}^{G_i}_{P^{op}_{S_i}})\otimes [1][|\cdot|^{\langle \rho_{G_i}, \mu_{S_i} - \mu_i \rangle}] \\
&= \mathrm{Ind}^G_{P_S} \circ [\mu] \circ (\delta_{P_S} \otimes \mathrm{Jac}^G_{P^{op}_S}) \otimes [1][| \cdot|^{\sum\limits^k_{i=1} \langle \rho_{G_i}, \mu_{S_i} - \mu_i \rangle}]\\
&=\mathrm{Ind}^G_{P_S} \circ [\mu] \circ (\delta_{P_S} \otimes \mathrm{Jac}^G_{P^{op}_S}) \otimes [1][| \cdot|^{ \langle \rho_{G}, \mu_S - \mu \rangle}]\\
&=[M_S, \mu_S].
\end{align*}
\end{proof}

For some finite subset $C \subset \mathcal{C}_G$, such that each $(M_S, \mu_S) \in C$ satisfies $\mu_S \sim_G \mu$, we would like to make sense of a sum
\[
\sum\limits_{(M_S, \mu_S) \in C} [M_S, \mu_S].
\]

This makes sense as a map $\mathrm{Groth}(G(\mathbb{Q}_p)) \to \mathrm{Groth}(G(\mathbb{Q}_p) \times W_E)$ where $W_E = \bigcap\limits_{{(M_S, \mu_S) \in C}} W_{E_{\{\mu_S\}_{M_S}}}$. However, for our purposes, we would like to understand when we can extend the image of this map to a representation in $\mathrm{Groth}(G(\mathbb{Q}_p) \times W_{E_{\{\mu\}_G}})$.
\begin{lemma}{\label{galact}}
Fix a pair $(G, \mu) \in \mathcal{C}_G$. Consider a finite subset $C \subset \mathcal{C}_G$ such that if $(M_S, \mu_S) \in C$ then $\mu_S \sim_G \mu$. Furthermore, suppose that for each $\gamma \in W_{E_{\{\mu\}_G}}$ and element $(M_S, \mu_S) \in C$, we have $(M_S, \gamma(\mu_S)) \in C$. Then 
\[
\sum\limits_{(M_S, \mu_S) \in C} [M_S, \mu_S],
\]
is a map
\[
\mathrm{Groth}(G(\mathbb{Q}_p)) \to \mathrm{Groth}(G(\mathbb{Q}_p) \times W_{E_{\{\mu\}_G}})
\]
in a natural way.
\end{lemma}
\begin{proof}
Our construction is analogous to that of Lemma \ref{indact}.
We fix $\rho \in \mathrm{Groth}(G(\mathbb{Q}_p))$ and give 
\[
V_C=\bigoplus\limits_{(M_S, \mu_S) \in C} [M_S,\mu_S](\rho),
\]
the structure of a $G(\mathbb{Q}_p) \times W_{E_{\{\mu\}_G}}$ representation. Suppose that $C=C_1 \coprod ... \coprod C_n$ where each $C_i$ is a single $W_{E_{\{\mu\}_G}}$-orbit. Then for each $i$, we give
\[
V_{C_i}=\bigoplus\limits_{(M_S, \mu_S) \in C_i} [M_S,\mu_S](\rho),
\]
the structure of a $G(\mathbb{Q}_p) \times W_{E_{\{\mu\}_G}}$-representation and then define the $G(\mathbb{Q}_p) \times W_{E_{\{\mu\}_G}}$-structure on $V_C$ to be the direct sum of the $V_{C_i}$.

Suppose now that $C$ contains a single $W_{E_{\{\mu\}_G}}$ orbit. In this case, we will show that 
\[
\bigoplus\limits_{(M_S, \mu_S) \in C} [M_S,\mu_S](\rho),
\]
can be given the structure of a $\mathrm{Groth}(G(\mathbb{Q}_p) \times W_{E_{\{\mu\}_G}})$ representation equal to 
\[
[\mathrm{Ind}^G_{P_S}(\delta_{P_S} \otimes \mathrm{Jac}^G_{P^{op}_S}(\rho))][r \circ LL(\delta_{P_S} \otimes \mathrm{Jac}^G_{P^{op}_S}(\rho))|_{W_{E_{\{\mu\}_G}}} \otimes | \cdot |^{-\langle \rho_G, \mu_S - \mu \rangle - \langle \rho_{M_S}, \mu_S \rangle}],
\]
where $r$ is the induced representation (\emph{not} parabolic induction) given by 
\[
\mathrm{Ind}^{\widehat{M_S} \rtimes W_{E_{\{ \mu \}_G}}}_{\widehat{M_S} \rtimes W_{E_{\{\mu_S\}_{M_S}}}} (r_{-\mu_S}),
\]
for a fixed choice of $(M_S, \mu_S) \in C$. The isomorphism class of $r$ will not depend on this choice.

We study the representation $r$. Fix representatives $\gamma_1,..., \gamma_k \in W_{E_{\{\mu\}_G}}/W_{E_{\{\mu_S\}_{M_S}}}$ so that $\gamma_1=1$. Then $r$ is defined to be the sum of $k$ copies of $r_{-\mu_S}$ indexed by the $\gamma_i$ and acted on by $W_{E_{\{\mu\}_G}}$ in the standard way. We check that the $i$th copy of $r_{-\mu_S}$ is a representation of $\widehat{M_S} \rtimes W_{E_{\{\gamma_i(\mu_S)\}_{M_S}}}$ and isomorphic to $r_{- \gamma_i(\mu_S)}$. Let $V_i$ be the underlying vector space of the $i$th copy of $r_{-\mu_S}$. Then $V_i$ is naturally a representation of $\widehat{M_S} \rtimes \gamma_i W_{E_{\{\mu_S\}_{M_S}}}\gamma^{-1}_i = \widehat{M_S} \rtimes W_{E_{\{\gamma_i(\mu_S) \}_{M_S}}}$. 

Now suppose $v \in V_1$ is a weight vector of $\widehat{T} \subset \widehat{M_S}$ of weight $\mu'$. Then we show that $(1, \gamma_i)v \in V_i$ has weight $\gamma_i(\mu')$. After all, for $t \in \widehat{T}$, we have 
\begin{align*}
    r((t,1))( (1, \gamma_i)v ) &=(t,\gamma_i)v \\
    &=(1, \gamma_i)(\gamma^{-1}_i(t),1)v\\
    &=(1, \gamma_i)r_{-\mu_S}((\gamma^{-1}_i(t),1))(v)\\
    &=(1,\gamma_i) \mu'(\gamma^{-1}_i(t))v\\ &=\gamma_i(\mu')(t) (1, \gamma_i)v.
\end{align*}
In particular, we have shown that $V_i$ is irreducible of extreme weight $-\gamma_i(\mu_S)$ as an $\widehat{M_S}$-representation (since $r_{-\mu_S}$ is irreducible of extreme weight $- \mu_S$ as an $\widehat{M_S}$-representation). It is a simple check similar to the above that $W_{E_{\{\gamma_i(\mu_S)\}_{M_S}}}$ acts trivially on the highest weight space of $V_i$. This proves that $V_i$ is isomorphic to $r_{- \gamma_i(\mu_S)}$.

In particular, this shows that we can give 
\[
\bigoplus\limits_{\gamma_i \in W_{E_{\{\mu\}_G}}/W_{E_{\{\mu_S\}_{M_S}}}} r_{-\gamma_i(\mu_S)} \circ LL(\delta_{P_S} \otimes \mathrm{Jac}^G_{P^{op}_S}(\rho))|_{W_{E_{\{\gamma_i(\mu_S)\}_{M_S}}}},
\]
the structure of a $W_{E_{\{\mu\}_G}}$ representation isomorphic to 
\[
r \circ LL(\delta_{P_S} \otimes \mathrm{Jac}^G_{P^{op}_S}(\rho))|_{W_{E_{\{\mu\}_G}}}.
\]
To conclude the proof, we just need to check that the $|\cdot|$ twists on each $[M_S, \gamma_i(\mu_S)]$-term are the same. This follows because $\rho_G$ and $\rho_{M_S}$ are both invariant by $W_{E_{\{\mu\}_G}}$.
\end{proof}

We would like to check the following:
\begin{lemma}{\label{galM}}
The sum $\mathcal{M}_{G, b, \mu}$ as in Definition \ref{MGBmu} gives a map
\[
[\mathcal{M}_{G, b, \mu}]:\mathrm{Groth}(G(\mathbb{Q}_p)) \to \mathrm{Groth}(G(\mathbb{Q}_p) \times W_{E_{\{\mu\}_G}}),
\]
where 
\begin{equation*}
    [\mathcal{M}_{G, b, \mu}] := \sum\limits_{(M_S, \mu_S) \in \mathcal{R}_{G,b,\mu}} (-1)^{L_{M_S, M_b}}[M_S, \mu_S].
\end{equation*}
\end{lemma}
\begin{proof}
By Lemma \ref{galact}, it suffices to show that $\mathcal{M}_{G,b, \mu}$ is invariant under the natural action of $W_{E_{\{\mu\}_G}}$ on $\mathbb{Z} \langle \mathcal{C}_G \rangle$. Pick $\gamma \in W_{E_{\{\mu\}_G}}$. Since the action of $\gamma$ on a cocharacter pair fixes the standard Levi subgroup in the first factor, signs will not be an issue and we will be done if we can check that $\mathcal{R}_{G,b,\mu}$ is $\gamma$-invariant. But if $(M_b, \mu_b) \in \mathcal{T}_{G,b,\mu}$ then it is a simple consequence of the definition of $\mathcal{T}$ that so is $(M_b, \gamma(\mu_b))$. Furthermore if $(M_S, \mu_S) \leq (M_b, \mu_b)$ then $(M_S, \gamma(\mu_S)) \leq (M_b, \gamma(\mu_b))$ by definition of the partial order relation (remarking that $\theta_{M_S}(\mu_S)=\theta_{M_S}(\gamma(\mu_S))$). This shows that $\mathcal{R}_{G, b, \mu}$ is $\gamma$-invariant as desired.
\end{proof}

If we combine the previous lemma with Proposition \ref{prodcochar}, and Lemma \ref{3.7} we get
\begin{align}{\label{combprod}}
\boxtimes^k_{i=1} [\mathcal{M}_{G_i, b_i, \mu_i}]=[\mathcal{M}_{G,b,\mu}].
\end{align}

We now prove the key result of this section which provides the connection between $\mathrm{Mant}$ and cocharacter pairs.
\begin{theorem}\label{3.8}
Assume that the Harris-Viehmann conjecture is true for the general linear groups we consider. \begin{enumerate}
    \item We have the following equality of morphisms $\mathrm{Groth}^2(G(\mathbb{Q}_p)) \to \mathrm{Groth}^2(G(\mathbb{Q}_p) \times W_{E_{\{\mu\}_G}})$:
\[
 \mathrm{Mant}_{G, b, \mu} \circ \mathrm{Red}_b= [\mathcal{M}_{G, b, \mu}].
\]
where $\mathrm{Groth}^2(G(\mathbb{Q}_p))$ is defined to be the span of the essentially square integrable representations in $\mathrm{Groth}(G(\mathbb{Q}_p))$.

\item Now assume further that Theorem \ref{3.2} holds for all admissible representations of $\mathrm{Groth}(G(\mathbb{Q}_p))$. Then the above equality holds as morphisms $\mathrm{Groth}(G(\mathbb{Q}_p)) \to \mathrm{Groth}(G(\mathbb{Q}_p) \times W_{E_{\{\mu\}_G}})$.
\end{enumerate}
\end{theorem}
\begin{proof}

We prove the second statement first. We prove this result by induction on the rank of $X_*(T)$. 

If the rank of $X_*(T)$ is $1$, then $\mathbf{B}(G, \mu)$ is a singleton and so the result follows from Theorem \ref{3.2}. 

Suppose the result holds for all non-basic $b \in \mathbf{B}(G, \mu)$ with $\mathrm{Rk}(X_*(T)) \leq r$. Then by Theorem \ref{3.2} and Theorem \ref{combinatorial sum}, the result holds for all $b \in \mathbf{B}(G ,\mu)$ with $\mathrm{Rk}(X_*(T)) \leq r$.

Finally, suppose the result holds for all $b \in \mathbf{B}(G, \mu)$ with $\mathrm{Rk}( X_*(T)) \leq r$. Then suppose $X_*(T)$ has rank $r+1$ and choose $b \in \mathbf{B}(G, \mu)$ such that $b$ is not basic. We write $M_b=M_{b_1} \times ... \times M_{b_k}$. By the Harris-Viehmann formula, 
\[
\mathrm{Mant}_{G, b, \mu} \circ \mathrm{Red}_b
\]
\[
=\sum\limits_{(M_b, \mu_b) \in \mathcal{I}^{G, \mu}_{M_b, b'}}(\mathrm{Ind}^G_{P_b} \circ  \otimes^k_{i=1} \mathrm{Mant}_{M_{b_i}, b'_i, \mu_{b i}} \circ \mathrm{Red}_b)\otimes [1][|\cdot|^{\langle \rho_G, \mu_b \rangle - \langle \rho_G, \mu \rangle}]
\]
\[
= \sum\limits_{(M_b, \mu_b) \in \mathcal{I}^{G, \mu}_{M_b, b'}} (\mathrm{Ind}^G_{P_b} \circ \otimes^k_{i=1} (\mathrm{Mant}_{M_{b_i}, b'_i, \mu_{b i}} \circ \mathrm{Red}_{b'_i}) \circ (\delta_{P_b} \otimes \mathrm{Jac}^G_{P^{op}_b}))\otimes [1][|\cdot|^{\langle \rho_G, \mu_b \rangle - \langle \rho_G, \mu \rangle}].
\]
By inductive assumption we get
\[
=\sum\limits_{(M_b, \mu_b) \in \mathcal{I}^{G, \mu}_{M_b, b'}} (\mathrm{Ind}^G_{P_b} \circ \otimes^k_{i=1} [\mathcal{M}_{M_{b_i}, b'_i, \mu_{b i}}]\circ (\delta_{P_b} \otimes \mathrm{Jac}^G_{P^{op}_b}))\otimes [1][|\cdot|^{\langle \rho_G, \mu_b \rangle - \langle \rho_G, \mu \rangle}],
\]
and now by Equation \eqref{combprod} 
\[
=\sum\limits_{(M_b, \mu_b) \in \mathcal{I}^{G, \mu}_{M_b, b'}} (\mathrm{Ind}^G_{P_b}\circ [\mathcal{M}_{M_b, b', \mu_b} ]\circ ( \delta_{P_b} \otimes \mathrm{Jac}^G_{P^{op}_b}))\otimes [1][|\cdot|^{\langle \rho_G, \mu_b \rangle - \langle \rho_G, \mu \rangle}].
\]
Finally, by Corollary \ref{combinatorial induction} and Lemma \ref{combind}
\[
=[\mathcal{M}_{G,b, \mu}].
\]
We must check that the $W_{E_{\{\mu\}_G}}$ structure coming from Remark \ref{Harris-Viehmann} is compatible with that of Lemma \ref{galact}. Pick $\rho \in \mathrm{Groth}(G(\mathbb{Q}_p))$. By inductive assumption and Lemma \ref{combind}, for each $(M_b, \mu_b) \in \mathcal{I}^{G, \mu}_{M_b, b'}$, the $W_{E_{\{\mu_b\}_{M_b}}}$-structures on 
\[
(\mathrm{Ind}^G_{P_b} \circ \mathrm{Mant}_{M_b, b', \mu_b} \circ \mathrm{Red}_{b'} \circ ( \delta_{P_b} \otimes \mathrm{Jac}^G_{P^{op}_b}))(\rho) \otimes [1][|\cdot|^{\langle \rho_G, \mu_b \rangle - \langle \rho_G, \mu \rangle}], 
\]
and 
\[
i^G_{M_b}([\mathcal{M}_{M_b, b', \mu_b}])(\rho),
\]
are the same. Thus by Lemma \ref{indact}, the $W_{E_{\{\mu\}_G}}$-structure on $\mathrm{Mant}_{G,b, \mu}(\mathrm{Red}_b(\rho))$ is a direct sum over the $W_{E_{\{\mu\}_G}}$-orbits of $\mathcal{I}^{G, \mu}_{M_b, b'}$ of induced representations of the form \[
\mathrm{Ind}^{W_{E_{\{\mu\}_G}}}_{W_{E_{\{\mu_b\}_{M_b}}}} i^G_{M_b}([\mathcal{M}_{M_b, b', \mu_b}])(\rho). 
\]
This $W_{E_{\{\mu\}_G}}$-structure matches the one on $[\mathcal{M}_{G,b, \mu}]$ (coming from Lemma \ref{galact}) by the transitivity of the induced representation construction (see Lemma \ref{indtrans} for instance).

We now prove the first statement of the theorem. To do so, we need to show that if we restrict ourselves to the span of the essentially square integrable representations $\mathrm{Groth}^2(G(\mathbb{Q}_p)) \subset \mathrm{Groth}(G(\mathbb{Q}_p))$, then we can remove the first assumption. In particular, these representations are accessible, so we have Theorem \ref{3.2} unconditionally. In the above proof we need only observe that the Jacquet module $\mathrm{Jac}^G_{P^{op}}(\rho)$ is a sum of essentially square integrable representations for $\rho \in \mathrm{Irr}^2(G(\mathbb{Q}_p))$. Thus, to get the result for $\mathrm{Groth}^2(G(\mathbb{Q}_p))$ by induction, our inductive assumption need only hold for all $\mathrm{Groth}^2(G'(\mathbb{Q}_p))$ for $rk G' < rk G$. This shows that under the condition that the Harris-Viehmann conjecture is true in the cases we consider, the theorem is true for essentially square integrable representations without any other assumptions.
\end{proof}

\section{Harris's Generalization of the Kottwitz Conjecture (proof of Theorem 1.5)} \label{section 4}

In this section, we discuss an explicit computation using the above results. In particular, we prove that Shin's formula for all admissible representations combined with the Harris-Viehmann conjecture proves Harris's conjecture for the general linear groups considered in \S3. This conjecture is distinct from the Harris-Viehmann conjecture and is \cite[Conj 5.4]{Har1}.

We begin by discussing the Kottwitz conjecture, which appears as \cite[Cor 7.7]{Shi1} in the cases we consider, and more generally as \cite[Conj 7.3]{RV1}. Fix $G$ as in section 3 of this paper and a cocharacter pair $(G, \mu)$ such that $\mu$ is minuscule. Let $b \in \mathbf{B}(G, \mu)$ be the unique basic element. Now, consider $\rho$ a representation of $J_b(\mathbb{Q}_p)$ such that $JL(\rho)$ is a supercuspidal representation of $G(\mathbb{Q}_p)$. Then
\[
\mathrm{Mant}_{G,b,\mu}(\mathrm{Red}_b(JL(\rho)))=\mathrm{Mant}_{G,b,\mu}(\rho),
\]
but by Theorem \ref{3.8}, the lefthand side equals 
\[
[\mathcal{M}_{G, b , \mu}](JL(\rho)).
\]
Now we see that since $JL(\rho)$ is supercuspidal, each term of the form $[M_S, \mu_S]( JL(\rho))$ is $0$ when $M_S$ is a proper Levi subgroup of $G$. Thus,
\[
\mathrm{Mant}_{G,b,\mu}(\rho)=[\mathcal{M}_{G, b , \mu}](JL(\rho))=[JL(\rho)][r_{-\mu} \circ LL(\rho)| \cdot|^{-\langle \rho_G, \mu \rangle}].
\]
This result is the Kottwitz conjecture for $G$. Alternatively, if $b \in \mathbf{B}(G, \mu)$ is not basic, then no cocharacter pairs with $G$ as the Levi subgroup will appear in $\mathcal{M}_{G,b,\mu}$ and so
\[
\mathrm{Mant}_{G,b,\mu}(\rho)=0.
\]
Of course, these results are already known by \cite{Shi1}, but we review them as motivation for Harris's conjecture.

We begin with the following useful definition. 

\begin{definition}
Fix $(G, \mu) \in \mathcal{C}_G$ and $b \in \mathbf{B}(G, \mu)$. Let $M_S$ be a standard Levi subgroup such that $M_S \subset M_b$. We define the subset $\mathrm{Rel}^{G, \mu}_{M_S, b} \subset \mathcal{C}_G$ as the set
\begin{align*}
   \{ (M_S, \mu_S) \in \mathcal{C}_G : \exists (M_b, \mu_b) \in \mathcal{T}_{G, b, \mu} \,\ \text{with} \,\ \theta_{M_b}(\mu_b)=\theta_{M_S}(\mu_S), \,\ \mu_b \sim_{M_b} \mu_S \}.
\end{align*}
The notation $\mu_S \sim_{M_b} \mu_b$ is defined to mean that $\mu_S$ and $\mu_b$ are conjugate in $M_b$. Note that we do not require $(M_S , \mu_S) \leq (G, \mu)$ or $(M_S, \mu_S) \leq (M_b, \mu_b)$.
\end{definition}
We record the following useful properties of $\mathrm{Rel}^{G, \mu}_{M_S, b}$.
\begin{lemma}{\label{indRel}}
We use the same notation as in the previous definition. Then
\[
\mathrm{Rel}^{G, \mu}_{M_S, b}= \coprod_{(M_b, \mu_b) \in \mathcal{I}^{G, 
\mu}_{M_b, b'}} \mathrm{Rel}^{M_b, \mu_b}_{M_S, b'}.
\]
\end{lemma}
\begin{proof}
If $(M_S, \mu_S) \in \mathrm{Rel}^{G, \mu}_{M_S, b}$, then there is an $(M_b, \mu_b) \in \mathcal{T}_{G,b, \mu}$ such that $\theta_{M_b}(\mu_b)=\theta_M(\mu_S)$ and $\mu_S \sim_{M_b} \mu_b$. Then by Proposition \ref{relTbu}, there is a unique $(M_b, \mu') \in \mathcal{I}^{G, \mu}_{M_b, b'}$ such that $(M_b, \mu_b) \in \mathcal{T}_{M_b, b', \mu'}$ and so $(M_S, \mu_S) \in \mathrm{Rel}^{M_b, \mu_b}_{M_S, b'}$. The reverse inclusion is analogous.
\end{proof}
\begin{lemma}{\label{galRel}}
The set $\mathrm{Rel}^{G, \mu}_{M_S, b}$ is invariant under the action of $W_{E_{\{\mu\}_G}}$.
\end{lemma}
\begin{proof}
If $(M_S, \mu_S) \in \mathrm{Rel}^{G, \mu}_{M_S, b}$ then we can find $(M_b, \mu_b) \in \mathcal{T}_{G, b, \mu}$ with $\theta_{M_b}(\mu_b)=\theta_{M_S}(\mu_S)$ and $\mu_b \sim_{M_b} \mu_S$. By a similar argument to Lemma \ref{galM}, we show that for each $\gamma \in W_{E_{\{\mu\}_G}}$, we have $(M_b, \gamma(\mu_b)) \in \mathcal{T}_{G, b, \mu}$ and $\theta_{M_S}(\gamma(\mu_S))=\theta_{M_b}(\gamma(\mu_b)) $ and $\gamma(\mu_S) \sim_{M_b} \gamma(\mu_b)$. This finishes the proof.
\end{proof}
Equipped with the above definition, we can now make the following restatement and slight generalization of \cite[Conj 5.4]{Har1} for the $G$ that we consider. Our statement is a generalization because we consider non-basic $b$ and do not assume the representation $I^G_{M_S}(\rho)$ is irreducible.

\begin{conjecture}[Harris]\label{Harris Conjecture}
Fix a $b \in \mathbf{B}(G, \mu)$ and a standard Levi subgroup $M_S \subset M_b$. Then for $\rho \in \mathrm{Groth}(M_S(\mathbb{Q}_p))$ a supercuspidal representation, the following representations are equal in  $\mathrm{Groth}(G(\mathbb{Q}_p) \times W_{E_{\{\mu\}_G}})$:
\[
\mathrm{Mant}_{G,b,\mu}(e(J_b)LJ( \delta^{\frac{1}{2}}_{G,P_b} \otimes I^{M_b}_{M_S}(\rho)))
\]
and
\[
[I^G_{M_S}(\rho)]\left[ \bigoplus_{(M_S, \mu_S) \in \mathrm{Rel}^{G, \mu}_{M_S, b}} r_{-\mu_S} \circ LL(\rho)|_{W_{E_{\{\mu_S\}_{M_S}}}}| \cdot|^{- \langle \rho_G, \mu \rangle}  \right].
\]
Here $r_{-\mu_S}$ is a representation of $\widehat{M_S} \rtimes W_{E_{\{\mu_S\}_{M_S}}}$ but the righthand side naturally acquires the structure of a $G(\mathbb{Q}_p) \times W_{E_{\{\mu\}_G}}$ representation from Lemma \ref{galRel} and the proof of Lemma \ref{galact}.

In particular, for $b$ basic, this says that
\[
\mathrm{Mant}_{G,b,\mu}(\mathrm{Red}_b(I^G_{M_S}(\rho))) 
=[I^G_{M_S}(\rho)] \left[ \bigoplus_{(M_S, \mu_S) \in \mathrm{Rel}^{G, \mu}_{M_S, b}} r_{-\mu_S} \circ LL(\rho)|_{W_{\{\mu_S\}_{M_S}}}| \cdot|^{- \langle \rho_G, \mu \rangle}  \right].
\]
\end{conjecture}
We will prove this conjecture assuming that Shin's formula (Theorem \ref{3.2} of this paper) holds for all admissible representations.

We proceed by induction on the rank of $T$. The key observation will be that Harris's conjecture is compatible with the Harris-Viehmann conjecture and Shin's formula. We will first assume that $I^G_{M_S}(\rho)$ is irreducible and later remove this assumption.

The following proposition shows that Conjecture \ref{Harris Conjecture} is compatible with the Harris-Viehmann conjecture (Conjecture  \ref{3.4}).

\begin{proposition} \label{4.3}
Fix $b \in \mathbf{B}(G, \mu)$ non-basic and fix a standard Levi subgroup $M_S$ of $G$ satisfying $M_S \subset M_b$. Pick $\rho \in \mathrm{Groth}(M_S(\mathbb{Q}_p))$ and suppose that $I^G_{M_S}(\rho)$ is irreducible. Suppose that Conjecture \ref{Harris Conjecture} for $\rho$ holds for $\mathrm{Mant}_{M_b, b' , \mu_b}$ for each $(M_b, \mu_b) \in \mathcal{I}^{G, \mu}_{M_b, b'}$. Then Conjecture \ref{Harris Conjecture} holds for $\mathrm{Mant}_{G, b, \mu}$.
\end{proposition}
\begin{proof}
We compute
\[
\mathrm{Mant}_{G,b, \mu}(e(J_b)LJ(\delta^{\frac{1}{2}}_{G,P_b} \otimes I^{M_b}_{M_S}(\rho)))
\]
\[
=\sum\limits_{(M_b, \mu_b) \in \mathcal{I}^{G, \mu}_{M_b, b'}} \mathrm{Ind}^G_{P_b} (\mathrm{Mant}_{M_b, b', \mu_b}(e(J_b)LJ( \delta^{\frac{1}{2}}_{G,P_b} \otimes I^{M_b}_{M_S}(\rho))))\otimes [1][|\cdot|^{\langle \rho_G, \mu_b - \mu \rangle}],
\]
so by assumption
\[
=\sum\limits_{(M_b, \mu_b) \in \mathcal{I}^{G, \mu}_{M_b, b'}} [\mathrm{Ind}^G_{P_b}(\delta^{\frac{1}{2}}_{G,P_b} \otimes I^{M_b}_{M_S}(\rho))]\left[ \bigoplus_{(M_S, \mu_S) \in \mathrm{Rel}^{M_b, \mu_b}_{M_S,b'}} r_{-\mu_S} \circ LL(I^{M_b}_{M_S}(\rho))|_{W_{E_{\{\mu_S\}_{M_S}}}}| \cdot|^{S} \right],
\]
where $S=-\langle \rho_{M_b}, \mu_b \rangle+ \langle \rho_G, \mu_b - \mu \rangle -\langle \frac{\det( Ad_{N_n}(M_b))|_T}{2}, \mu_b \rangle= - \langle \rho_G, \mu \rangle$ (following the discussion in Remark \ref{twister}). Now simplifying the above expression, we get 
\[
=\sum\limits_{(M_b, \mu_b) \in \mathcal{I}^{G, \mu}_{M_b, b'}} [I^G_{M_S}(\rho)]\left[ \bigoplus_{(M_S, \mu_S) \in \mathrm{Rel}^{M_b, \mu_b}_{M_S,b'}} r_{-\mu_S} \circ LL(I^{G}_M(\rho))|_{W_{E_{\{\mu_S\}_{M_S}}}}| \cdot|^{-\langle \rho_G, \mu \rangle} \right].
\]
Thus, we are reduced to showing that 
\[
\mathrm{Rel}^{G, \mu}_{M_S, b}= \coprod_{(M_b, \mu_b) \in \mathcal{I}^{G, 
\mu}_{M_b, b'}} \mathrm{Rel}^{M_b, \mu_b}_{M_S,b'}.
\]
This is just Lemma \ref{indRel}.
\end{proof}
With Proposition \ref{4.3} in hand, it remains to show that if Conjecture \ref{Harris Conjecture} holds for all non-basic $b \in \mathbf{B}(G, \mu)$ then it holds for the basic $b$. The key to proving this is Theorem \ref{3.2}.

We begin by making some observations about $r_{-\mu}$. Since we assumed $I^G_{M_S}(\rho)$ is irreducible, we have $LL(I^G_{M_S}(\rho))=LL(\rho)$ and the image of this representation lies inside $^LM_S \subset {}^LG$. Thus, the term $[r_{-\mu} \circ LL(I^G_{M_S}(\rho))|_{W_{E_{\{\mu\}_G}}}]$ depends only on the restriction $r_{-\mu}|_{\widehat{M_S} \rtimes W_{E_{\{\mu\}_G}}}$. Since $\mu$ is assumed to be minuscule, we have the following equality of $\widehat{M_S}$ representations. 
\begin{align}{\label{rmurest}}
r_{- \mu}|_{\widehat{M_S}}=\bigoplus_{(M_S, \mu_S) \in \mathcal{C}_G, \mu_S \sim_G \mu} r_{-\mu_S}|_{\widehat{M_S}}.
\end{align}
We further note that each $r_{-\mu_S}$ is a representation of $\widehat{M_S} \rtimes W_{E_{\{\mu_S\}_{M_S}}}$. Since $\{ (M_S, \mu_S) \in \mathcal{C}_G : \mu_S \sim_G \mu \}$ is invariant under the natural action of $W_{E_{\{\mu\}_G}}$, it follows from the proof of \ref{galact} that the right-hand side of the above equation can be promoted to a representation of $\widehat{M_S} \rtimes W_{E_{\{ \mu\}_G}}$ so that \ref{rmurest} is an equality of $W_{E_{\{\mu\}_G}}$ representations. 

Now we recall the following subsets of $W^{\mathrm{rel}}$ defined in $\S 2.11$ of \cite{Ber1}.

\begin{definition} \label{4.4}
Let $M_S, N_S$ be standard Levi subgroups of $G$. We define
\[
W^{M_S}=\{w \in W^{\mathrm{rel}}: w(M_S \cap B) \subset B \},
\]
\[
W^{M_S,N_S}=\{w \in W^{\mathrm{rel}}: w(M_S \cap B) \subset B, w^{-1}(N_S \cap B) \subset B \}
\]
\end{definition}
We record the following lemma:
\begin{lemma}{\cite[Lem 2.11]{Ber1} \label{BZleviprop}}
Suppose $M_S, N_S$ are standard Levi subgroups of $G$ and $w \in W^{M_S, N_S}$. Then $w(M_S) \cap N_S$ and $w^{-1}(N_S) \cap M_S$ are standard Levi subgroups.
\end{lemma}
\begin{lemma} \label{4.5}
Suppose $M_S$ is a standard Levi subgroup of $G$. Then $W^{M_S}$ contains a unique representative of each left coset of $W^{\mathrm{rel}}_{M_S}$. Equivalently, $(W^{M_S})^{-1}$ contains a unique representative of each right coset of $W^{\mathrm{rel}}_{M_S}$.
\end{lemma}
\begin{proof}
Suppose $w \in W^{\mathrm{rel}}$. Then $B'=w^{-1}(B)$ is a Borel subgroup of $G$ containing the maximal torus $T$. Since $B'$ contains exactly one of each root and its negative, $B' \cap M_S$ is a Borel subgroup of $M_S$. In particular, since $B' \cap M_S, B \cap M_S$ are both Borel subgroups of $M_S$ containing $T$, there exists a $w_m \in W^{\mathrm{rel}}_{M_S}$ so that 
\[
w_m(B \cap M_S)=B' \cap M_S.
\]
Then $ww_m(B \cap M_S)=B \cap M_S \subset B$, so that $ww_m \in W^{M_S}$. Thus the coset $wW^{\mathrm{rel}}_{M_S}$ contains at least one element of $W^{M_S}$.

Suppose $ww_m, ww'_m \in wW^{\mathrm{rel}}_{M_S} \cap W^{M_S}$. In particular, $ww'_m=(ww_m)(w^{-1}_mw'_m)$. But $ww_m$ takes all positive roots of $M_S$ to positive roots of $G$, and equivalently, negative roots of $M_S$ to negative roots of $G$. Thus, if $w^{-1}_mw'_m$ takes any positive root of $M_S$ to a negative root of $M_S$, then $ww'_m$ cannot be an element of $W^{M_S}$. In particular, this implies that $w^{-1}_mw'_m=1$ which shows uniqueness.
\end{proof}
\begin{lemma}{\label{wstablem}}
Suppose $M_S$ is a standard Levi subgroup of $G$ and $x \in \mathfrak{A}^+_{\mathbb{Q}, M_S}$ and $w \in W^{rel}$. Then  $w(x)=x$ if and only if $w \in W^{\mathrm{rel}}_{M_S}$.
\end{lemma}
%% this lemma and the next proposition are a little fishy a) because they might conflate W-invariants and centralizers and b) because it's not even clear what the centralizer of a cocharacter in X_*(T)_Q even means
%% Actually I think it's ok because the image of a cocharacter is connected anyway, so given a rational multiple of a cocharacter, we can always multiply it by some constant to get an honest cocharacter. Also, the centralizer of a torus (i.e. image of $\mu$) is connected and a Levi.
\begin{proof}
Recall that by assumption, $G$ is quasi-split over $\mathbb{Q}_p$ and $A$ is a split torus of $G$ of maximal rank. Pick $g \in N_G(A)(\overline{\mathbb{Q}_p})$ so that $g$ projects to $w \in W^{\mathrm{rel}}=N_G(A)(\overline{\mathbb{Q}_p})/Z_G(A)(\overline{\mathbb{Q}_p})$. Then the equation $w(x)= x$ implies that $g \in Z_G(x)(\overline{\mathbb{Q}_p})$. The centralizer of a cocharacter is a Levi subgroup, and since $x \in \mathfrak{A}^+_{\mathbb{Q}, M_S}$, we have $Z_G(x)=M_S$. In particular, $g \in N_{M_S}(A)(\overline{\mathbb{Q}_p})$ and so $w \in W^{\mathrm{rel}}_{M_S}$. 

We remark that $x$ is not a cocharacter, but that $Z_G(x)$ still makes sense as there is an induced action of $G$ on $X_*(A)_{\mathbb{Q}}$.
\end{proof}
We can now prove the following key proposition.
\begin{proposition} \label{4.7}
Fix $(G, \mu) \in \mathcal{C}_G$ and suppose $(M_S, \mu_S) \in \mathcal{C}_G$ satisfies $\mu_S \sim_G \mu$. Then there exists a unique $b \in \mathbf{B}(G, \mu)$ and a unique $w \in W^{M_S, M_b}$ so that $(w(M_S), w(\mu_S)) \in \mathrm{Rel}^{G, \mu}_{w(M_S), b}$.
\end{proposition}
\begin{proof}
We first discuss uniqueness. By assumption, $w(M_S)$ is a standard Levi subgroup. Then $w$ induces an equality $wW^{\mathrm{rel}}_{M_S}w^{-1}=W^{\mathrm{rel}}_{w(M_S)}$. In particular, $W^{\mathrm{rel}}$ acts on $X_*(T)$ through Corollary \ref{absrelweyl} and it follows that
\[
w(\theta_{M_S}(\mu_S))=\theta_{w(M_S)}(w(\mu_S)). 
\]
Since $(w(M_S), w(\mu_S)) \in \mathrm{Rel}^{G, \mu}_{w(M_S), b}$, it follows that $\theta_{w(M_S)}(w(\mu_S))$ is dominant in the relative root system. In particular, $\theta_{w(M_S)}(w(\mu_S))$ must be equal to the unique element $x$ in the $W^{\mathrm{rel}}$ orbit of $\theta_{M_S}(\mu_S)$ which is dominant in $\mathfrak{A}_{\mathbb{Q}}$. Now $x \in \mathfrak{A}^+_{M_{S'}, \mathbb{Q}}$ for a unique $M_{S'}$. Since any $(M_b, \mu_b) \in \mathcal{T}_{G, b, \mu}$ is definitionally strictly decreasing, it follows that even though we can't yet conclude the uniqueness of $b$, we have shown that any other $b_1$ must satisfy $M_{b_1}=M_b=M_{S'}$.

Now, suppose we had $w, w' \in W^{M_S, M_b}$ such that 
\[
w(\theta_{M_S}(\mu_S))=x=w'(\theta_{M_S}(\mu_S)).
\]
 Then in particular, $w'w^{-1}$ stabilizes $x$ and so by Lemma \ref{wstablem}, $w'w^{-1} \in W^{\mathrm{rel}}_{M_b}$. So $w$ and $w'$ are in the same right coset $W^{\mathrm{rel}}_{M_b}w$. However, $ W^{M_S, M_b} \subset (W^{M_b})^{-1}$. By Lemma \ref{4.5}, $(W^{M_b})^{-1}$ contains a unique representative of each right coset of $(W^{M_b})^{-1}$ and so there is a unique $w \in (W^{M_b})^{-1}$ satisfying $w(\theta_{M_S}(\mu_S))=x$. In particular, this implies that $w=w'$. Thus, we have shown that $w$ is unique, if it exists. There is exactly one $\mu' \in X_*(T)$ such that $\mu' \sim_{M_b} w(\mu)$ and $\mu'$ is dominant in $M_b$. Then $(M_b, \mu') \in \mathcal{T}_{G, b, \mu}$ for at most one $b \in \mathbf{B}(G, \mu)$. This shows uniqueness.

To prove existence, we again define $x$ to be the unique dominant element in the $W^{\mathrm{rel}}$-orbit of $\theta_{M_S}(\mu_S)$. Define $M_{S'}=Z_G(x)$ and take the unique $w \in (W^{M_{S'}})^{-1}$ such that $w(\theta_{M_S}(\mu_{M_S}))=x$.  We would like to show that $w \in W^{M_S, M_{S'}}$. 

By definition,
\[
w(M_S) \subset w(Z_G(\theta_{M_S}(\mu_S))) =Z_G(x)=M_{S'}.
\]
Suppose it is not the case that $w(M_S \cap B) \subset B$. In particular, $w$ maps a positive root $r$ of $M_S$ to a root $w(r)$ of $M_{S'}$ which is not positive. In particular, $-w(r)$ is positive and so $w^{-1}(-w(r))=-r$ is positive (since $w \in (W^{M_{S'}})^{-1}$). But this is clearly a contradiction. Thus, in fact $w \in W^{M_S, M_{S'}}$.

By Lemma \ref{BZleviprop}, $w(M_S) \cap M_{S'}=w(M_S)$ is a standard Levi. It remains to show that $(w(M_S), w(\mu_S))$ is a cocharacter pair and an element of $\mathrm{Rel}^{G, \mu}_{w(M_S),b}$. Now if $r$ is a positive root in the absolute root system of $w(M_{S})$, then $\langle r, w(\mu_S) \rangle= \langle w^{-1}(r), \mu_S \rangle  \geq 0$ (since $(M_S, \mu_S)$ is a cocharacter pair and $w^{-1}(r)$ is a positive root of $M_S$). Thus, $(w(M_S), w(\mu_S))$ is a cocharacter pair. By construction, $x=\theta_{w(M_S)}(w(\mu_S))=\theta_{M_{S'}}(w(\mu_S))$. Suppose $\mu' \in X_*(T)$ is the unique cocharacter conjugate to $w(\mu_S)$ in $M_{S'}$ and dominant in $M_{S'}$. Then by Corollary \ref{conjinv}, $(M_{S'}, \mu')$ is strictly decreasing and therefore $(M_{S'}, \mu') \in \mathcal{T}_{G,b,\mu}$ for some $b$ and so $(w(M_S), w(\mu_S)) \in \mathrm{Rel}^{G, \mu}_{w(M_S),b}$.
\end{proof}
\begin{corollary}{\label{sumrel}}
Fix a cocharacter pair $(G, \mu) \in \mathcal{C}_G$ and a standard Levi subgroup $M_S$ of $G$. For $b \in \mathbf{B}(G, \mu)$, define $W_b$ by $\{w \in W^{M_S, M_b} : w(M_S) \subset M_b \}$. Then the previous lemma gives a bijection
\[
\{(M_S, \mu_S) \in \mathcal{C}_G: \mu_S \sim_G \mu\} \cong \coprod\limits_{b \in \mathbf{B}(G, \mu)} \coprod\limits_{w \in W_b} \mathrm{Rel}^{G, \mu}_{w(M_S),b}.
\]
\end{corollary}
\begin{proof}
By the construction in the previous proposition, it is clear that given an $(M_S, \mu_S) \in \mathcal{C}_G$ we get an element of the right-hand side of the above equation. Conversely, an element $(w(M_S), \mu')$ of the right-hand side comes with a fixed $w \in W_b$ and so we can recover $(M_S, w^{-1}(\mu'))$ on the left-hand side.
\end{proof}

We are now ready to finish the proof of Conjecture \ref{Harris Conjecture}. By inductive assumption we assume we've shown Conjecture \ref{Harris Conjecture} for $G$ with maximal torus of rank less than $n$. Then Proposition \ref{4.3} implies that Conjecture \ref{Harris Conjecture} holds for $G$ with maximal torus of rank $n$ in the case where $b$ is not basic. It remains to prove the basic case, for which it suffices to show that Theorem \ref{3.2} is compatible with Conjecture \ref{Harris Conjecture}.
We have
\[
\sum\limits_{b \in \mathbf{B}(G, \mu)} \mathrm{Mant}_{G, b, \mu}(\mathrm{Red}_b(I^G_{M_S}(\rho)))
\]
\[
= \sum\limits_{b \in \mathbf{B}(G, \mu) } \mathrm{Mant}_{G, b, \mu}(e(J_b)LJ(\delta^{\frac{1}{2}}_{P_b} \otimes J^G_{P^{op}_b} I^G_{M_S}( \rho))).
\]
By the geometric lemma of \cite{Ber1} and noting that $W^{M_S, M_b}$ defined with respect to $B$ is equal to the analogous set defined with respect to $B^{op}$, we have
\[
J^G_{P^{op}_b} I^G_{M_S}( \rho)=\sum\limits_{w \in W^{M_S, M_b}} I^{M_b}_{M'_b}(w(J^{M_S}_{{P'}^{op}_S}(\rho))),
\]
where $M'_S=M_S \cap w^{-1}(M_b), M'_b=w(M_S) \cap M_b$. By the assumption that $\rho$ is supercuspidal we must have $M'_S=M_S$ and $M'_b=w(M_S)$. In this case, we have from the geometric lemma that $w(M_S)$ is a standard Levi subgroup. Thus we get that the previous expression is equal to
\[
\sum\limits_{b \in \mathbf{B}(G, \mu) } \mathrm{Mant}_{G, b, \mu}(e(J_b)\sum\limits_{w \in W_b} LJ(\delta^{\frac{1}{2}}_{P_b} \otimes I^{M_b}_{w(M_S)}( w(\rho))),
\]
where $W_b \subset W^{M_S, M_b}$ is the subset of $w$ such that $w(M_S) \subset M_b$. We now apply Corollary \ref{Harris Conjecture} by inductive assumption to get
\[
\sum\limits_{b \in \mathbf{B}(G, \mu)} \sum\limits_{w \in W_b} [I^G_{w(M_S)}(w(\rho))]\left[ \bigoplus_{(w(M_S), \mu') \in \mathrm{Rel}^{G, \mu}_{w(M_S), b}} r_{-\mu'} \circ LL(I^G_{w(M_S)}(w(\rho)))|_{W_{E_{\{\mu'\}_{w(M_S)}}}}| \cdot |^{-\langle \rho_G, \mu \rangle} \right].
\]
By \cite[Thm 2.9]{Ber1}, we have that 
\[
[I^G_{w(M_S)}(w(\rho))]=[I^G_{M_S}(\rho)],
\]
and since $I^G_{M_S}(\rho)$ is assumed to be irreducible, we have 
\[
LL(I^G_{M_S}(\rho))=LL(\rho).
\]
Finally, we note that $W_{E_{\{w^{-1}(\mu')\}_{M_S}}}=W_{E_{\{\mu'\}_{w(M_S)}}}$ and we have an equality
\[
[r_{- \mu'} \circ LL(w(\rho))|_{W_{E_{\{\mu'\}_{w(M_S)}}}}]=[r_{-w^{-1} (\mu')} \circ LL(\rho)|_{W_{E_{\{w^{-1}(\mu')\}_{M_S}}}}].
\]

Thus the above expression becomes
\[
\sum\limits_{b \in \mathbf{B}(G, \mu)} \sum\limits_{w \in W_b} [I^G_{M_S}(\rho)]\left[ \bigoplus_{(w(M_S), \mu') \in \mathrm{Rel}^{G, \mu}_{w(M_S), b}} r_{-w^{-1}(\mu')} \circ LL(\rho)|_{W_{E_{\{w^{-1}(\mu')\}_{M_S}}}}| \cdot |^{-\langle \rho_G, \mu \rangle}\right].
\]
By Corollary \ref{sumrel} this equals
\[
[I^G_{M_S}(\rho)][ \bigoplus_{(M_S, \mu_S): \mu_S \sim_G \mu} r_{-\mu_S} \circ LL(\rho)|_{W_{E_{\{\mu_S\}_{M_S}}}}| \cdot |^{-\langle \rho_G, \mu \rangle}].
\]
Finally, we apply the decomposition given by Equation \eqref{rmurest} to get
\[
[I^G_{M_S}(\rho)][r_{-\mu}|_{\widehat{M_S} \rtimes W_{E_{\{\mu\}_G}}} \circ LL(\rho)|_{W_{E_{\{\mu\}_{G}}}} | \cdot |^{- \langle \rho_G, \mu \rangle}],
\]
which is the desired result.

Finally, we show that Conjecture \ref{Harris Conjecture} holds even if $I^G_{M_S}(\rho)$ is not irreducible. Our verification that Conjecture \ref{Harris Conjecture} is compatible with the Harris-Viehmann conjecture did not rely on the irreducibility of $I^G_{M_S}(\rho)$. Thus in the case where we do not assume $I^G_{M_S}(\rho)$ is irreducible, it would suffice to show that Conjecture \ref{Harris Conjecture} is true in the case where $b$ is basic. If $b$ is basic, then $M_b=G$ so we have
\[
\mathrm{Mant}_{G,b,\mu}(e(J_b)LJ( \delta^{\frac{1}{2}}_{G,P_b}I^{M_b}_{M_S}(\rho)))=\mathrm{Mant}_{G,b, \mu}(\mathrm{Red}_b(I^G_{M_S}(\rho))).
\]
This can now be computed by cocharacter pairs using the results of \S3. If $I^G_{M_S}(\rho)$ is assumed to be irreducible, then for each cocharacter pair $(M_{S'}, \mu_{S'})$ of $G$, we have
\[
[M_{S'}, \mu_{S'}](I^G_{M_S}(\rho))=(\mathrm{Ind}^G_{P_{S'}} \circ [\mu_{S'}] )( \delta^{\frac{1}{2}}_{P_S} \otimes J^G_{P^{op}_{S'}}I^G_{M_S}(\rho))\otimes [1][|\cdot|^{\langle \rho_G, \mu_{S'}-\mu \rangle}]
\]
\[
=(\mathrm{Ind}^G_{P_{S'}} \circ [\mu_{S'}] )(\bigoplus_{w \in W_{\rho}} \delta^{\frac{1}{2}}_{P_{S'}} \otimes I^{M_{S'}}_{w(M_S)} (w(\rho)))\otimes [1][|\cdot|^{\langle \rho_G, \mu_{S'}-\mu \rangle}],
\]
where $W_{\rho}$ is the subset of $w \in W^{M_S, M_{S'}}$ such that $w(M_S) \subset M_{S'}$. Then the above equals
\[
[I^G_{M_S}(\rho)]\left[ \bigoplus_{w \in W_{\rho}} r_{-\mu_{S'}} \circ LL(w(\rho))|\cdot|^{- \langle \rho_G, \mu \rangle } \right].
\]
Thus we see that applying various $[M_{S'}, \mu_{S'}]$ to $I^G_{M_S}(\rho)$ in the irreducible case will always yield the same term of $\mathrm{Groth}(G(\mathbb{Q}_p))$ (namely $[I^G_{M_S}(\rho)]$) and so when evaluating $\mathrm{Mant}_{G,b,\mu}(\mathrm{Red}_b(I^G_{M_S}(\rho))$ as a sum of cocharacter pairs, the different Galois terms must cancel to give Conjecture \ref{Harris Conjecture}. Thus, if we can show that in the reducible case, the $\mathrm{Groth}(G(\mathbb{Q}_p))$ part of each $[M_{S'}, \mu_{S'}](I^G_{M_S}(\rho))$ is fixed and the Galois part is identical to the irreducible case, then Conjecture \ref{Harris Conjecture} must hold for this case as well.

The first part of our previous computation did not depend on the irreducibility of $I^G_{M_S}(\rho)$ so we still have
\[
[M_{S'}, \mu_{S'}](I^G_{M_S}(\rho))=(\mathrm{Ind}^G_{P_{S'}} \circ [\mu_{S'}] )(\bigoplus_{w \in W_{\rho}} \delta^{\frac{1}{2}}_{P_{S'}} \otimes I^{M_{S'}}_{w(M_S)} (w(\rho)))\otimes [1][|\cdot|^{\langle \rho_G, \mu_{S'}-\mu \rangle}].
\]
Suppose now that $I^{M_{S'}}_{w(M_S)}(w(\rho))=\pi_1 \oplus ... \oplus \pi_k$. Then using that for all $i$, we have $LL(\pi_i)=LL(w(\rho))$,
\[
[\mu_{S'}](I^{M_{S'}}_{w(M_S)}(w(\rho)))=\oplus^k_{i=1}[\pi_i][r_{-\mu_{S'}} \circ LL(\pi_i)\otimes | \cdot |^{ - \langle \rho_{M_{S'}}, \mu_{S'} \rangle}]
\]
\[
=\oplus^k_{i=1}[\pi_i][r_{-\mu_{S'}} \circ LL(w(\rho))\otimes | \cdot |^{ - \langle \rho_{M_{S'}}, \mu_{S'} \rangle}]
\]
\[
=[I^{M_{S'}}_{w(M_S)}(w(\rho))][r_{-\mu_{S'}} \circ LL(w(\rho))\otimes | \cdot |^{ - \langle \rho_{M_{S'}}, \mu_{S'} \rangle}]
\]
Thus, the expression for $[M_{S'}, \mu_{S'}](I^G_{M_S}(\rho))$ becomes 
\[
[I^G_{M_S}(\rho)] \left[ \bigoplus_{w \in W^{M_S, M_{S'}}} r_{-\mu_{S'}} \circ LL( w(\rho))| \cdot |^{ - \langle \rho_G, \mu \rangle} \right],
\]
as desired.

\appendix

 \section{Examples}{\label{A}} \label{examples}
In this section, we give an example to show that even in the unramified EL-type case, we do not get an expression as simple as Harris's conjecture for $\mathrm{Mant}_{G,b, \mu}(\rho)$ for general $\rho$. We generally use the same notation as in the computation in Example \ref{HVex}.

Let $G=\mathrm{GL}_4$, suppose $\mu$ has weights $(1^2,0^2)$, and take $b$ basic. Let $T$ be the diagonal maximal torus and $B$ be the Borel subgroup of upper triangular matrices. Then the set of cocharacter pairs less than or equal to $(G, \mu)$ is as follows.

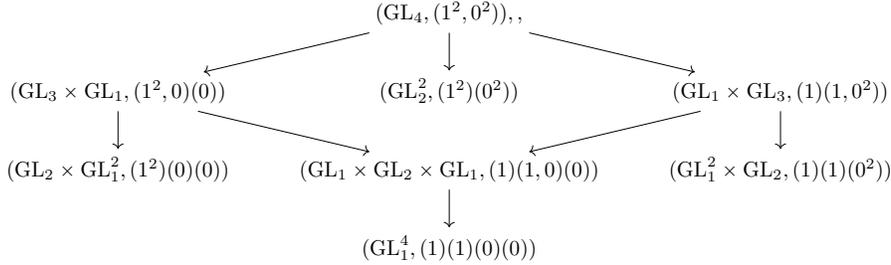
\begin{adjustbox}{max size={.95\textwidth}{.8\textheight}}
\begin{tikzcd}
& (\mathrm{GL}_4, (1^2,0^2)) \arrow[ld] \arrow[d] \arrow[rd] & \\
(\mathrm{GL}_3 \times \mathrm{GL}_1, (1^2,0)(0)) \arrow[d] \arrow[rd] & (\mathrm{GL}^2_2, (1^2)(0^2)) & (\mathrm{GL}_1 \times \mathrm{GL}_3, (1)(1,0^2)) \arrow[ld] \arrow[d] \\
(\mathrm{GL}_2 \times \mathrm{GL}^2_1, (1^2)(0)(0) ) & (\mathrm{GL}_1 \times \mathrm{GL}_2 \times \mathrm{GL}_1 , (1)(1,0)(0)) \arrow[d] & (\mathrm{GL}^2_1 \times \mathrm{GL}_2, (1)(1)(0^2))\\
& (\mathrm{GL}^4_1, (1)(1)(0)(0) ) \\
\end{tikzcd}
\end{adjustbox}\\
Let $\rho \in \mathrm{Groth}(\mathrm{GL}_1(\mathbb{Q}_p))$ and consider $\pi$ the unique essentially square integrable quotient of $I^G_{\mathrm{GL}^4_1}(\rho \boxtimes \rho(1) \boxtimes \rho(2) \boxtimes \rho(3))$. We want to compute $\mathrm{Mant}_{G,b,\mu}(\mathrm{Red}_b(\pi))$.

We introduce some notation which will allow us to describe the answer to this question. The results of $\S 2$ of  \cite{Zel1} show that  $I^G_{\mathrm{GL}^4_1}(\rho \boxtimes \rho(1) \boxtimes \rho(2) \boxtimes \rho(3))$ has exactly $8$ irreducible subquotients. If $\pi'$ is one such subquotient, then $J^G_{B^{op}}(\pi')$ will be a finite sum of representations of the form $\rho(\lambda(0))\boxtimes \rho(\lambda(1)) \boxtimes \rho(\lambda(2)) \boxtimes \rho(\lambda(3))$ where $\lambda$ is a permutation of $\{0 , 1 ,2, 3 \}$. In particular, if $\Omega$ denotes the set of all such permutations of $\rho \boxtimes \rho(1) \boxtimes \rho(2) \boxtimes \rho(3)$, then each permutation lies in the Jacquet module of exactly one irreducible subquotient of $I^G_{\mathrm{GL}^4_1}(\rho \boxtimes \rho(1) \boxtimes \rho(2) \boxtimes \rho(3))$ so that the irreducible subquotients correspond to a partition of $\Omega$. We use the following shorthand: we define the notation $(0123)$ to refer to the representation $\rho(0)\boxtimes \rho(1) \boxtimes \rho(2) \boxtimes \rho(3)$.  Following Zelevinsky, our $8$ irreducible subquotients naturally correspond to vertices of a 3-dimensional cube, and so we denote them by binary strings of length $3$. Then if we denote the subset of $\Omega$ corresponding to some subquotient $\pi'$ by $\Omega(\pi')$,we have
\begin{align*}
\Omega([000])&=\{(3210)\} \\
\Omega([100])&=\{(2310), (2130), (2103) \} \\
\Omega([010])&=\{ (3120), (1320), (1302), (3102), (1032)\} \\
\Omega([001])&=\{(3201), (3021), (0321)\} \\
\Omega([110])&=\{(1203), (1023), (1230) \} \\
\Omega([101])&=\{(2013), (2031), (0213), (0231), (2301) \} \\
\Omega([011])&=\{(3012), (0312), (0132) \} \\
\Omega([111])&=\{(0123)\}
\end{align*}
In particular, our representation $\pi$ corresponds to $[111]$ under the above notation. A tedious computation using Theorem \ref{3.8} yields the following
\begin{proposition}
\begin{align*}
\mathrm{Mant}_{G,b,\mu}(\mathrm{Red}_b(\pi))
&=[111][\widecheck{LL(\rho)}^2(-7)+\widecheck{LL(\rho)}^2(-6)] \\
&-([110][\widecheck{LL(\rho)}^2(-5)]+[011][\widecheck{LL(\rho)}^2(-5)])\\
&+[010][\widecheck{LL(\rho)}^2(-4)]\\
&-[000][\widecheck{LL(\rho)}^2(-3)]
\end{align*}
\end{proposition}
We finish by remarking that the set of cocharacter pairs less than or equal to $(G, \mu)$ has some special properties in the above case that make the general case more complicated. 

For instance, each $\mathcal{T}_{G,b, \mu}$ has at most a single element. However, if $G$ has a nontrivial action by $\Gamma$, this need not be the case.

In the case we consider, we have a single cocharacter pair for each Levi subgroup. In general, this need not be the case. For instance, if $G=\mathrm{GL}_5, \mu=(1^3,0^2)$, then $(\mathrm{GL}_3 \times \mathrm{GL}_2, (1^3)(0^2)), (\mathrm{GL}_3\times \mathrm{GL}_2, (1^2,0)(1,0))$ are both less than $(G, \mu)$.

Further, in the above example, each cocharacter pair $(M_S, \mu_S)$ had the property that $\mu_S$ was dominant as a cocharacter of $G$ relative to $B$. In general this need not be the case. In fact, $(\mathrm{GL}^5_1, (1)(1)(0)(1)(0)) \leq (\mathrm{GL}_5, (1^3,0^2))$.

\section{Relative Root Systems and Weyl Chambers}
In this section we prove a fact about root systems that is needed in the text (for instance in the proof of Proposition \ref{2.16}). We assume that $G$ is a quasisplit group over a field $k$ of characteristic $0$ and pick a separable closure $k^{sep}$. We fix a split $k$-torus $A$ of maximal rank in $G$ and choose a maximal torus $T$ and Borel subgroup $B$ both defined over $k$ and such that $A \subset T \subset B$. Associated to this data, we have an absolute root datum 
\[
(X^*(T), \Phi^*(G,T), X_*(T), \Phi_*(G,T) ),
\]
and a relative root datum 
\[
(X^*(A), \Phi^*(G, A), X_*(A), \Phi_*(G, A)).
\]
Our choice of $B$ also gives sets $\Delta$ of absolute simple roots and $_k\Delta$ of relative simple roots. Note that we also have a natural restriction map 
\[
\mathrm{res}: X^*(T) \to X^*(A),
\]
and that by definition an absolute root in $\Phi^*(G,T)$ restricts to an element of $\Phi^*(G,A) \cup \{0\}$.

We record two standard consequences of our assumption that $G$ is quasisplit.
\begin{proposition}{\label{rootrestrict}}
Let $G$ be quasisplit and use the notations as above. Then,
\begin{enumerate}
    \item The centralizer $Z_G(A)=T$,
    \item We have $\mathrm{res}(\Delta)= \,\ _k\Delta$. The key point being that no absolute simple root restricts to the trivial character.
\end{enumerate}
\end{proposition}
We have the following easy consequence on the structure of the Weyl group of the relative root system. Recall that the absolute Weyl group $W $equals
\[
N_G(T)(k^{sep})/Z_G(T)(k^{sep}),
\]
and the relative Weyl group $W^{\mathrm{rel}}$ is $N_G(A)(k)/Z_G(A)(k)$.
\begin{corollary}{\label{absrelweyl}}
We have the following equality: $W^{\mathrm{rel}}=W^{\Gamma}$, where $\Gamma= \mathrm{Gal}(k^{sep}/ k)$.
\end{corollary}
\begin{proof}
It suffices to show that $Z_G(A)=Z_G(T)$ and that $N_G(A)(k)=N_G(T)(k)$. For the first equality, we note that by the quasisplit assumption, $Z_G(A)=T=Z_G(T)$. For the second equality, we note that any $g \in N_G(A)(k)$ must also normalize the centralizer of $A$ which is $T$. Conversely, if $g \in N_G(T)(k)$ then $g$ normalizes the unique maximal $k$-split sub-torus of $T$ which is $A$.
\end{proof}

Define the absolute Weyl chamber $\overline{C}^*_{\mathbb{Q}} \subset X^*(T)_{\mathbb{Q}}$ by $\{ x \in X^*(T)_{\mathbb{Q}}: \langle \check{\alpha}, x \rangle \geq 0, \alpha \in \Delta\}$ and define the relative Weyl chamber $_k\overline{C}^*_{\mathbb{Q}} \subset X^*(A)_{\mathbb{Q}}$ analogously. The key result of this section is that 
\[
\mathrm{res}(\overline{C}^*_{\mathbb{Q}}) \,\ = \,\ _k\overline{C}^*_{\mathbb{Q}}.
\]
Despite its simple statement, the author has been unable to locate a convenient reference of this fact. For $x \in X^*(T)_{\mathbb{Q}}$ and $\alpha \in \Delta$, we need to relate $\langle \widecheck{\alpha}, x \rangle$ and $\langle \widecheck{\mathrm{res}(\alpha)}, \mathrm{res}(x) \rangle$. If we let $\sigma_{\alpha} \in W$ be the reflection corresponding to the root $\alpha$, then we have
\begin{align}{\label{pairw}}
x- \sigma_{\alpha}(x)=\langle \widecheck{\alpha}, x \rangle \alpha.
\end{align}
and analogously for $\widecheck{\mathrm{res}(\alpha)}$. Thus it will suffice to relate $\sigma_{\alpha}$ and $\sigma_{\mathrm{res}(\alpha)}$.

Note that since $B$ is defined over $k$, we have $\gamma(\Delta)=\Delta$ for every $\gamma \in \Gamma$.  Moreover, for each $\alpha \in \Delta$, we have $\mathrm{res}(\gamma(\alpha))=\mathrm{res}(\alpha)$. After all, $\Gamma$ acts trivially on $X^*(A)_{\mathbb{Q}}$ and the restriction map is $\Gamma$-equivariant.

Now fix $\alpha \in \Delta$  and let $W_{\alpha}$ be the subgroup of $W$ generated by the elements $\sigma_{\gamma(\alpha)}$ for each $\gamma \in \Gamma$. We claim that if we can find a nontrivial $\Gamma$-invariant element of $W_{\alpha}$, then it must equal $\sigma_{\mathrm{res}(\alpha)}$. To prove this, we first recall the construction of $\sigma_{\alpha}$ and $\sigma_{\mathrm{res}(\alpha)}$ (see \cite[pg 230]{Bor1}) for instance). Given a root $\alpha \in \Phi^*(G,T)$ we can define a group $G_{\alpha}=Z_G(T_{\alpha})$ where $T_{\alpha}= \mathrm{ker}(\alpha)^0 \subset T$. Then $N_{G_{\alpha}}(T)(k^{sep})/Z_{G_{\alpha}}(T)(k^{sep})$ embeds into $W$ and has a unique nontrivial element which is $\sigma_{\alpha}$. Analogously, we define $A_{\mathrm{res}(\alpha)}$ and $G_{\mathrm{res}(\alpha)} =Z_G(A_{\mathrm{res}(\alpha)})$. Then  $N_{G_{\mathrm{res}(\alpha)}}(A)(k)/Z_{G_{\mathrm{res}(\alpha)}}(A)(k)$ embeds into $W^{\mathrm{rel}}$ and has a unique nontrivial element that is identified with $\sigma_{\mathrm{res}(\alpha)}$. 

Now, by Corollary \ref{absrelweyl} we have
\[
N_{G_{\mathrm{res}(\alpha)}}(A)(k)/Z_{G_{\mathrm{res}(\alpha)}}(A)(k)=N_{G_{\mathrm{res}(\alpha)}}(T)(k)/Z_{G_{\mathrm{res}(\alpha)}}(T)(k).
\]
Thus to complete the proof of the claim, we need to show that 
\begin{align}{\label{weyleqn}}
N_{G_{\alpha}}(T)(k^{sep})/Z_{G_{\alpha}}(T)(k^{sep}) \hookrightarrow N_{G_{\mathrm{res}(\alpha)}}(T)(k^{sep})/Z_{G_{\mathrm{res}(\alpha)}}(T)(k^{sep}).
\end{align}
After all, the unique nontrivial $\Gamma$-invariant element of the group on the right is $\sigma_{\mathrm{res}(\alpha)}$ and the group on the left contains $\sigma_{\alpha}$. Since we get the same equation if we replace $\alpha$ everywhere with $\gamma(\alpha)$, this will imply that 
\[
W_{\alpha} \subset N_{G_{\mathrm{res}(\alpha)}}(T)(k^{sep})/Z_{G_{\mathrm{res}}}(T)(k^{sep}).
\]
Now, Equation \eqref{weyleqn} follows from the fact that
\[
Z_{G_{\alpha}}(T)=Z_{G_{\mathrm{res}(\alpha)}}(T)=T
\]
and 
\[
N_{G_{\alpha}}(T) \subset N_{G_{\mathrm{res}(\alpha)}}(T).
\]

We are now interested in finding a nontrivial $\Gamma$-invariant element of the group $W_{\alpha}$ defined above. In fact, $W_{\alpha}$ will be a finite Coxeter group and the element we seek is the unique element of longest length. We need to compute this element explicitly, which we now do. We treat two cases. Suppose first that the elements of the $\Gamma$-orbit of $\sigma_{\alpha}$ commute pairwise. Then clearly the product $\prod\limits_{\gamma \in \Gamma/ \mathrm{stab}(\sigma_{\alpha})} \sigma_{\gamma(\alpha)}$ is $\Gamma$-invariant.

In the second case, suppose that the $\Gamma$-orbit of $\sigma_{\alpha}$ has precisely two elements which we denote $X$ and $Y$. Then we have $(XY)^k=1$ for some $k \geq 2$ which we assume to be minimal. If $k$ is even, then $(XY)^{k/2}$ is invariant and nontrivial and if $k$ is odd, then $Y(XY)^{(k-1)/2}$ is invariant and nontrivial.

We now prove that any $\Gamma$ action on the simple roots $\Delta$ of $G$ is a combination of these cases. The action of $\Gamma$ on $\Delta$ induces an action on the associated (not necessarily connected) Dynkin diagram $D$. Each $\gamma \in \Gamma$ maps connected components of $D$ to connected components and so there is an induced action of $\Gamma$ on the set of connected components $\pi_0(D)$.

Now fix an $\alpha \in \Delta$ and consider the $\Gamma$-orbit $\Gamma \alpha$ of $\alpha$. Suppose $D^i$ is a connected component of $D$ such that $D^i \cap \Gamma \alpha \neq \emptyset$. Then via the classification of connected Dynkin diagrams, we see that $\Gamma \alpha \cap D^i$ contains either a single node, $2$ non-adjacent nodes, $2$ adjacent nodes, or $3$ nodes where no two are adjacent. In particular, these are all covered by the cases we considered above, so we can find an element $w_i$ of $W_{\alpha}$ that is invariant by the action of $\mathrm{stab}(D^i) \subset \Gamma$. Then $\Gamma \alpha$ consists of finitely many disjoint copies of one of the above possibilities and so we see that $\prod\limits_i w_i$ is $\Gamma$-invariant and an element of $W_{\alpha}$ and therefore equal to $\sigma_{\mathrm{res}(\alpha)}$. Equipped with this description, we now give a proof of the main result of this section.
\begin{proposition}{\label{chamberresprop}}
We continue to observe the assumptions made above. In particular, $G$ is a quasisplit group over $k$. Then the map $\mathrm{res}: X^*(T) \twoheadrightarrow X^*(A)$ induces an equality 
\[
\mathrm{res}(\overline{C}^*_{\mathbb{Q}})=  \,\ _k\overline{C}^*_{\mathbb{Q}}.
\]
\end{proposition}
\begin{proof}
We first show that $\mathrm{res}(\overline{C}^*_{\mathbb{Q}}) \subset \,\ _k\overline{C}^*_{\mathbb{Q}}$. Pick $x \in \overline{C}^*_{\mathbb{Q}}$ and $\alpha \in \Delta$. Then we need to show that 
\[
\langle \widecheck{\mathrm{res}(\alpha)}, \mathrm{res}(x) \rangle \geq 0
\]
or equivalently, that
\[
\mathrm{res}(x)- \sigma_{\mathrm{res}(\alpha)}(\mathrm{res}(x))
\]
is a non-negative multiple of $\mathrm{res}(\alpha)$. Note that $\mathrm{res}$ is $W^{\Gamma}$-equivariant (where $W^{\Gamma}$ acts as $W^{\mathrm{res}}$ on $X^*(A)$). Thus, it suffices to show that
\[
\mathrm{res}(x- \sigma_{\mathrm{res}(\alpha)}(x))
\]
is a non-negative multiple of $\mathrm{res}(\alpha)$. Thus, we need to compute $x- \sigma_{\mathrm{res}(\alpha)}(x)$. We do so using our description of $\sigma_{\mathrm{res}(\alpha)}$. 

We first consider the case where the $\Gamma$-orbit of $\sigma_{\alpha}$ consists of pairwise commuting elements. Equivalently, the elements of $\Gamma \alpha$ are pairwise orthogonal. Then 
\[
\sigma_{\mathrm{res}(\alpha)}=\sigma_{\alpha_n} \circ ... \circ \sigma_{\alpha_1}
\]
for $\{\alpha_1, ..., \alpha_n\}=\Gamma\alpha$. Since $x$ is dominant in the absolute root system, we have
\[
x- \sigma_{\alpha_i}(x)= a_i \alpha_i
\]
for some $a_i \geq 0$. Then since $\alpha_i$ is orthogonal to $\alpha_j$ for $i \neq j$, we have $\sigma_{\alpha_i}(\alpha_j)=\alpha_j$. Thus,
\begin{align*}
x-\sigma_{\mathrm{res}(\alpha)}(x)&= \sum\limits^n_{i=1} (\sigma_{\alpha_1} \circ ... \circ \sigma_{\alpha_{i-1}})(x) - (\sigma_{\alpha_1} \circ ... \circ \sigma_{\alpha_i})(x)\\
&=\sum\limits^n_{i=1} (\sigma_{\alpha_1} \circ ... \circ \sigma_{\alpha_{i-1}})(x- \sigma_{\alpha_i}(x))\\
&=\sum\limits^n_{i=1} (\sigma_{\alpha_1} \circ ... \circ \sigma_{\alpha_{i-1}})(a_i \alpha_i)\\
&=\sum\limits^n_{i=1} a_i \alpha_i.
\end{align*}
Thus in this case,
\[
\mathrm{res}(x- \sigma_{\mathrm{res}(\alpha)}(x))=(a_1+...+a_n) \mathrm{res}(\alpha)
\]
and $a_1+...+a_n \geq 0$ as desired.

Now we consider the case where $\Gamma \alpha=\{ \alpha, \beta \}$ and $\alpha$ and $\beta$ are adjacent in $D$ and connected by a single edge. Then $\sigma_{\alpha}(\beta)=\alpha+\beta= \sigma_{\beta}(\alpha)$. In this case, $\sigma_{\mathrm{res}(\alpha)}=\sigma_{\beta} \circ \sigma_{\alpha} \circ \sigma_{\beta}$. By assumption, we have that $x- \sigma_{\alpha}(x)=a\alpha$ and $x- \sigma_{\beta}(x)=b \beta$ for $a$ and $b$ non-negative. Thus,
\begin{align*}
    x-\sigma_{\mathrm{res}(\alpha)}(x) &= (x- \sigma_{\beta}(x))+ \sigma_{\beta}(x- \sigma_{\alpha}(x)) + (\sigma_{\beta} \circ \sigma_{\alpha})(x- \sigma_{\beta}(x))\\
    &=b\beta+a(\alpha+\beta)+b\alpha\\
    &=(a+b)(\alpha+\beta),
\end{align*}
which projects to $2(a+b) \mathrm{res}(\alpha)$ and $2(a+b) \geq 0$ as desired.

Finally, we must consider the case where $\Gamma \alpha$ equals $\{ \alpha_1, \beta_1, ..., \alpha_n, \beta_n\}$ such that $\alpha_i$ and $\beta_i$ are connected by a single edge in $D$ but for $i \neq j$, neither $\alpha_i$ nor $\beta_i$  are connected to either $\alpha_j$ or $\beta_j$. We compute $x- (\sigma_{\beta_i} \circ \sigma_{\alpha_i} \circ \sigma_{\beta_i})(x)$ as in the previous paragraph. Then if we let $w_i=\sigma_{\beta_i} \circ \sigma_{\alpha_i} \circ \sigma_{\beta_i}$, we have 
\[
\sigma_{\mathrm{res}(\alpha)}=w_1 \circ .. \circ w_n.
\]
Now we can compute $x- \sigma_{\mathrm{res}(\alpha)}(x)$ as in the commuting case, substituting $w_i$ for $\sigma_{\alpha_i}$. We see in this case that 
\[
\mathrm{res}(x- \sigma_{\mathrm{res}(\alpha)}(x))=2(a_1+b_1+...+a_n+b_n)\mathrm{res}(\alpha).
\]
This concludes the proof that $\mathrm{res}(\overline{C}^*_{\mathbb{Q}}) \subset \,\ _k\overline{C}^*_{\mathbb{Q}}$.

It remains to show that we actually have equality. We claim it suffices to show that the fundamental weight $\delta_{\mathrm{res}(\alpha)}$ is an element of $\mathrm{res}(\overline{C}^*_{\mathbb{Q}})$. Recall that $\delta_{\mathrm{res}(\alpha)}$ is the element in the $\mathbb{Q}$-span of the relative roots defined so that the pairing with $\widecheck{\mathrm{res}(\alpha)}$ is $1$ and the pairing is $0$ with all the other relative simple coroots. To show the claim proves our result, we note there is a natural isomorphism $X^*(A)_{\mathbb{Q}}  \cong X^*(A_0)_{\mathbb{Q}} \times X^*(A')_{\mathbb{Q}}$ where $A_0$ is the maximal $k$-split central torus and $A'$ is the identity component of the intersection of $A$ with the derived subgroup of $G$. Then  $_k\overline{C}^*_{\mathbb{Q}}$ corresponds under this identification to the product of $X^*(A_0)_{\mathbb{Q}}$ with the projection of $ _k\overline{C}^*_{\mathbb{Q}}$ to $X^*(A')$. Then we have a natural map $X^*(Z(G)^0)_{\mathbb{Q}} \twoheadrightarrow X^*(A_0)_{\mathbb{Q}}$ where $Z(G)^0$ is the identity component of the center of $G$ and $X^*(Z(G)^0)_{\mathbb{Q}} \subset \overline{C}^*_{\mathbb{Q}}$. Thus it suffices to show that $\mathrm{res}(\overline{C}^*_{\mathbb{Q}})$ surjects onto the projection of $ _k\overline{C}^*_{\mathbb{Q}}$ to $X^*(A')$. This latter space is identified with the set of non-negative linear combinations of the fundamental relative weights, thus proving the claim.

To prove that $\delta_{\mathrm{res}(\alpha)}$ is an element of $\mathrm{res}(\overline{C}^*_{\mathbb{Q}})$, we make use of an equivalent description of $\delta_{\mathrm{res}(\alpha)}$. It is the unique element in the $\mathbb{Q}$-span of the relative roots so that $\sigma_{\mathrm{res}(\beta)}(\delta_{\mathrm{res}(\alpha)})=\delta_{\mathrm{res}(\alpha)}$ for $\mathrm{res}(\alpha)$ and $\mathrm{res}(\beta)$ distinct simple roots and $\sigma_{\mathrm{res}(\beta)}(\delta_{\mathrm{res}(\alpha)})= \delta_{\mathrm{res}(\alpha)} - \mathrm{res}(\beta)$ when $\mathrm{res}(\alpha)=\mathrm{res}(\beta)$.

 In the case where the elements of $\Gamma \alpha$ are mutually orthogonal, we have by the above characterization of fundamental weights that the absolute fundamental weight $
 \delta_{\alpha}$ restricts to $\delta_{\mathrm{res}(\alpha)}$. In the case where $\Gamma \alpha$ has two elements that are connected in $D$, then $\delta_{\alpha}$ restricts to $2 \delta_{\mathrm{res}(\alpha)}$. In the final case, $\delta_{\alpha}$ restricts to $2 \delta_{\mathrm{res}(\alpha)}$. Thus, in all cases, we can find an element of $X^*(T)_{\mathbb{Q}}$ that restricts to $\delta_{\mathrm{res}(\alpha)}$. This completes the proof.
\end{proof}
We record an important corollary of this proposition.
\begin{corollary}{\label{posres}}
Suppose $\mu, \mu' \in X_*(T)_{\mathbb{Q}}$ and $\mu \succeq \mu'$. Let $\mu^{\Gamma}$ be the average of $\mu$ over its $\Gamma$ orbit. Then $\mu^{\Gamma} \succeq \mu'^{\Gamma}$ in $X_*(A)_{\mathbb{Q}}$. We caution that the first inequality means that $\mu - \mu'$ is a non-negative combination of absolute simple coroots, while the second means that $\mu^{\Gamma}- \mu'^{\Gamma}$ is a non-negative combination of relative simple coroots.
\end{corollary}
\begin{proof}
Recall that the action of $\Gamma$ stabilizes $\widecheck{\Delta}$. Thus for each $\gamma \in \Gamma$, we have $\gamma(\mu) \succeq \gamma(\mu')$ and so also $\mu^{\Gamma} \succeq \mu'^{\Gamma}$ in the absolute root system. Thus, we are reduced to showing that if $x \in X_*(T)^{\Gamma}_{\mathbb{Q}}$ is a non-negative combination of simple absolute coroots, then it is also a non-negative combination of simple relative coroots (under the identification $X_*(A)_{\mathbb{Q}}=X_*(T)^{\Gamma}_{\mathbb{Q}}$).

Equivalently, we need to show that if $x$ has non-negative pairing with every element of $\overline{C}^*_{\mathbb{Q}}$, then $x$ has non-negative pairing with every element of $_k\overline{C}^*_{\mathbb{Q}}$. This is indeed equivalent because $x$ has non-negative pairing with each element of $\overline{C}^*_{\mathbb{Q}}$ if and only if it has non-negative pairing with each fundamental weight $\delta_{\alpha}$ and this is the case if and only if $x$ is a non-negative combination of simple roots. 

Finally, this equivalent statement is an immediate consequence of the proposition.
\end{proof}

\printbibliography

\end{document}